\documentclass[preprint, 12pt]{elsarticle}
\usepackage{amssymb, amsmath, amsthm, exscale, relsize}
\usepackage{graphicx, subfigure, color, xcolor, cases, algorithm,algorithmic, multirow}
\usepackage[colorlinks=true, linkcolor=blue, anchorcolor=blue, citecolor=blue,]{hyperref}

\newtheorem{remark}{Remark}

\numberwithin{figure}{subsection}
\numberwithin{table}{subsection}

\allowdisplaybreaks

\begin{document}
\begin{frontmatter}
\title{When surface evolution meets Fokker-Planck equation: a novel tangential velocity model for uniform parametrization}

\author[label1]{Jiangong Pan}
\ead{mathpjg@sina.com}
\author[label2]{Guozhi Dong}
\ead{guozhi.dong@csu.edu.cn}
\author[label3]{Hailong Guo}
\ead{hailong.guo@unimelb.edu.au}
\author[label4,label5]{Zuoqiang Shi\corref{CA}}
\ead{zqshi@tsinghua.edu.cn}

\affiliation[label1]{organization={Department of Mathematical Sciences},
            addressline={Tsinghua University},
            city={Beijing},
            postcode={100084},
            country={China}}
\affiliation[label2]{organization={School of Mathematics and Statistics, HNP-LAMA},
            addressline={Central South University},
            city={Changsha},
            postcode={410083},
            country={China}}
\affiliation[label3]{organization={School of Mathematics and Statistics},
            addressline={The University of Melbourne},
            city={Parkville},
            postcode={VIC 3010},
            country={Australia}}
\affiliation[label4]{organization={Yau Mathematical Sciences Center},
            addressline={Tsinghua University},
            city={Beijing},
            postcode={100084},
            country={China}}
\affiliation[label5]{organization={YanqiLake Beijing Institute of Mathematical Sciences and Applications},
            city={Beijing},
            postcode={101408},
            country={China}}
\cortext[CA]{Corresponding author.}

\begin{abstract}
A common issue in simulating geometric evolution of surfaces is unexpected clustering of points that may cause numerical instability.
We propose a novel artificial tangential velocity method for this matter. 
The artificial tangential velocity is generated from a surface density field governed by a Fokker–Planck equation to guide the point distribution. 
A target distribution matching algorithm is developed leveraging the surface Kullback-Leibler divergence of density functions. 
The numerical method is formulated within a fully meshless framework using the moving least squares approximation, thereby eliminating the need for mesh generation and allowing flexible treatment of unstructured point cloud data. 
Extensive numerical experiments are conducted to demonstrate the robustness, accuracy, and effectiveness of the proposed approach across a variety of surface evolution problems, including the mean curvature flow.
\end{abstract}

\begin{highlights}
\item Proposes a novel artificial tangential velocity model derived from a surface Fokker–Planck equation to prevent unexpected point clustering during surface evolution.

\item Introduces flexible point redistribution algorithms that could match target densities, e.g. uniform or curvature-adapted, leveraging the surface Kullback-Leibler divergence.

\item Formulates a fully meshless numerical framework via moving least squares on point clouds, and combines BDF-k schemes in temporal discretization. Numerical results demonstrate the convergence and efficiency of the proposed algorithms.
\end{highlights}

\end{frontmatter}

\section{Introduction}\label{Sec-Intro}
We investigate the evolution of a closed surface $\Gamma(t) \subset \mathbb{R}^3$ governed by a prescribed velocity field $\boldsymbol{v} \in \mathbb{R}^3 \times [0, T]$. The surface evolution is described by the following system:
\begin{equation}\label{E-eq}
\begin{aligned}
\frac{d }{d t}\boldsymbol{X}(\boldsymbol{x}, t) &= \boldsymbol{v}(\boldsymbol{X}(\boldsymbol{x}, t), t), && \boldsymbol{x} \in \Gamma_{0}, \quad t \geq 0, \\
\boldsymbol{X}(\boldsymbol{x}, 0) &= \boldsymbol{x}, && \boldsymbol{x} \in \Gamma_{0},
\end{aligned}
\end{equation}
where $\Gamma_{0}$ denotes a closed surface embedded in $\mathbb{R}^3$, and $\Gamma(t)$ is the image of the flow map $\boldsymbol{X}(\cdot, t)$ at time $t$. 
One of the major challenges in geometric surface flow lies in the degradation of mesh quality over time, as node clustering and mesh distortion may occur, eventually leading to numerical breakdown. 
This issue affects both mesh-based and mesh-free approaches. 
To address this, B\"ansch et al. introduced a novel mesh reconstruction technique that reinitializes the mesh when excessive node aggregation or deformation is detected \cite{bansch2005finite}. 
In mesh-based methods, another effective strategy involves leveraging parametric geometric flows guided by harmonic maps of a reference surface with well-distributed mesh points. 
Related methodologies can also be found in \cite{marchandise2011high, remacle2010high}.

Among mesh-based approaches, a seminal contribution was made by Dziuk in 1990, who pioneered the numerical treatment of surface evolution under geometric flows \cite{dziuk1990algorithm}. 
He introduced the (parametric) finite element method (FEM) for approximating the mean curvature flow (MCF) of closed surfaces in three-dimensional space. 
Given an approximate surface $\Gamma_{h}(t_{j}) \subset \mathbb{R}^3$ represented via triangular faces of a polyhedron, the parametric FEM computes a parametrization $u_{h}^{j+1}: \Gamma_{h}(t_{j}) \rightarrow \mathbb{R}^3$ of the evolved surface mesh $\Gamma_{h}(t_{j+1}) = u_{h}^{j+1} (\Gamma_{h}(t_{j}))$ by solving the following weak formulation on the known surface $\Gamma_{h}(t_{j})$: Find $u_{h}^{j+1}$ in the three-dimensional vector-valued finite element space $S_{h}(\Gamma_{h}(t_{j}))^{3}$ such that
\begin{equation*}\label{Dzuik-FEM}
\begin{aligned}
\int_{\Gamma_{h}(t_{j})} \frac{u_h^{j+1} - \mathrm{id}}{\tau} \cdot \chi_h + \int_{\Gamma_{h}(t_{j})} \nabla_{\Gamma_{h}(t_{j})} u_h^{j+1} \cdot \nabla_{\Gamma_{h}(t_{j})} \chi_h = 0, \quad \forall \chi_h \in S_h\left(\Gamma_{h}(t_{j})\right)^3.
\end{aligned}
\end{equation*}
The discrete flow map $X_h^{j+1}: \Gamma_h(0) \rightarrow \Gamma_h(t_{j+1})$ is then updated via composition: $X_h^{j+1} = u_h^{j+1} \circ X_h^j$. 
Since its inception, this method has been widely adopted for simulating surface evolution under various geometric flows, such as mean curvature flow and Willmore flow \cite{bansch2005finite, bonito2010parametric, dziuk2008computational}. 
The development of numerical methods for solving partial differential equations on evolving surfaces continues to receive significant attention \cite{cheung2015localized, deckelnick2018stability, dziuk2007finite, dziuk2012fully, lehrenfeld2018stabilized, li2018direct, petras2019least, barrett2020review}

In a series work \cite{barrett2007parametric, barrett2008parametric, barrett2008parametric1}, Barrett, Garcke, and N\"urnberg introduced a novel variational formulation for the normal component of the velocity equation, permitting tangential motion of the approximated surface. 
This approach, known as the BGN method, implicitly defines the tangential velocity by enforcing that the mapping from $\Gamma_{h}(t_{j})$ to $\Gamma_{h}(t_{j+1})$ is a discrete harmonic map. 
This formulation significantly enhances mesh quality and numerical robustness without resorting to explicitly redistribute the mesh. 
The BGN scheme for MCF is stated as follows: find $u_{h}^{j+1} \in S_{h}(\Gamma_{h}(t_{j}))^{3}$ and $H_{h}^{j+1} \in S_{h}(\Gamma_{h}(t_{j}))$ satisfying the weak formulation:
\begin{equation*}\label{BGN-FEM}
\begin{aligned}
\left( \frac{u_h^{j+1} - \mathrm{id}}{\tau} \cdot \hat{\boldsymbol{n}}_h^j, \xi_h \right){\Gamma_h(t_j)}^h + \left( H_h^{j+1}, \xi_h \right){\Gamma_h(t_j)}^h &= 0, \quad \forall \xi_h \in S\left(\Gamma_h(t_j)\right), \\
\left( H_h^{j+1} \hat{\boldsymbol{n}}_h^j, \chi_h \right){\Gamma_h(t_j)}^h - \int{\Gamma_h(t_j)} \nabla_{\Gamma_h(t_j)} u_h^{j+1} \cdot \nabla_{\Gamma_h(t_j)} \chi_h &= 0, \quad \forall \chi_h \in S\left(\Gamma_h(t_j)\right)^3,
\end{aligned}
\end{equation*}
where the superscript $h$ denotes a mass-lumped inner product on the discrete surface $\Gamma_h(t_j)$, and $\hat{\boldsymbol{n}}_h^j$ is a vertex-wise averaged unit normal vector on $\Gamma_h(t_j)$. 
Since then, the BGN method has been extended and adapted to various applications. 
For example, Bao et al. \cite{bao2022volume, bao2021structure, zhao2021energy} developed parametric FEM schemes for surface diffusion incorporating artificial tangential velocities, energy stability, and volume preservation. 
These methods have been employed to simulate interface dynamics in two-phase incompressible Navier-Stokes flows \cite{barrett2013eliminating, barrett2015stable, fu2020arbitrary, ganesan2017ale}, particularly in the contexts involving contact line migration and axisymmetric geometric evolution. 

In \cite{m2017approximations}, Elliott and Fritz introduced the DeTurck flow, which reparameterizes the original geometric flow by coupling it with the harmonic heat flow in a reference domain. 
This reparameterization introduces an implicit tangential velocity component that enhances the distribution of mesh points on the evolving surface. 
In \cite{hu2022evolving}, Hu and Li demonstrated that as the time step approaches zero, the velocity produced by the BGN method asymptotically converges to a limiting velocity field $\boldsymbol{w}$ characterized by
\begin{equation}\label{HL-FEM}
\begin{aligned}
\boldsymbol{w} \cdot \boldsymbol{n} = \boldsymbol{v} \cdot \boldsymbol{n}, \quad \Delta_{\Gamma} \boldsymbol{w} = \kappa \boldsymbol{n},
\end{aligned}
\end{equation}
where $\boldsymbol{v}$ is the velocity associated with the original geometric flow and $\kappa$ is an auxiliary scalar function. 
Their analysis reveals that the tangential component defined by \eqref{HL-FEM} minimizes the instantaneous surface deformation rate, thereby offering a theoretical explanation for the consistently high mesh quality observed in practice. 
Building on this line of research, Duan and Li recently proposed an artificial tangential motion strategy in \cite{duan2024new}, which enforces that the flow map $\boldsymbol{X}(\cdot, t): \Gamma(0) \rightarrow \Gamma(t)$ minimizes the deformation energy
\begin{equation*}\label{DL-FEM}
\begin{aligned}
E[\boldsymbol{X}(\cdot, t)] = \frac{1}{2} \int_{\Gamma(0)} \left| \nabla_{\Gamma(0)} \boldsymbol{X}(\cdot, t) \right|^2,
\end{aligned}
\end{equation*}
subject to the constraint $\left( \frac{\partial \boldsymbol{X}}{\partial t} \circ \boldsymbol{X}^{-1} \right) \cdot \boldsymbol{n} = \boldsymbol{v} \cdot \boldsymbol{n}$. 
The corresponding Euler–Lagrange system reveals that this formulation enforces $\boldsymbol{X}(\cdot, t)$ to be a harmonic map with minimal deformation, thus effectively suppressing mesh distortion during surface evolution. 
Li et al. also presented a series of convergence analysis results based on the grid method \cite{li2021convergence, bai2024convergent, bai2024new, bai2024convergence}.

In this work, we are interested in the numerical simulation of geometric flows with unstructured point cloud data.
Fokker-Planck equation (FPE) has been a tool to describe the density evolution of particles \cite{risken1991fokker}. Given the velocity field $\boldsymbol{v}$ of a particle system in a Euclidean space, the trajectory of its density function $\rho$ follows the FPE:
$$\frac{\partial \rho}{\partial t}+\nabla \cdot (\rho \boldsymbol{v})=0.$$
With this inspiration, through a coupling of geometric flows and specific FPEs, we propose a novel artificial tangential velocity model for stable numerical simulations of surface evolutions. 
In the Euclidean space, giving some target density $p$, it is known that if the velocity is chosen to be 
$$\boldsymbol{v}=\nabla \log p-\nabla \log \rho,$$ then the density $\rho(t)$ will converge to the target density $p$ along the trajectory of the FPE. 
This result can be generalized to general surfaces and helps to give appropriate tangential velocity to control the density of the points. 
To the end, we derive a coupled system of surface evolution with the extra tangential velocity and a transform of some FPE to avoid undesired point clustering during surface evolution. The details can be found in Section \ref{Sec-TV-method}. 
Distinguished from traditional approaches such as the surface finite element method, which rely on structured meshes and incur significant mesh generation costs, the proposed numerical methods fit well to a meshless framework based on the moving least squares approximation. 
This avoids the effort for mesh construction and enables flexible handling of unstructured point cloud data. 
Extensive numerical experiments demonstrate the robustness and effectiveness of the proposed methods across various point cloud scenarios.

The remainder of this paper is structured as follows. In Section \ref{Sec-TV-method}, we derive the novel artificial tangential velocity by constructing a surface density diffusion equation. 
Section \ref{Sec-Discret} presents the temporal and spatial discretization of the proposed model, including a brief overview of the moving least squares method for approximating differential operators on point clouds. 
Finally, we demonstrate the effectiveness and convergence of our approach through a series of numerical experiments in Section \ref{Sec-Numerical}, highlighting its capability to characterize the evolution of surfaces given by point clouds and its robust application to MCF.

\section{Tangential velocity model inspired from Fokker-Planck equation}\label{Sec-TV-method}

\subsection{Briefs of Fokker-Planck equations in Euclidean domain}\label{sec:FPE}
Here, we give a brief introduction to FPEs \cite{risken1991fokker}. 
To simplify the notation, we restrict ourselves to the Euclidean domain in this subsection. 
For a particle system governed by the velocity field $\boldsymbol{v}$,
$$\frac{d }{d t}\boldsymbol{X}(\boldsymbol{x}, t) = \boldsymbol{v}(\boldsymbol{X}(\boldsymbol{x}, t),t),\quad \boldsymbol{x}\in \mathbb{R}^d,\quad t>0,$$
it is known that the associate density $\rho$ obeys the following FPE
\begin{equation}\label{eq:FPE}
    \frac{\partial \rho}{\partial t}+\nabla \cdot (\rho\boldsymbol{v})=0.
\end{equation}
Based on this FPE, with a given target density $p(\boldsymbol{x})$, we can design velocity field such that the density converges to the target density. One popular choice is to let
\begin{equation}\label{eq:velocity-KL}
    \boldsymbol{v}(\boldsymbol{x},t)=\nabla \log p(\boldsymbol{x})-\nabla \log \rho(\boldsymbol{x},t).
\end{equation}
This velocity field is derived from the gradient flow of Kullback-Leibler (KL) divergence
$$D_{KL}(\rho(t)||p)=\int_{\mathbb{R}^d} \rho(\boldsymbol{x},t)\log\frac{\rho(\boldsymbol{x},t)}{p(\boldsymbol{x})}d \boldsymbol{x}.$$
Calculating its time derivative and taking into account the equations \eqref{eq:FPE} and \eqref{eq:velocity-KL}, we have
\begin{align*}
    \frac{d}{dt} D_{KL}(\rho(t)||p)&=\int_{\mathbb{R}^d} \frac{\partial}{\partial t} \rho(\boldsymbol{x},t)\log\frac{\rho(\boldsymbol{x},t)}{p(\boldsymbol{x})} + \rho(\boldsymbol{x},t)\frac{\partial}{\partial t}\log\rho(\boldsymbol{x},t) d \boldsymbol{x}\\
    &=\int_{\mathbb{R}^d} -\nabla \cdot(\rho(\boldsymbol{x},t) \boldsymbol{v}(\boldsymbol{x},t))(1+\log\rho(\boldsymbol{x},t)-\log p(\boldsymbol{x})) d \boldsymbol{x}\\
    &=\int_{\mathbb{R}^d} \rho(\boldsymbol{x},t) \boldsymbol{v}(\boldsymbol{x},t)\cdot(\nabla\log\rho(\boldsymbol{x},t)-\nabla\log p(\boldsymbol{x})) d \boldsymbol{x}\\
    &=-\int_{\mathbb{R}^d} \rho(\boldsymbol{x},t) |\nabla\log\rho(\boldsymbol{x},t)-\nabla\log p(\boldsymbol{x})|^2 d \boldsymbol{x}.
\end{align*}
As a density function, $\rho$ is supposed to be nonnegative. This implies that with the velocity field in \eqref{eq:velocity-KL}, $D_{KL}(\rho||p)$ is monotonically decreasing along the trajectory of \eqref{eq:FPE}.

Actually, given \eqref{eq:velocity-KL}, the FPE \eqref{eq:FPE} becomes
\begin{equation*}
    \frac{\partial \rho}{\partial t}+\nabla \cdot (\rho \nabla \log p)=\Delta \rho.
\end{equation*}
In the next subsection, we will see that it is useful to introduce an auxiliary variable $s:=\log \rho$, and subsequently we derive the following system of equations:
\begin{align*}
    \frac{d }{d t}\boldsymbol{X}(\boldsymbol{x}, t) =& (\nabla \log p-\nabla s)(\boldsymbol{X}(\boldsymbol{x}, t), t), \\
    \frac{d}{dt}s(\boldsymbol{X}(\boldsymbol{x}, t),t)=&-\Delta \log p+\Delta s.
\end{align*}
If the target distribution is set to be uniform, i.e. $p=\text{const}$, then the above equations can be simplified to
\begin{align*}
    \frac{d }{d t}\boldsymbol{X}(\boldsymbol{x}, t) =& -\nabla s(\boldsymbol{X}(\boldsymbol{x}, t), t),\\
    \frac{d}{dt}s(\boldsymbol{X}(\boldsymbol{x}, t),t)=&\Delta s.
\end{align*}
\subsection{Tangential velocity model on surfaces}
Now we go back to our surface evolution problems with the motivation of introducing tangential velocity field based on the surface diffusion of density function. 
Suppose that we add an extra tangential velocity field $\boldsymbol{v}_{T}$ to the surface velocity on the right-hand side of \eqref{E-eq}, then the surface evolution becomes 
\begin{equation}\label{ES-vt-eq}
    \begin{aligned}
        \frac{d }{d t}\boldsymbol{X}(\boldsymbol{x}, t) &= \boldsymbol{v}(\boldsymbol{X}(\boldsymbol{x}, t), t) + \boldsymbol{v}_{T}(\boldsymbol{X}(\boldsymbol{x}, t), t), && \boldsymbol{x}\in \Gamma_{0}, \quad t\geq 0,\\
        \boldsymbol{X}(\boldsymbol{x}, 0) &= \boldsymbol{x}, && \boldsymbol{x}\in \Gamma_{0}, 
    \end{aligned}
\end{equation}
where $\boldsymbol{v}_{T}(\boldsymbol{X}(\boldsymbol{x}, t), t)\in \mathcal{T}_{\Gamma(t)}(\boldsymbol{X}(\boldsymbol{x}, t))$, $\mathcal{T}_{\Gamma(t)}(\boldsymbol{X}(\boldsymbol{x}, t))$ denotes the tangential space of $\Gamma(t)$ at $\boldsymbol{X}(\boldsymbol{x}, t)$ and $\Gamma(t)$ is the surface corresponding to $\boldsymbol{X}(\cdot, t)$. 
Note that this extra tangential velocity will not change the shape of $\Gamma(t)$ but only the distribution of points comparing to the one associated with \eqref{E-eq}.

Let $\rho(\boldsymbol{X}(\boldsymbol{x}, t), t)$ be the density function of $\Gamma(t)$, which describes the distribution of the points on the surface $\Gamma(t)$. 
With the given velocity field in \eqref{ES-vt-eq}, $\rho$ is evolved following the FPE: 
\begin{equation*}\label{C-eq}
    \begin{aligned}
        \partial^{*}_{t}\rho(\boldsymbol{X}(\boldsymbol{x}, t), t) &= -\rho(\boldsymbol{X}(\boldsymbol{x}, t), t) \text{div}_{\Gamma(t)}(\boldsymbol{v}(\boldsymbol{X}(\boldsymbol{x}, t), t) + \boldsymbol{v}_{T}(\boldsymbol{X}(\boldsymbol{x}, t), t)), && \boldsymbol{x}\in \Gamma_{0}, \quad t\geq 0, \\
        \rho(\boldsymbol{x}, 0) &= \rho_{0}(\boldsymbol{x}), && \boldsymbol{x}\in \Gamma_{0}, 
    \end{aligned}
\end{equation*}
where $\text{div}_{\Gamma(t)}$ is the divergence operator on $\Gamma(t)$, and $\partial^{*}_{t}\rho$ is the material derivative, i.e. 
\[\partial^{*}_{t}\rho(\boldsymbol{X}(\boldsymbol{x}, t), t) = \frac{d }{dt}\rho(\boldsymbol{X}(\boldsymbol{x}, t), t).\]
Here $\boldsymbol{X}(\cdot, t)$ is a diffeomorphism between $\Gamma_0$ and $\Gamma(t)$. Dividing by $\rho$ on both sides, we have
\begin{equation*}\label{C-eq-new}
    \begin{aligned}
        \frac{d }{dt}\log\rho(\boldsymbol{X}(\boldsymbol{x}, t), t) &= -\text{div}_{\Gamma(t)}(\boldsymbol{v}(\boldsymbol{X}(\boldsymbol{x}, t), t) + \boldsymbol{v}_{T}(\boldsymbol{X}(\boldsymbol{x}, t), t)), && \boldsymbol{x}\in \Gamma_{0}, \quad t\geq 0, \\
        \rho(\boldsymbol{x}, 0) &= \rho_{0}(\boldsymbol{x}), && \boldsymbol{x}\in \Gamma_{0}.
    \end{aligned}
\end{equation*}
Denote 
\begin{equation*}\label{s-rho}
    \begin{aligned}
        s(\boldsymbol{X}(\boldsymbol{x}, t), t) := \log\rho(\boldsymbol{X}(\boldsymbol{x}, t), t),
    \end{aligned}
\end{equation*}
then $s$ satisfies
\begin{equation}\label{s-C-eq-new}
    \begin{aligned}
        \frac{d }{dt}s(\boldsymbol{X}(\boldsymbol{x}, t), t) &= -\text{div}_{\Gamma(t)}(\boldsymbol{v}(\boldsymbol{X}(\boldsymbol{x}, t), t) + \boldsymbol{v}_{T}(\boldsymbol{X}(\boldsymbol{x}, t), t)), && \boldsymbol{x}\in \Gamma_{0}, \quad t\geq 0, \\
        s(\boldsymbol{x}, 0) &= \log\rho_{0}(\boldsymbol{x}), && \boldsymbol{x}\in \Gamma_{0}.
    \end{aligned}
\end{equation}

The key idea is to choose the extra tangential velocity field being
\begin{equation*}\label{vt-eq}
    \begin{aligned}
        \boldsymbol{v}_{T}(\boldsymbol{X}(\boldsymbol{x}, t), t) = -\eta\nabla_{\Gamma(t)}s(\boldsymbol{X}(\boldsymbol{x}, t), t), \quad \eta>0,
    \end{aligned}
\end{equation*}
where $\nabla_{\Gamma(t)}$ is the gradient operator on $\Gamma(t)$ and $\eta>0$ is a parameter. 
With the above tangential velocity field, the equation \eqref{s-C-eq-new} becomes
\begin{equation}\label{s-vt-C-eq-new}   
    \begin{aligned}
        \frac{d }{dt}s(\boldsymbol{X}(\boldsymbol{x}, t), t) &= -\text{div}_{\Gamma(t)}(\boldsymbol{v}(\boldsymbol{X}(\boldsymbol{x}, t), t)) + \eta\Delta_{\Gamma(t)}s(\boldsymbol{X}(\boldsymbol{x}, t), t), && \boldsymbol{x}\in \Gamma_{0}, \quad t\geq 0, \\
        s(\boldsymbol{x}, 0) &= \log\rho_{0}(\boldsymbol{x}), && \boldsymbol{x}\in \Gamma_{0},
    \end{aligned}
\end{equation}
where $\Delta_{\Gamma(t)}$ is the Laplace-Beltrami operator on $\Gamma(t)$. 
Equation \eqref{s-vt-C-eq-new} is a diffusion equation with a source term. 
When the velocity $\boldsymbol{v}=0$ (no source anymore) and the surface $\Gamma(t)$ is fixed, its solution tends to a constant on a closed smooth surface. 
This inspires us to couple the evolutionary equation \eqref{s-vt-C-eq-new} with the geometric evolution for more uniform distribution, as well as the redistribution algorithm in the next subsection.

In summary, we propose the following coupled system of equations of which the numerical discretizations will lead to more stable numerical simulation of \eqref{E-eq}:
\begin{equation}\label{ES-rho-eq}\left\{
    \begin{aligned}
        \frac{d }{d t}\boldsymbol{X}(\boldsymbol{x}, t) &= \boldsymbol{v}(\boldsymbol{X}(\boldsymbol{x}, t), t) -\eta\nabla_{\Gamma(t)}s(\boldsymbol{X}(\boldsymbol{x}, t), t), && \boldsymbol{x}\in \Gamma_{0}, \quad t\geq 0,\\
        \frac{d}{d t} s(\boldsymbol{X}(\boldsymbol{x}, t), t) &= - \text{div}_{\Gamma(t)}\boldsymbol{v}(\boldsymbol{X}(\boldsymbol{x}, t), t) + \eta\Delta_{\Gamma(t)}s(\boldsymbol{X}(\boldsymbol{x}, t), t), && \boldsymbol{x}\in \Gamma_{0}, \quad t\geq 0,\\
        \boldsymbol{X}(\boldsymbol{x}, 0) &= \boldsymbol{x}, \quad s(\boldsymbol{x}, 0) = \log\rho_{0}(\boldsymbol{x}), && \boldsymbol{x}\in \Gamma_{0}.
    \end{aligned}\right.
\end{equation}

\subsection{Points redistribution} \label{sec:redistribution} 
In the case that the distribution of the points becomes non-uniform even with the extra tangential velocity, we propose to decouple the surface evolution and the FPE, and add a point redistribution step by fixing the surface. More precisely, we set $\boldsymbol{v}$ to be zero in \eqref{ES-rho-eq} and the surface is fixed in this case. Then we solve the FPE on this fixed surface and make the points uniformly distributed over the surface. After this redistribution step, we go back to \eqref{ES-rho-eq} and resume the evolution of the surface with the velocity at the last stopping time. 

Moreover, in the redistribution step, we can choose other target distribution beside the uniform distribution. In some applications, it is better to have the distribution dependent on the geometrical structure of the surface, e.g., curvature. Our model is capable of providing this flexibility. Denote the target distribution to be $p$, then the redistribution step is to solve the following equation:
\begin{equation}\label{ES-rho-eq-nonuniform}
    \begin{aligned}
        \frac{d }{d t}\boldsymbol{X}(\boldsymbol{x}, t) &= \eta\nabla_{\Gamma(t)} \log p(\boldsymbol{X}(\boldsymbol{x}, t)) -\eta\nabla_{\Gamma(t)}s(\boldsymbol{X}(\boldsymbol{x}, t), t), && \boldsymbol{x}\in \Gamma_{t}, \quad t\geq 0,\\
        \frac{d}{d t} s(\boldsymbol{X}(\boldsymbol{x}, t), t) &= -\eta\Delta_{\Gamma(t)} \log p (\boldsymbol{X}(\boldsymbol{x}, t))+\eta\Delta_{\Gamma(t)}s(\boldsymbol{X}(\boldsymbol{x}, t), t), && \boldsymbol{x}\in \Gamma_{t}, \quad t\geq 0.
    \end{aligned}
\end{equation}
\begin{remark}
    With fixed surface $\Gamma$ and target distribution $p$, \eqref{ES-rho-eq-nonuniform} is known as the gradient flow associate with the KL divergence 
    $$D_{KL}(q||p)=\int_\Gamma q(z)\log\frac{q(z)}{p(z)}d z,$$
    for $q=\exp s$. See the analogous discussion in Euclidean domain from Section \ref{sec:FPE}. 
\end{remark}

\section{Discretization over point cloud}\label{Sec-Discret}
Since the surfaces are presented by point clouds, we employ the moving least square method to approximate the differential operator on point clouds and use backward differentiation formula (BDF) schemes in temporal direction.
\subsection{The moving least square method}\label{Sec-Derivat}
We first discuss derivatives of functions defined on hypersurfaces.
Let $\Gamma\subset\mathbb{R}^3$ be a hypersurface and suppose that it is locally parameterized by $(\alpha, \beta)\in \mathbb{R}^2$. For a given patch of $\Gamma$, it can be written as $\Gamma(\alpha, \beta) = (x(\alpha, \beta), y(\alpha, \beta), z(\alpha, \beta))$. The metric tensor $G=[g_{ij}]$ is represented by $g_{ij}=<\Gamma_{\alpha}, \Gamma_{\beta}>$, where $\Gamma_{\alpha}=(x_{\alpha}, y_{\alpha}, z_{\alpha})$ and $\Gamma_{\beta}=(x_{\beta}, y_{\beta}, z_{\beta})$ are two tangent vectors at $\boldsymbol{x}\in \Gamma$. The tangent space $T_{\boldsymbol{x}}\Gamma$ at $\boldsymbol{x}\in\Gamma$ is then spanned by $\Gamma_{\alpha}(\boldsymbol{x})$ and $\Gamma_{\beta}(\boldsymbol{x})$. Let $f:\Gamma \to \mathbb{R}$ and $f\in C^2(\Gamma)$. Under this parameterization, one can calculate the gradient of $f$  using the formula below \cite{do1976differential}
\begin{equation*}\label{nabla-f} 
    \begin{aligned}
        \nabla_{\Gamma}f = \left[\Gamma_{\alpha}, \Gamma_{\beta}\right]G^{-1} \nabla f = \left(g^{11}\frac{\partial f}{\partial\alpha} + g^{12}\frac{\partial f}{\partial\beta}\right)\Gamma_{\alpha} + \left(g^{21}\frac{\partial f}{\partial\alpha} + g^{22}\frac{\partial f}{\partial\beta}\right)\Gamma_{\beta},
    \end{aligned}
\end{equation*}
where $g^{ij}$ are the components of $G^{-1}$, the inverse of the metric tensor $G$. Denote $g={\rm det}(G)$, the Laplace-Beltrami operation to $f$ is
\begin{equation}\label{triangle-f}
    \begin{aligned}
        \Delta_{\Gamma}f = \frac{1}{\sqrt{g}}\left(
        \frac{\partial}{\partial \alpha}\left(\sqrt{g}g^{11}\frac{\partial f}{\partial \alpha}\right)\right. &+ 
        \frac{\partial}{\partial \alpha}\left(\sqrt{g}g^{12}\frac{\partial f}{\partial \beta}\right) \\
        & \left.+ \frac{\partial}{\partial \beta}\left(\sqrt{g}g^{21}\frac{\partial f}{\partial \alpha}\right) + 
        \frac{\partial}{\partial \beta}\left(\sqrt{g}g^{22}\frac{\partial f}{\partial \beta}\right)\right).
    \end{aligned}
\end{equation}

In this article, only point clouds sampled from hypersurfaces are available, but not the surfaces themselves.
Thus, in the rest of this section, we briefly review the moving least squares method (MLS) \cite{liang2013solving} to approximate the surfaces in a local coordinate system and then compute the metric tensors as well as the derivatives of functions at each point. 
First of all, we use principal component analysis to construct the local coordinate system. 
Given a point cloud $\boldsymbol{P}=\{\boldsymbol{p}_{i}|i=1,2,\ldots,N_{\boldsymbol{x}}\}$ with $N_{\boldsymbol{x}}$ points sampled from a two-dimensional manifold in $\mathbb{R}^{3}$. 
Define $\Lambda(\boldsymbol{p}_{i})$ to be the set of adjacent points to $\boldsymbol{x}_{i}$ obtained by the K-nearest-neighbor (KNN) method. Define the covariance matrix $P_{i}$ at $\boldsymbol{x}_{i}$, $P_{i} = \sum_{k\in\Lambda(\boldsymbol{p}_{i})}(\boldsymbol{p}_{k}-\boldsymbol{c}_{i})^T(\boldsymbol{p}_{k}-\boldsymbol{c}_{i})$, where $\boldsymbol{c}_{i}$ is the local barycenter $\boldsymbol{c}_{i}=\frac{1}{K}\sum_{k\in\Lambda(\boldsymbol{p}_{i})}\boldsymbol{p}_{k}$. 
Then through $P_{i}$, we can obtain the sorted eigenvalues $\lambda_{i,1}>\lambda_{i,2}>\lambda_{i,3}$ and the corresponding eigenvectors $(\boldsymbol{e}_{i,1}, \boldsymbol{e}_{i,2}, \boldsymbol{e}_{i,3})$. 

Since the surface we considering is two-dimensional, in the local coordinate system $\{\boldsymbol{x}_{i};\boldsymbol{e}_{i,1}, \boldsymbol{e}_{i,2}, \boldsymbol{e}_{i,3}\}$ of the point $\boldsymbol{x}_{i}$ and local coordinates $(\alpha_{i}, \beta_{i}, \gamma_{i})$, the local second-order binary polynomial $\gamma_{i}(\alpha, \beta)$ is obtained by minimizing the following weighted sum:
\begin{equation*}\label{mls-surface}
    \begin{aligned}
\sum\limits_{k\in\Lambda(\boldsymbol{x}_{i})}\omega(\|\boldsymbol{x}_{k}-\boldsymbol{x}_{i}\|)(\gamma_{i}(\alpha_{i,k}, \beta_{i,k}) - \gamma_{i,k})^2,
    \end{aligned}
\end{equation*}
where $\Lambda(\boldsymbol{x}_{i})$ is the set of adjacent points to $\boldsymbol{x}_{i}$, and $\omega(d_{k})$ is the weight coefficient with $d_{k} = \|\boldsymbol{x}_{k}-\boldsymbol{x}\|$ being the Euclidean distances between $\boldsymbol{x}$ and the position of data point $\boldsymbol{x}_{k}$. 
Thus, $\Gamma_{i}=(\alpha, \beta, \gamma_{i}(\alpha, \beta))$ is a smooth approximation of some underlying surface  near the point $\boldsymbol{x}_i$ under the local coordinate system $\{\boldsymbol{x}_{i}; \boldsymbol{e}_{i,1}, \boldsymbol{e}_{i,2}, \boldsymbol{e}_{i,3}\}$. 

\begin{remark}\label{remark1}
Although the weight function in MLS does not affect the calculation accuracy, the stability will be challenged. There are many alternative functions in \cite{do1976differential}. We use the most popular Wendland function
\begin{equation*}\label{Wendland}
    \begin{aligned}
        \omega(d)=\left(1-\frac{d}{D}\right)^{4}\left(\frac{4d}{D}+1\right),
    \end{aligned}
\end{equation*}
which is defined on the interval $d\in[0, D]$ and $\omega(0)=1$, $\omega(D)=0$, $\omega^{\prime}(0)=0$ and $\omega^{\prime\prime}(0)=0$.
\end{remark}

Particularly, we use the local binary second-order polynomials $\gamma_{i}(\alpha, \beta) = c_{0} + c_{1}\alpha + c_{2}\beta + c_{3}\alpha^{2} + c_{4}\alpha\beta + c_{5}\beta^{2}$. The local basis consists of two tangent vectors given by $\Gamma_{\alpha}(\boldsymbol{x}_{i})= (1, 0, \frac{\partial\gamma_{i}}{\partial\alpha}) = (1, 0, c_{1})$ and $\Gamma_{\beta}(\boldsymbol{x}_{i})= (0, 1, \frac{\partial\gamma_{i}}{\partial\beta}) = (0, 1, c_{2})$. Based on this local polynomial surface, we compute the metric tensor $G(\boldsymbol{x})$ and its inverse $G^{-1}(\boldsymbol{x})$ which are functions dependent on $\gamma_{i}$. Then, the gradient in local coordinate system can be computed as follows
\begin{equation}\label{mls-nabla}
    \begin{aligned}
        \nabla_{\Gamma}f(\boldsymbol{x}_{i}) = \left(g^{11}(\boldsymbol{x}_{i})\frac{\partial f}{\partial\alpha}(\boldsymbol{x}_{i})\right. & \left.+ g^{12}(\boldsymbol{x}_{i})\frac{\partial f}{\partial\beta}(\boldsymbol{x}_{i})\right)\Gamma_{\alpha}(\boldsymbol{x}_{i}) \\
        &+ \left(g^{21}(\boldsymbol{x}_{i})\frac{\partial f}{\partial\alpha}(\boldsymbol{x}_{i}) + g^{22}(\boldsymbol{x}_{i})\frac{\partial f}{\partial\beta}(\boldsymbol{x}_{i})\right)\Gamma_{\beta}(\boldsymbol{x}_{i}),
    \end{aligned}
\end{equation}
and Laplace-Beltrami operation is realized via
\begin{equation}\label{mls-laplace}
    \begin{aligned}
        \Delta_{\Gamma}f(\boldsymbol{x}_{i}) = A_{0}(\boldsymbol{x}_{i})\frac{\partial f}{\partial\alpha}(\boldsymbol{x}_{i}) + A_{1}(\boldsymbol{x}_{i})\frac{\partial f}{\partial\beta}(\boldsymbol{x}_{i}) &+ A_{2}(\boldsymbol{x}_{i})\frac{\partial^{2} f}{\partial\alpha^{2}}(\boldsymbol{x}_{i}) \\
        & + A_{3}(\boldsymbol{x}_{i})\frac{\partial^{2} f}{\partial\alpha\partial\beta}(\boldsymbol{x}_{i}) + A_{4}(\boldsymbol{x}_{i})\frac{\partial^{2} f}{\partial\beta^{2}}(\boldsymbol{x}_{i}),
    \end{aligned}
\end{equation}
where 
\begin{equation}\label{mls-detail}
    \begin{aligned}
        &g^{11}(\boldsymbol{x}_{i}) = \frac{1+c_{2}^2}{1+c_{1}^2+c_{2}^2}, \quad  g^{12}(\boldsymbol{x}_{i}) = \frac{-c_{1}c_{2}}{1+c_{1}^2+c_{2}^2}, \\
        &g^{21}(\boldsymbol{x}_{i}) = \frac{-c_{1}c_{2}}{1+c_{1}^2+c_{2}^2}, \quad  g^{22}(\boldsymbol{x}_{i}) = \frac{1+c_{1}^2}{1+c_{1}^2+c_{2}^2}. \\
        &A_0(\boldsymbol{x}_{i}) = \frac{1}{\sqrt{g(\boldsymbol{x}_{i})}}\left((\sqrt{g})_{\alpha}(\boldsymbol{x}_{i})g^{11}(\boldsymbol{x}_{i}) + (\sqrt{g})_{\beta}(\boldsymbol{x}_{i})g^{21}(\boldsymbol{x}_{i})\right) + g^{11}_{\alpha}(\boldsymbol{x}_{i}) + g^{21}_{\beta}(\boldsymbol{x}_{i}), \\
        &A_1(\boldsymbol{x}_{i}) = \frac{1}{\sqrt{g(\boldsymbol{x}_{i})}}\left((\sqrt{g})_{\alpha}(\boldsymbol{x}_{i})g^{12}(\boldsymbol{x}_{i}) + (\sqrt{g})_{\beta}(\boldsymbol{x}_{i})g^{22}(\boldsymbol{x}_{i})\right) + g^{12}_{\alpha}(\boldsymbol{x}_{i}) + g^{22}_{\beta}(\boldsymbol{x}_{i}), \\ 
        &A_2(\boldsymbol{x}_{i}) = g^{11}(\boldsymbol{x}_{i}), \quad A_3(\boldsymbol{x}_{i}) = \left(g^{12}(\boldsymbol{x}_{i})+g^{21}(\boldsymbol{x}_{i})\right), \quad A_4(\boldsymbol{x}_{i}) = g^{22}(\boldsymbol{x}_{i}).
    \end{aligned}
\end{equation}
The formulas in \eqref{mls-laplace} are derived from Equation \eqref{triangle-f}. Note that these approximated quantities depend only on local polynomial coefficients $c_{1}, c_{2},..., c_{5}$. 

We are now ready to reconstruct the function and its derivatives under these local coordinate systems with the quantities estimated in \eqref{mls-detail}. Similarly, we reconstruct the function $f_{s}$ locally from a set of function values $\{f_{i,k}\}$ using second-order binary polynomial $f_{i}(\alpha, \beta)$ in the local coordinate system $\{\boldsymbol{x}_{i};\boldsymbol{e}_{i,1}, \boldsymbol{e}_{i,2}, \boldsymbol{e}_{i,3}\}$ via minimizing the following weighted sum
\begin{equation*}\label{mls-function}
    \begin{aligned} \sum\limits_{k\in\Lambda(\boldsymbol{x}_{i})}\omega(\|\boldsymbol{p}_{k}-\boldsymbol{x}_{i}\|)(f_{i}(\alpha_{i,k}, \beta_{i,k}) - f_{i,k})^2.
    \end{aligned}
\end{equation*}
Then suppose $f_{i}(\alpha, \beta) = \bar{c}_{0} + \bar{c}_{1}\alpha + \bar{c}_{2}\beta + \bar{c}_{3}\alpha^{2} + \bar{c}_{4}\alpha\beta + \bar{c}_{5}\beta^{2}$.
We update the formulas in \eqref{mls-nabla} and \eqref{mls-laplace} taking into account that the derivatives at $\boldsymbol{x}_{i}$ is at the origin in the local coordinate system, i.e., $\alpha=0$ and $\beta=0$. Then we end up with
\begin{equation}\label{mls-new}
    \begin{aligned}
        \nabla_{\Gamma}f(\boldsymbol{x}_{i}) &= \left(g^{11}(\boldsymbol{x}_{i})\bar{c}_{1} + g^{12}(\boldsymbol{x}_{i})\bar{c}_{2}\right)(\boldsymbol{e}_{i,1} + c_{1}\boldsymbol{e}_{i,3}) + \left(g^{21}(\boldsymbol{x}_{i})\bar{c}_{1} + g^{22}(\boldsymbol{x}_{i})\bar{c}_{2}\right)(\boldsymbol{e}_{i,2} + c_{2}\boldsymbol{e}_{i,3}), \\
        \Delta_{\Gamma}f(\boldsymbol{x}_{i}) &= A_{0}(\boldsymbol{x}_{i})\bar{c}_{1} + A_{1}(\boldsymbol{x}_{i})\bar{c}_{2} + 2A_{2}(\boldsymbol{x}_{i})\bar{c}_{3} + A_{3}(\boldsymbol{x}_{i})\bar{c}_{4} + 2A_{4}(\boldsymbol{x}_{i})\bar{c}_{5},
    \end{aligned}
\end{equation}
where $g^{ij}, A_i$ are quantities provided in \eqref{mls-detail}.

\subsection{Time discretization}
In this part, we consider time semi-discretization of the system of equations \eqref{ES-rho-eq}. 
Let time steps $t_i=(i-1)\Delta t, \, i=1, 2, \cdots, N_t$, $\tau = \frac{T}{N_t-1}$. 
The time derivative $\frac{d \boldsymbol{X}(\cdot, t)}{d t}$ is approximated using BDF:
\begin{equation}\label{Dis-x-eq}
    \begin{aligned}
        \frac{d }{d t}\boldsymbol{X}(\cdot, t_{k}) \sim \left\{
        \begin{array}{l}
            \frac{\boldsymbol{X}(\cdot, t_{k}) - \boldsymbol{X}(\cdot, t_{k-1})}{\tau} + \mathcal{O}(\tau^{1}),\\
            \frac{\frac{3}{2}\boldsymbol{X}(\cdot, t_{k}) - 2\boldsymbol{X}(\cdot, t_{k-1}) + \frac{1}{2}\boldsymbol{X}(\cdot, t_{k-2})}{\tau} + \mathcal{O}(\tau^{2}), \\ 
            \frac{\frac{11}{6}\boldsymbol{X}(\cdot, t_{k}) - 3\boldsymbol{X}(\cdot, t_{k-1}) + \frac{3}{2}\boldsymbol{X}(\cdot, t_{k-2}) - \frac{1}{3}\boldsymbol{X}(\cdot, t_{k-3})}{\tau} + \mathcal{O}(\tau^{3}),  
        \end{array}\right.
    \end{aligned}
\end{equation}
and the material derivative of $s$ is approximately computed as
\begin{equation}\label{Dis-s-eq}
    \begin{aligned}
        \frac{d }{d t}s(\boldsymbol{X}_{k}, t_{k}) \sim \left\{
        \begin{array}{l}
        \frac{s(\boldsymbol{X}_{k}, t_{k}) - s(\boldsymbol{X}_{k-1}, t_{k-1})}{\tau} + \mathcal{O}(\tau^{1}),\\
        \frac{\frac{3}{2}s(\boldsymbol{X}_{k}, t_{k}) - 2s(\boldsymbol{X}_{k-1}, t_{k-1}) + \frac{1}{2}s(\boldsymbol{X}_{k-2}, t_{k-2})}{\tau} + \mathcal{O}(\tau^{2}),\\
        \frac{\frac{11}{6}s(\boldsymbol{X}_{k}, t_{k}) - 3s(\boldsymbol{X}_{k-1}, t_{k-1}) + \frac{3}{2}s(\boldsymbol{X}_{k-2}, t_{k-2}) - \frac{1}{3}s(\boldsymbol{X}_{k-3}, t_{k-3})}{\tau} + \mathcal{O}(\tau^{3}).
        \end{array}\right.
    \end{aligned}
\end{equation}
To simplify the notation, we use $\boldsymbol{X}_{k} = \boldsymbol{X}(\cdot, t_{k})$ and $s_{k} = s(\boldsymbol{X}_{k}, t_{k})$.

Combining \eqref{Dis-s-eq} and \eqref{Dis-x-eq}, we get the time semi-discrete format of \eqref{ES-rho-eq}: 
\begin{equation}\label{SemiDis-ES-rho-eq}
    \begin{aligned}
        \frac{a\boldsymbol{X}_{k} - \hat{\boldsymbol{X}}_{k-1}}{\tau} &= \boldsymbol{v}(\bar{\boldsymbol{X}}_{k}, t_{k}) -\eta\nabla_{\bar{\Gamma}_{k}}s_{k}, \\
        \frac{as_{k} - \hat{s}_{k-1}}{\tau} &= - \text{div}_{\bar{\Gamma}_{k}}\boldsymbol{v}(\bar{\boldsymbol{X}}_{k}, t_{k}) + \eta\Delta_{\bar{\Gamma}_{k}}s_{k}, \\
        \boldsymbol{X}(\boldsymbol{x}, 0) &= \boldsymbol{x}, \quad s(\boldsymbol{x}, 0) = \log\rho_{0}(\boldsymbol{x}),
    \end{aligned}
\end{equation}
where $a$, $\hat{\boldsymbol{X}}_{k-1}$ and $\hat{s}_{k-1}$ are defined as
\begin{equation}\label{BDF-eq}
    \begin{aligned}
        & \text{BDF1:}\ a = 1, && \hat{\boldsymbol{X}}_{k-1} = \boldsymbol{X}_{k-1}, && \hat{s}_{k-1} = s_{k-1}, \\
        & \text{BDF2:}\ a = \frac{3}{2}, && \hat{\boldsymbol{X}}_{k-1} = 2\boldsymbol{X}_{k-1} - \frac{1}{2}\boldsymbol{X}_{k-2}, && \hat{s}_{k-1} = 2s_{k-1} - \frac{1}{2}s_{k-2}, \\
        & \text{BDF3:}\ a = \frac{11}{6}, && \hat{\boldsymbol{X}}_{k-1} = 3\boldsymbol{X}_{k-1} - \frac{3}{2}\boldsymbol{X}_{k-2} + \frac{1}{3}\boldsymbol{X}_{k-3}, && \hat{s}_{k-1} = 3s_{k-1} - \frac{3}{2}s_{k-2} + \frac{1}{3}s_{k-3}.
    \end{aligned}
\end{equation}
$\bar{\boldsymbol{X}}_{k}$ in \eqref{SemiDis-ES-rho-eq} is the BDF-k extrapolation, e.g.,
\begin{equation}\label{BDF-explo-eq}
    \begin{aligned}
        & \text{BDF1:} && \bar{\boldsymbol{X}}_{k} = \boldsymbol{X}_{k-1}, \\
        & \text{BDF2:} && \bar{\boldsymbol{X}}_{k} = 2\boldsymbol{X}_{k-1} - \boldsymbol{X}_{k-2}, \\
        & \text{BDF3:} && \bar{\boldsymbol{X}}_{k} = 3\boldsymbol{X}_{k-1} - 3\boldsymbol{X}_{k-2} + \boldsymbol{X}_{k-3}. 
    \end{aligned}
\end{equation}
Moreover, when applying higher-order BDF methods (e.g., BDF2 and BDF3), suitable initial values should be provided to reach the desired temporal accuracy. 
To mitigate the impact of low-accuracy initial approximations, we adopt the fourth-order Runge–Kutta (RK4) method \cite{rosser1967runge}. 
The implementation details with respect to BDF2 and BDF3 schemes are provided in Algorithm \ref{BDF2-Algorithm} and Algorithm \ref{BDF3-Algorithm}, respectively. 
The BDF1 scheme, which is rather straightforward, is omitted for brevity. 
\begin{algorithm}[!ht]
    \renewcommand{\algorithmicrequire}{\textbf{Input:}}
	\renewcommand{\algorithmicensure}{\textbf{Output:}}
	\caption{BDF2 Algorithm}
    \label{BDF2-Algorithm}
    \begin{algorithmic}[1]
        \REQUIRE Sample points $\boldsymbol{X}_{0}=\boldsymbol{x}$ on initial surface $\Gamma(0)$, Initial value $s_{0}$, Time step $\tau$, End time $T$; 
	    \ENSURE Discrete solution $\boldsymbol{X}_{T/\tau}$, $s_{T/\tau}$ at end time;
        
        \STATE Taking $\boldsymbol{X}_{0}$, $s_{0}$, $\tau$ and $\bar{\Gamma}_{1}=\{\boldsymbol{X}_{0}^{i}|i=1, \ldots, N_{\boldsymbol{x}}\}$ as inputs, calculate $\boldsymbol{X}_{1}$, $s_{1}$ using the BDF1 format of \eqref{SemiDis-ES-rho-eq}; 
        \STATE Let $k = 1$;

        \WHILE {$k < T/\tau$}
            \STATE Taking $\boldsymbol{X}_{k-1}$ and $\boldsymbol{X}_{k-2}$ as inputs, calculate $\bar{\boldsymbol{X}}_{k}$ using the BDF2 format of \eqref{BDF-explo-eq};
            \STATE Taking $\boldsymbol{X}_{k-1}$, $\boldsymbol{X}_{k}$, $s_{k-1}$, $s_{k}$, $\tau$ and $\bar{\Gamma}_{k+1}=\{\boldsymbol{X}_{k}^{i}|i=1, \ldots, N_{\boldsymbol{x}}\}$ as inputs, calculate $\boldsymbol{X}_{k+1}$, $s_{k+1}$ using the BDF2 format of \eqref{SemiDis-ES-rho-eq};
            \STATE Let $k = k+1$;
        \ENDWHILE
        
        \STATE \textbf{return} $\boldsymbol{X}_{T/\tau}, s_{T/\tau}$.
    \end{algorithmic}
\end{algorithm}

\begin{algorithm}[!ht]
    \renewcommand{\algorithmicrequire}{\textbf{Input:}}
	\renewcommand{\algorithmicensure}{\textbf{Output:}}
	\caption{BDF3 Algorithm}
    \label{BDF3-Algorithm}
    \begin{algorithmic}[1]
        \REQUIRE Sample points $\boldsymbol{X}_{0}=\boldsymbol{x}$ on initial surface $\Gamma(0)$, Initial value $s_{0}$, Time step $\tau$, End time $T$; 
	    \ENSURE Discrete solution $\boldsymbol{X}_{T/\tau}$, $s_{T/\tau}$ at end time;
        
        \STATE Taking $\boldsymbol{X}_{0}$, $s_{0}$, $\tau$ and $\bar{\Gamma}_{1}=\{\boldsymbol{X}_{0}^{i}|i=1, \ldots, N_{\boldsymbol{x}}\}$ as inputs, calculate $\boldsymbol{X}_{1}$, $s_{1}$ using the RK4 method; 
        \STATE Taking $\boldsymbol{X}_{1}$, $s_{1}$, $\tau$ and $\bar{\Gamma}_{2}=\{\boldsymbol{X}_{1}^{i}|i=1, \ldots, N_{\boldsymbol{x}}\}$ as inputs, calculate $\boldsymbol{X}_{2}$, $s_{2}$ using the BDF2 format of \eqref{SemiDis-ES-rho-eq}; 
        \STATE Let $k = 2$;

        \WHILE {$k < T/\tau$}
            \STATE Taking $\boldsymbol{X}_{k-1}$, $\boldsymbol{X}_{k-2}$ and $\boldsymbol{X}_{k-3}$ as inputs, calculate $\bar{\boldsymbol{X}}_{k}$ using the BDF3 format of \eqref{BDF-explo-eq};
            \STATE Taking $\boldsymbol{X}_{k-2}$, $\boldsymbol{X}_{k-1}$, $\boldsymbol{X}_{k}$, $s_{k-2}$, $s_{k-1}$, $s_{k}$, $\tau$ and $\bar{\Gamma}_{k+1}=\{\boldsymbol{X}_{k}^{i}|i=1, \ldots, N_{\boldsymbol{x}}\}$ as inputs, calculate $\boldsymbol{X}_{k+1}$, $s_{k+1}$ using the BDF3 format of \eqref{SemiDis-ES-rho-eq};
            \STATE Let $k = k+1$;
        \ENDWHILE
        
        \STATE \textbf{return} $\boldsymbol{X}_{T/\tau}, s_{T/\tau}$.
    \end{algorithmic}
\end{algorithm}

\subsection{Fully discrete schemes on point clouds}
For the spatial discretization, we consider a point cloud $\boldsymbol{P}=\{\boldsymbol{x}_{i}|i=1,2,\ldots,N_{\boldsymbol{x}}\}$ with $N_{\boldsymbol{x}}$ points sampled from a two-dimensional manifold in $\mathbb{R}^{3}$. 
The scalar function $f_{s}:\Gamma \to \mathbb{R}$ is approximated  by $F=[f_{1}, f_{2}, \ldots, f_{N_{\boldsymbol{x}}}]^T$ in the sense of $f_{i}\sim f_{s}(\boldsymbol{x}_i)$ and \eqref{mls-new}  becomes
\begin{equation}\label{mls-new-matrix}
    \begin{aligned}
        \nabla_{\Gamma}f_{s} \sim [M_{G1}\ M_{G2}\ M_{G3}]F = M_{G}F, \quad \Delta_{\Gamma}f_{s} \sim M_{L}F.
    \end{aligned}
\end{equation}
In \eqref{mls-new-matrix}, $M_{G1}, M_{G2}, M_{G3}$ and $M_{L}$ are $N_{\boldsymbol{x}}\times N_{\boldsymbol{x}}$ matrices which are typically sparse. 
Then, we use \eqref{mls-new-matrix} to obtain the full-discrete form of \eqref{ES-rho-eq}: 
\begin{equation}\label{fullDis-ES-rho-eq}
    \begin{aligned}
        \frac{a}{\tau}\boldsymbol{X}_{k} + \eta M_{G}^{T}(\bar{\boldsymbol{X}}_{k})s_{k} &= \frac{1}{\tau}\hat{\boldsymbol{X}}_{k-1} + \boldsymbol{v}(\bar{\boldsymbol{X}}_{k}, t_{k}), \\
        \frac{a}{\tau}s_{k} - \eta M_{L}(\bar{\boldsymbol{X}}_{k})s_{k} &= \frac{1}{\tau}\hat{s}_{k-1} - M_{G}(\bar{\boldsymbol{X}}_{k})\boldsymbol{v}(\bar{\boldsymbol{X}}_{k}, t_{k}), \\
        \boldsymbol{X}(\boldsymbol{x}, 0) &= \boldsymbol{x}, \quad s(\boldsymbol{x}, 0) = \log\rho_{0}(\boldsymbol{x}),
    \end{aligned}
\end{equation}
where $a$, $\hat{\boldsymbol{X}}_{k-1}$, $\hat{S}_{k-1}$ and $\bar{\boldsymbol{X}}_{k}$ are defined in \eqref{BDF-eq}. 
In the case of point clouds, in order to estimate the density $\hat{\rho}_{0}(\boldsymbol{x}_{i})$ of the underlying surfaces, we rely on some reconstructed triangle template $\Lambda_{\triangle}(\boldsymbol{x}_{i})$ out of the point clouds. 
The construction procedure for these triangles is detailed in \cite{lai2013local}, which involves forming Delaunay triangles from points $\Lambda(\boldsymbol{x}_{i})$ in the vicinity of $\boldsymbol{x}_{i}$ and subsequently excluding any triangles that do not include $\boldsymbol{x}_{i}$ as a vertex. 
Then, the initial density $\rho_{0}(\boldsymbol{x})$ is estimated via the following formula:
\begin{equation}\label{rho-0}
    \begin{aligned}
        \rho_{0}(\boldsymbol{x}_{i}) \propto \hat{\rho}_{0}(\boldsymbol{x}_{i}) = \frac{1}{S(\boldsymbol{x}_i)} = \frac{3}{\sum\limits_{K_{j}\in\Lambda_{\triangle}(\boldsymbol{x}_{i})}(Area(K_{j}))}.
    \end{aligned}
\end{equation}

Since the velocity field $\boldsymbol{v}(\boldsymbol{X}(\boldsymbol{x},t), t)$ is given arbitrarily, and we use a constant $\eta$ in our model for simplicity, the tangential velocity field $\boldsymbol{v}_{T}(\boldsymbol{X}(\boldsymbol{x},t), t)$ might still be insufficient to counteract the component of the given velocity $\boldsymbol{v}(\boldsymbol{X}(\boldsymbol{x},t), t)$ in the tangential direction. 
The aggregation effect of points may still exist in such cases. 
This phenomenon has been observed in some of our subsequent experiments, cf. Figure~\ref{imp-dumbbell} for details. 
In this case, we can add point redistribution step as mentioned in Section \ref{sec:redistribution}. 
The core idea is to halt the normal motion by setting $\boldsymbol{v} = 0$ in \eqref{fullDis-ES-rho-eq}, and solve the following equation:
\begin{equation}\label{Correction-eq}
    \begin{aligned}
        \frac{a}{\tau}\boldsymbol{X}_{k} + \eta M_{G}^{T}(\bar{\boldsymbol{X}}_{k})s_{k} &= \frac{1}{\tau}\hat{\boldsymbol{X}}_{k-1}, \\
        \frac{a}{\tau}s_{k} - \eta M_{L}(\bar{\boldsymbol{X}}_{k})s_{k} &= \frac{1}{\tau}\hat{s}_{k-1}, \\
        \boldsymbol{X}(\boldsymbol{x}, T) = \boldsymbol{X}_{T}, \quad s(\boldsymbol{x}, T) &= \log \rho_{T}(\boldsymbol{x}).
    \end{aligned}
\end{equation}
$\rho_T$ obtained via \eqref{rho-0} at time $T$. 
The computational procedure is outlined in Algorithm \ref{RA-Algorithm}. 
Naturally, the redistribution can be invoked once the variation in $s$ exceeds a predefined threshold, thereby maintaining accuracy throughout the surface evolution process. 

Once a target distribution $p(\boldsymbol{x})$ is desired to be matched, the equation \eqref{ES-rho-eq-nonuniform} is discretized as follow:
\begin{equation}\label{Correction-eq-p}
    \begin{aligned}
        \frac{a}{\tau}\boldsymbol{X}_{k} + \eta M_{G}^{T}(\bar{\boldsymbol{X}}_{k})s_{k} &= \frac{1}{\tau}\hat{\boldsymbol{X}}_{k-1} + \eta M_{G}^{T}(\bar{\boldsymbol{X}}_{k})\log{p(\bar{\boldsymbol{X}}_{k})}, \\
        \frac{a}{\tau}s_{k} - \eta M_{L}(\bar{\boldsymbol{X}}_{k})s_{k} &= \frac{1}{\tau}\hat{s}_{k-1} - \eta M_{L}(\bar{\boldsymbol{X}}_{k})\log{p(\bar{\boldsymbol{X}}_{k})}, \\
        \boldsymbol{X}(\boldsymbol{x}, T) = \boldsymbol{X}_{T}, \quad s(\boldsymbol{x}, T) &= \log\rho_{T}(\boldsymbol{x}).
    \end{aligned}
\end{equation}
The computational procedure is outlined in Algorithm \ref{Target-RA-Algorithm}. 

\begin{algorithm}[!ht]
    \renewcommand{\algorithmicrequire}{\textbf{Input:}}
	\renewcommand{\algorithmicensure}{\textbf{Output:}}
	\caption{Re-distribution Algorithm}
    \label{RA-Algorithm}
    \begin{algorithmic}[1]
        \REQUIRE Evolution points $\boldsymbol{X}_{T}$, Initial value $s_{T}$, Estimated value $\hat{\rho}_{T}$, Time step $\tau$, Small time step $\tilde{\tau}\sim\tau^{2}$, Threshold $\varepsilon_{0}=5e-5$; 
	    \ENSURE Discrete solution $\boldsymbol{X}_{c}$, $s_{c}$ at end time;

        \STATE Correction $\hat{s}_{T}=\log\hat{\rho}_{T}$;
        \STATE Taking $\boldsymbol{X}_{T}$, $\hat{s}_{T}$, $\tau$ and $\bar{\Gamma}_{T+1}=\{\boldsymbol{X}_{T}^{i}|i=1, \ldots, N_{\boldsymbol{x}}\}$ as inputs, calculate $\boldsymbol{X}_{T+1}$, $s_{T+1}$ using the BDF1 format of \eqref{Correction-eq}; 
        \STATE Let $k = 1$ and $\varepsilon_{s}=2$;

        \WHILE {$\varepsilon_{s} > \varepsilon_{0}$}
            \STATE Taking $\boldsymbol{X}_{T+k-1}$ and $\boldsymbol{X}_{T+k-2}$ as inputs, calculate $\bar{\boldsymbol{X}}_{T+k}$ using the BDF2 format of \eqref{BDF-explo-eq};
            \STATE Taking $\boldsymbol{X}_{T+k-1}$, $\boldsymbol{X}_{T+k}$, $s_{T+k-1}$, $s_{T+k}$, $\tau$ and $\bar{\Gamma}_{T+k+1}=\{\boldsymbol{X}_{T+k}^{i}|i=1, \ldots, N_{\boldsymbol{x}}\}$ as inputs, calculate $\boldsymbol{X}_{T+k+1}$, $s_{T+k+1}$ using the BDF2 format of \eqref{Correction-eq};
            \STATE Calculate $\varepsilon_{s} = |\max(s_{T+k+1}) - \min(s_{T+k+1})|$;
            \STATE Let $k = k+1$;
        \ENDWHILE
        
        \STATE \textbf{return} $\boldsymbol{X}_{c}, s_{c}$.
    \end{algorithmic}
\end{algorithm}

\begin{algorithm}[!ht]
    \renewcommand{\algorithmicrequire}{\textbf{Input:}}
	\renewcommand{\algorithmicensure}{\textbf{Output:}}
	\caption{Target distribution matching Algorithm}
    \label{Target-RA-Algorithm}
    \begin{algorithmic}[1]
        \REQUIRE Evolution points $\boldsymbol{X}_{T}$, Initial value $s_{T}$, Target distribution $p$, Estimated value $\hat{\rho}_{T}$, Time step $\tau$, Small time step $\tilde{\tau}\sim\tau^{2}$, Threshold $\varepsilon_{0}=5e-5$; 
	    \ENSURE Discrete solution $\boldsymbol{X}_{c}$, $s_{c}$ at end time;

        \STATE Correction $\hat{s}_{T}=\log\hat{\rho}_{T}$;
        \STATE Taking $\boldsymbol{X}_{T}$, $\hat{s}_{T}$, $\tau$ and $\bar{\Gamma}_{T+1}=\{\boldsymbol{X}_{T}^{i}|i=1, \ldots, N_{\boldsymbol{x}}\}$ as inputs, calculate $\boldsymbol{X}_{T+1}$, $s_{T+1}$ using the BDF1 format of \eqref{Correction-eq-p}; 
        \STATE Let $k = 1$ and $\varepsilon_{s}=2$;

        \WHILE {$\varepsilon_{s} > \varepsilon_{0}$}
            \STATE Taking $\boldsymbol{X}_{T+k-1}$ and $\boldsymbol{X}_{T+k-2}$ as inputs, calculate $\bar{\boldsymbol{X}}_{T+k}$ using the BDF2 format of \eqref{BDF-explo-eq};
            \STATE Taking $\boldsymbol{X}_{T+k-1}$, $\boldsymbol{X}_{T+k}$, $s_{T+k-1}$, $s_{T+k}$, $\tau$ and $\bar{\Gamma}_{T+k+1}=\{\boldsymbol{X}_{T+k}^{i}|i=1, \ldots, N_{\boldsymbol{x}}\}$ as inputs, calculate $\boldsymbol{X}_{T+k+1}$, $s_{T+k+1}$ using the BDF2 format of \eqref{Correction-eq-p};
            \STATE Calculate $\varepsilon_{s} = |\max(s_{T+k+1}) - \min(s_{T+k+1})|$;
            \STATE Let $k = k+1$;
        \ENDWHILE
        
        \STATE \textbf{return} $\boldsymbol{X}_{c}, s_{c}$.
    \end{algorithmic}
\end{algorithm}

\section{Numerical Experiments}\label{Sec-Numerical}
In this section, we present results of numerical experiments to show the convergence, feasibility and efficiency of the proposed method, with particular focus on the improvement to the distribution of points in the evolution. 
Moreover, we apply the proposed method to simulating the MCF motion by point clouds. 
All source code and datasets used in the experiments are available at \url{https://github.com/Poker-Pan/EMSL-SE}. 
We emphasis that in all the experiments, only point clouds are used and the triangles are merely for the display of our results. 

In addition, in the subsequent experiments, to demonstrate the performance of our method, we use the area element of the point, which is defined as follows
\begin{equation}\label{eq:Sp}
    \begin{aligned}
        S_{p}(\boldsymbol{x}_i) = \frac{1}{3}\sum\limits_{K_{j}\in\Lambda_{\triangle}(\boldsymbol{x}_{i})}(Area(K_{j})), 
    \end{aligned}
\end{equation}
where $\Lambda_{\triangle}(\boldsymbol{x}_{i})$ is a set composed of triangles. It is a Delaunay triangle formed by the points in $\Lambda(\boldsymbol{x}_{i})$, and it excludes triangles that do not contain fixed points in $\boldsymbol{x}_{i}$.

\subsection{Validation of temporal convergence}\label{sec:converge}
We test the convergence rate of the proposed method using a sphere $\Gamma(0)=\{\boldsymbol{x}\in\mathbb{R}^3, |\boldsymbol{x}|=\frac{1}{2}\}$ as the initial surface where the point cloud is sampled. The velocity field $\boldsymbol{v}(\boldsymbol{X}(\boldsymbol{x}, t),t)$ is given as follows:
\begin{equation*}\label{Time-con-v}
    \begin{aligned}
        \boldsymbol{v}(\boldsymbol{X}(\boldsymbol{x}, t),t) = \boldsymbol{X}(\boldsymbol{x}, t)(1 - |\boldsymbol{X}(\boldsymbol{x}, t)|).
    \end{aligned}
\end{equation*}
The error of the numerical solution at time $t = T$, as provided by the proposed method with step size $\tau$, number of vertices $N_{\boldsymbol{x}}$, is measured by
\begin{equation*}\label{Error-time}
    \begin{aligned}
        \varepsilon^{t}_{\boldsymbol{x}}(N_{\boldsymbol{x}}, \tau) =& \max\limits_{j=1, \dots, N_{\boldsymbol{x}}}  |\boldsymbol{X}_{N_{\boldsymbol{x}}, \tau}(\boldsymbol{x}_{j}, T) - \boldsymbol{X}_{N_{\boldsymbol{x}}, 2\tau}(\boldsymbol{x}_{j}, T)|, \\ 
        \varepsilon^{t}_{s}(N_{\boldsymbol{x}}, \tau) =& \max\limits_{j=1, \dots, N_{\boldsymbol{x}}}  |S_{N_{\boldsymbol{x}}, \tau}(\boldsymbol{x}_{j}, T) - S_{N_{\boldsymbol{x}}, 2\tau}(\boldsymbol{x}_{j}, T)|.
    \end{aligned}
\end{equation*}

To evaluate the temporal convergence, we use a point cloud sufficiently sampled on the initial surface with a total number of points $N_{\boldsymbol{x}} = 30054$. 
Additionally, the initial condition is set to $s(\boldsymbol{x}, 0)=0$. 
The experimental results are presented in Table~\ref{tab:Error-con-time}. 
The temporal convergence rate of our method is consistent with the theoretical results corresponding to the BDF scheme. 
\begin{table}[htbp]
    \centering
    \caption{The time convergence test of the spatial point $\boldsymbol{x}$ and $s$ with $T = 1.0$, $N_{\boldsymbol{x}}=30054$ and $\eta=100$.}
    \label{tab:Error-con-time}
    \begin{tabular}{ccccccc}
        \hline
         & \multicolumn{2}{c}{BDF1} & \multicolumn{2}{c}{BDF2} & \multicolumn{2}{c}{BDF3}\\
        $1/\tau$ & $\varepsilon^{t}_{\boldsymbol{x}}$ & ord & $\varepsilon^{t}_{\boldsymbol{x}}$ &ord & $\varepsilon^{t}_{\boldsymbol{x}}$ & ord \\
        \hline
        $160$  & $1.47\times10^{-4}$ & -      & $4.59\times10^{-6}$ & -      & $6.25\times10^{-8}$  & -    \\
        $320$  & $7.37\times10^{-5}$ & $1.00$ & $1.16\times10^{-6}$ & $1.98$ & $7.76\times10^{-9}$  & 3.01 \\
        $640$  & $3.68\times10^{-5}$ & $1.00$ & $2.92\times10^{-7}$ & $1.99$ & $9.70\times10^{-10}$  & 3.00 \\
        $1280$ & $1.84\times10^{-5}$ & $1.00$ & $7.34\times10^{-8}$ & $2.00$ & $1.20\times10^{-10}$ & 3.01 \\
        \hline
        $1/\tau$ & $\varepsilon^{t}_{s}$ & ord & $\varepsilon^{t}_{s}$ &ord & $\varepsilon^{t}_{s}$ & ord \\
        \hline
        $160$  & $1.33\times10^{-3}$ & -      & $4.03\times10^{-5}$ & -      & $4.67\times10^{-7}$  & -    \\
        $320$  & $6.68\times10^{-4}$ & $1.00$ & $1.00\times10^{-5}$ & $2.00$ & $5.86\times10^{-8}$  & 3.00 \\
        $640$  & $3.33\times10^{-4}$ & $1.00$ & $2.52\times10^{-6}$ & $2.00$ & $7.34\times10^{-9}$  & 3.00 \\
        $1280$ & $1.66\times10^{-4}$ & $1.00$ & $6.31\times10^{-7}$ & $2.00$ & $9.09\times10^{-10}$ & 3.01 \\
        \hline
    \end{tabular}
\end{table}

Under this velocity, the surface $\Gamma(t)$ is a sphere centered at the origin with radius $r$. 
Considering only the radial component and neglecting tangential motion of points on the spherical surface, the radius $r(t)$ satisfies the following equation:
\begin{equation*}\label{Exact-time}
    \begin{aligned}
        \frac{d r}{dt} = r(1-r).
    \end{aligned}
\end{equation*}
Then the solution satisfying the initial value condition $r(0)=\frac{1}{2}$ is $r(t) = \frac{1}{1+e^{-t}}$. 
The error of radius at time $t = T$, as provided by the proposed method with step size $\tau$ and number of vertices $N_{\boldsymbol{x}}$, is measured by
\begin{equation*}\label{Error-time-r}
    \begin{aligned}
        \varepsilon^{t}_{r}(N_{\boldsymbol{x}}, \tau) = \max\limits_{j=1, \dots, N_{\boldsymbol{x}}} \left| \|\boldsymbol{X}_{N_{\boldsymbol{x}}, \tau}(\boldsymbol{x}_{j}, T)\|_{2} - r(T)\right|.
    \end{aligned}
\end{equation*}
The experimental results are presented in Table~\ref{tab:Error-con-time-r}. 
The temporal convergence of the radius is also consistent with the theoretical analysis of the BDF scheme. 
\begin{table}[htbp]
    \centering
    \caption{The time convergence test of radius $r$ with $T = 1.0$, $N_{\boldsymbol{x}}=30054$ and $\eta=100$.}
    \label{tab:Error-con-time-r}
    \begin{tabular}{ccccccc}
        \hline
         & \multicolumn{2}{c}{BDF1} & \multicolumn{2}{c}{BDF2} & \multicolumn{2}{c}{BDF3}\\
        $1/\tau$ & $\varepsilon^{t}_{r}$ & ord & $\varepsilon^{t}_{r}$ &ord & $\varepsilon^{t}_{r}$ & ord \\
        \hline
        $80$   & $294\times10^{-4}$ & -      & $614\times10^{-6}$ & -      & $7.13\times10^{-8}$  & -    \\
        $160$  & $147\times10^{-4}$ & $1.00$ & $155\times10^{-6}$ & $1.98$ & $8.84\times10^{-9}$  & 3.01 \\
        $320$  & $737\times10^{-5}$ & $1.00$ & $390\times10^{-7}$ & $1.99$ & $1.08\times10^{-9}$  & 3.03 \\
        $640$  & $368\times10^{-5}$ & $1.00$ & $979\times10^{-8}$ & $2.00$ & $1.35\times10^{-10}$ & 2.99 \\
        $1280$ & $184\times10^{-5}$ & $1.00$ & $245\times10^{-8}$ & $2.00$ & $1.68\times10^{-11}$ & 3.01 \\
        \hline
    \end{tabular}
\end{table}

\subsection{Validation of spatial convergence}
We evaluate the spatial error and convergence rate of the proposed method using a torus as the initial surface. 
The velocity field $\boldsymbol{v}(\boldsymbol{X}(\boldsymbol{x}, t), t)$ is defined as follows:
\begin{equation*}\label{Space-con-v}
    \begin{aligned}
        \boldsymbol{v}(\boldsymbol{X}(\boldsymbol{x}, t),t) = 
        500\left(\begin{array}{c}
             \cos{(\pi x)}\sin{(\pi y)}\sin{(\pi z)}  \\
             \sin{(\pi x)}\cos{(\pi y)}\sin{(\pi z)}  \\
             \sin{(\pi x)}\sin{(\pi y)}\cos{(\pi z)}
        \end{array}\right), 
    \end{aligned}
\end{equation*}
initial density value $\rho_{0}(\boldsymbol{x})=1$.

The error of the numerical solution at time $t = T$, as provided by the proposed method with step size $\tau$, number of vertices $N_{\boldsymbol{x}}$, is measured by
\begin{equation*}\label{Error-space}
    \begin{aligned}
        \varepsilon^{s}_{\boldsymbol{x}}(N_{\boldsymbol{x}}, \tau) =& \max\limits_{j=1, \dots, N_{\boldsymbol{x}}}  |\boldsymbol{X}_{N_{\boldsymbol{x}}, \tau}(\boldsymbol{x}_{j}, T) - \boldsymbol{X}_{4N_{\boldsymbol{x}}, \tau}(\boldsymbol{x}_{j}, T)|, \\ 
        \varepsilon^{s}_{s}(N_{\boldsymbol{x}}, \tau) =& \max\limits_{j=1, \dots, N_{\boldsymbol{x}}}  |S_{N_{\boldsymbol{x}}, \tau}(\boldsymbol{x}_{j}, T) - S_{4N_{\boldsymbol{x}}, \tau}(\boldsymbol{x}_{j}, T)|. 
    \end{aligned}
\end{equation*}
To evaluate spatial convergence, we employ a sufficiently small time step size $\tau = 1\times10^{-4}$ to ensure that the temporal discretization error is negligible. 
\begin{table}[htbp]
    \small
    \centering
    \setlength{\tabcolsep}{0.14cm}
    \caption{The spatial convergence test of the spatial point $\boldsymbol{x}$ and $s$ with $T = 1.0$, $\tau = 1\times10^{-4}$ by Algorithm~\ref{BDF2-Algorithm}.}
    \label{tab:Error-con-space}
    \begin{tabular}{ccccccccccc}
        \hline
        \multicolumn{2}{c}{$\eta$} & \multicolumn{3}{c}{1} & \multicolumn{3}{c}{10} & \multicolumn{3}{c}{100} \\
        \hline
        $N_{\boldsymbol{x}}$ & $h$ & $\varepsilon^{s}_{\boldsymbol{x}}$ & $ord_{N_{\boldsymbol{x}}}$ & $ord_{h}$ & $\varepsilon^{s}_{\boldsymbol{x}}$ & $ord_{N_{\boldsymbol{x}}}$ & $ord_{h}$ & $\varepsilon^{s}_{\boldsymbol{x}}$ & $ord_{N_{\boldsymbol{x}}}$ & $ord_{h}$  \\
        \hline
        $256$   & $0.450$ & $7.77\times10^{-5}$ & -      & -      & $7.21\times10^{-4}$ & -      & -      & $4.27\times10^{-3}$ & -      & -      \\
        $1024$  & $0.237$ & $5.06\times10^{-5}$ & $0.31$ & $0.67$ & $5.00\times10^{-4}$ & $0.26$ & $0.57$ & $2.98\times10^{-3}$ & $0.26$ & $0.56$ \\
        $4096$  & $0.120$ & $1.62\times10^{-5}$ & $0.82$ & $1.68$ & $1.30\times10^{-4}$ & $0.97$ & $1.98$ & $5.47\times10^{-4}$ & $1.22$ & $2.50$ \\
        $16384$ & $0.060$ & $5.14\times10^{-6}$ & $0.83$ & $1.66$ & $4.17\times10^{-5}$ & $0.82$ & $1.65$ & $1.56\times10^{-4}$ & $0.90$ & $1.81$ \\
        \hline
        $N_{\boldsymbol{x}}$ & $h$ & $\varepsilon^{s}_{s}$ & $ord_{N_{\boldsymbol{x}}}$ & $ord_{h}$ & $\varepsilon^{s}_{s}$ & $ord_{N_{\boldsymbol{x}}}$ & $ord_{h}$ & $\varepsilon^{s}_{s}$ & $ord_{N_{\boldsymbol{x}}}$ & $ord_{h}$  \\
        \hline
        $256$   & $0.450$ & $6.45\times10^{-2}$ & -      & -      & $6.53\times10^{-2}$ & -      & -      & $6.70\times10^{-2}$ & -      & -      \\
        $1024$  & $0.237$ & $1.97\times10^{-2}$ & $0.85$ & $1.85$ & $2.02\times10^{-2}$ & $0.85$ & $1.83$ & $1.81\times10^{-2}$ & $0.94$ & $2.04$ \\
        $4096$  & $0.120$ & $4.37\times10^{-3}$ & $1.09$ & $2.22$ & $3.24\times10^{-3}$ & $1.32$ & $2.70$ & $2.37\times10^{-3}$ & $1.47$ & $3.00$ \\
        $16384$ & $0.060$ & $1.11\times10^{-3}$ & $0.99$ & $1.99$ & $7.40\times10^{-4}$ & $1.06$ & $2.14$ & $6.61\times10^{-4}$ & $0.92$ & $1.85$ \\
        \hline
    \end{tabular}
\end{table}

Table~\ref{tab:Error-con-space} presents the errors and convergence rate for the spatial variables $\boldsymbol{x}$ and $s$ under various choices of $\eta$. 
As the number of spatial points increases, the errors consistently decrease. 
The estimated convergence rate with respect to the step size $h$ is approximately twice that of the spatial rate in $\boldsymbol{x}$, aligning well with theoretical expectations. 
However, it is observed that the convergence rate with respect to $h$ fluctuates around $2$, which is reasonable taking into account that the unstructured points allocation may cause perturbation of errors in computations.

\subsection{Validation of point distribution improvement}
\begin{figure}[htbp]
    \begin{center}
    \subfigure[Initial surface]{\includegraphics[width=3.2cm]{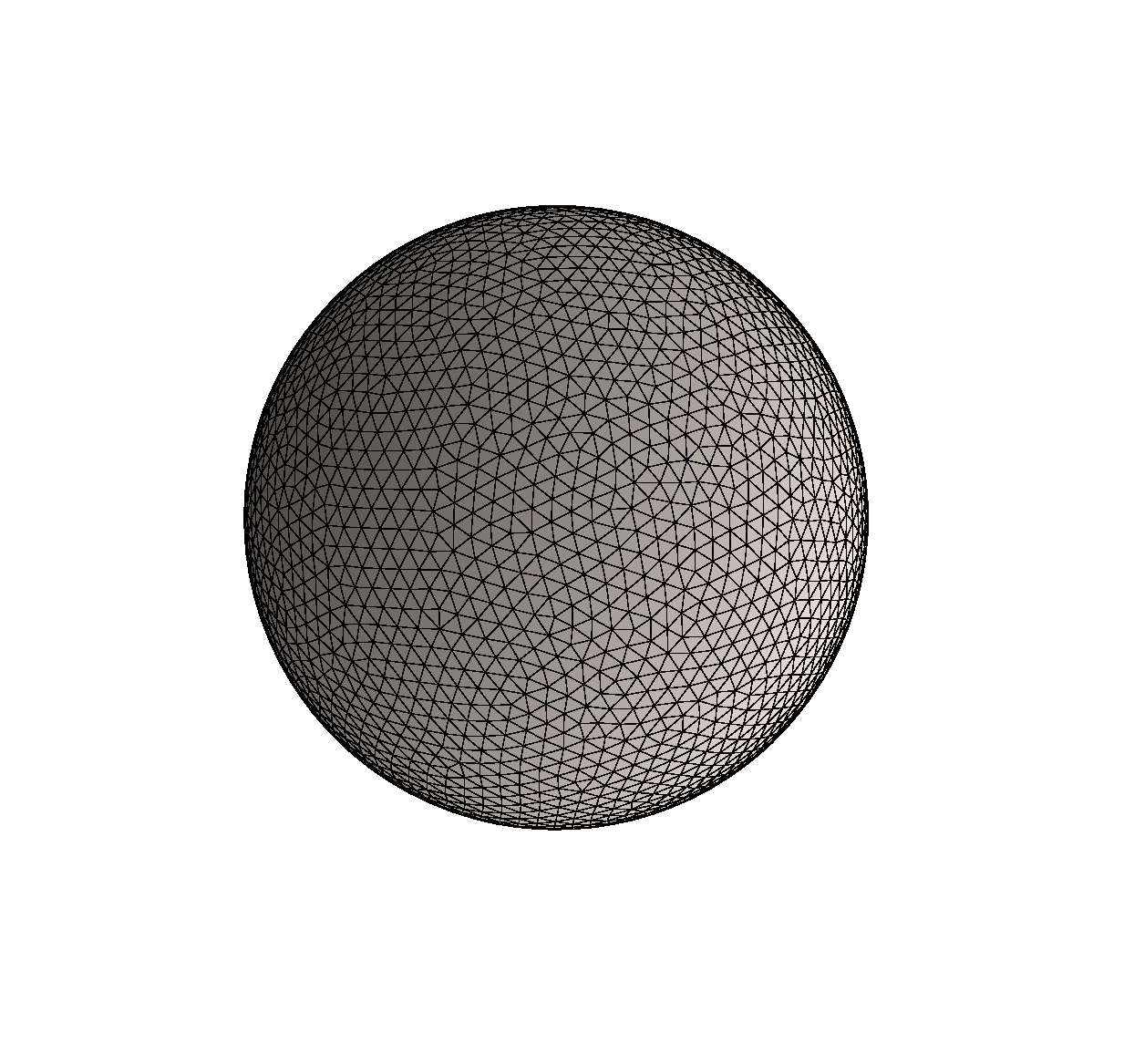}}
    \subfigure[$\boldsymbol{v}(\boldsymbol{x},t)$]{\includegraphics[width=3.2cm]{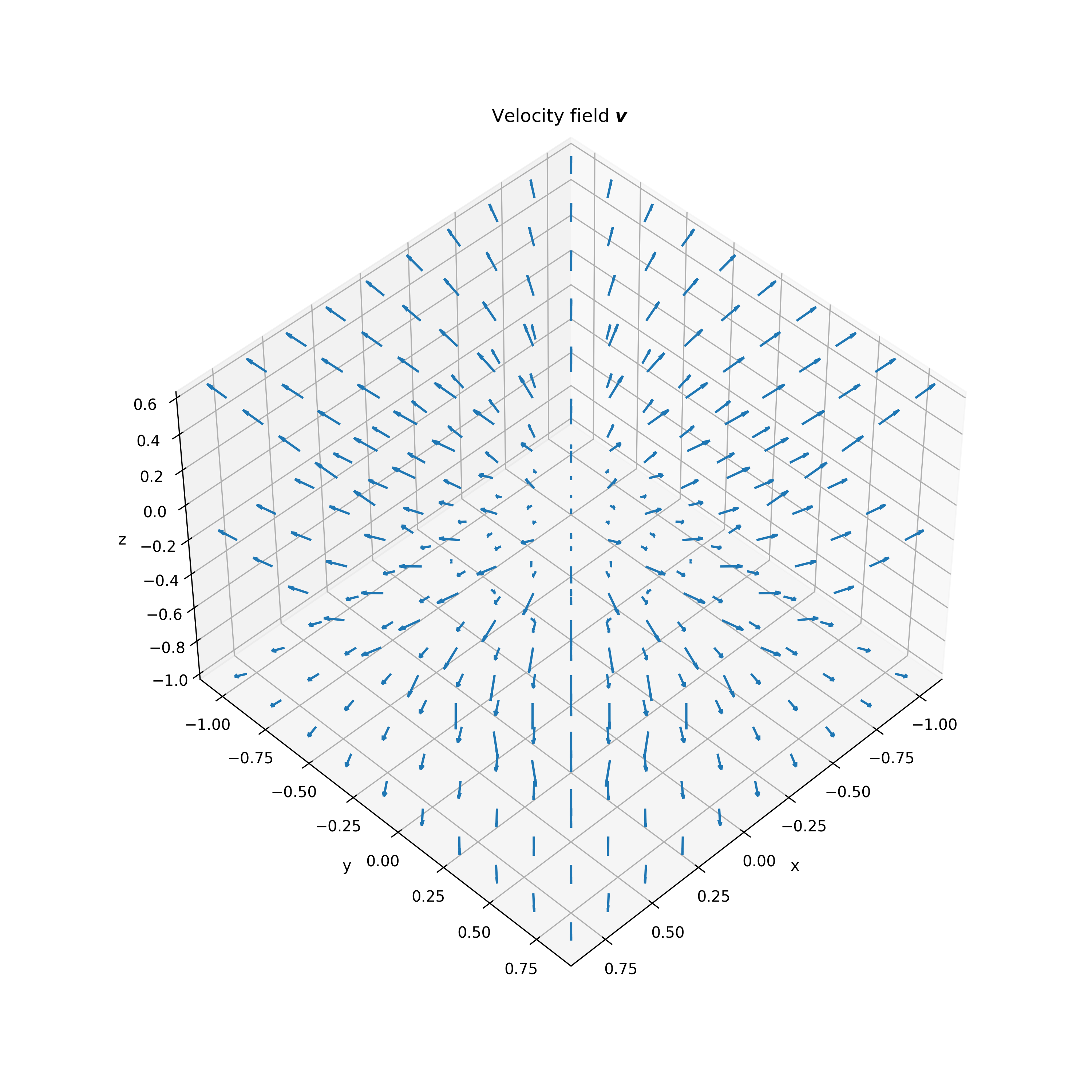}}\\
    \subfigure[$\eta=0$]{\includegraphics[width=3.2cm]{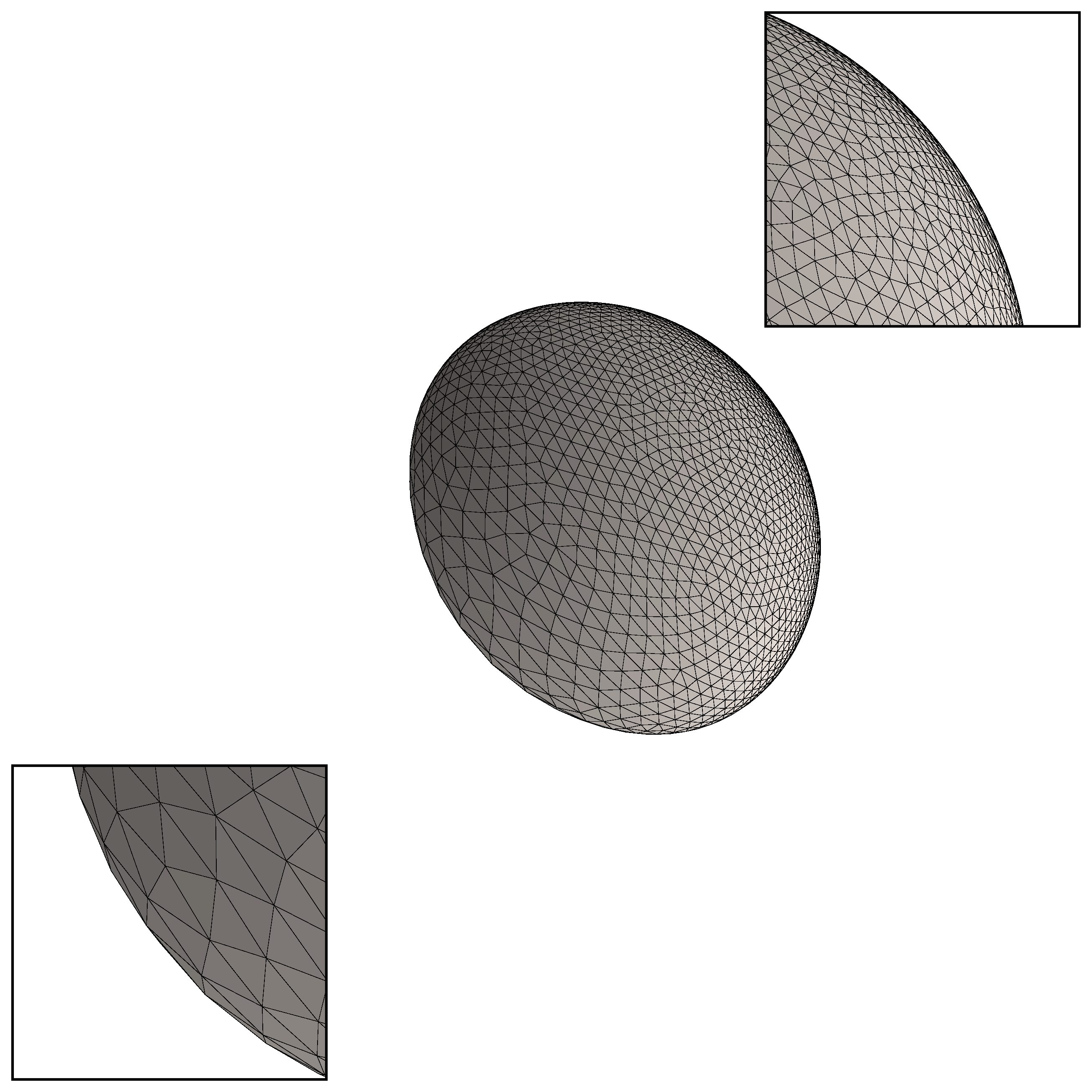}}
    \subfigure[$\eta=1$]{\includegraphics[width=3.2cm]{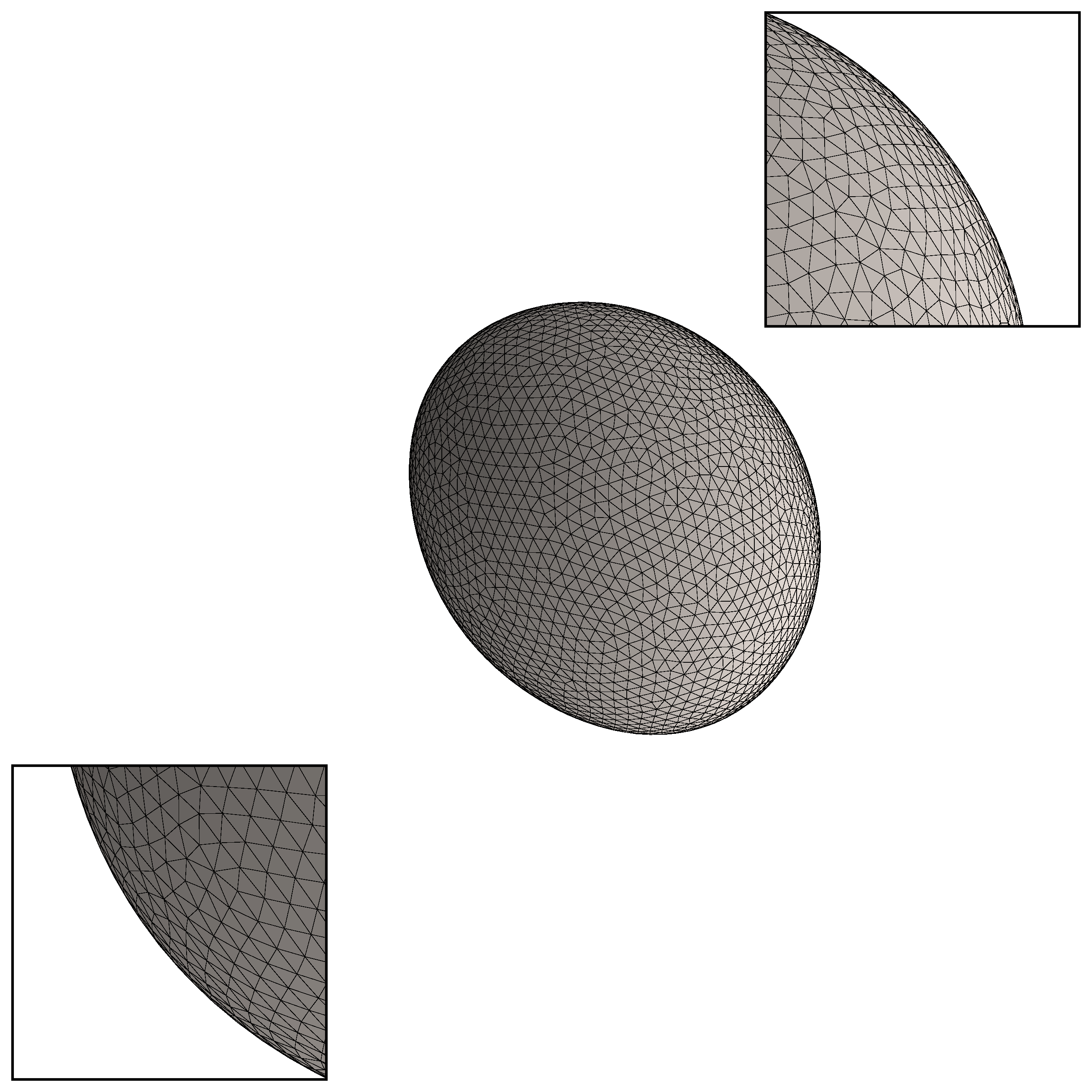}}
    \subfigure[$\eta=10$]{\includegraphics[width=3.2cm]{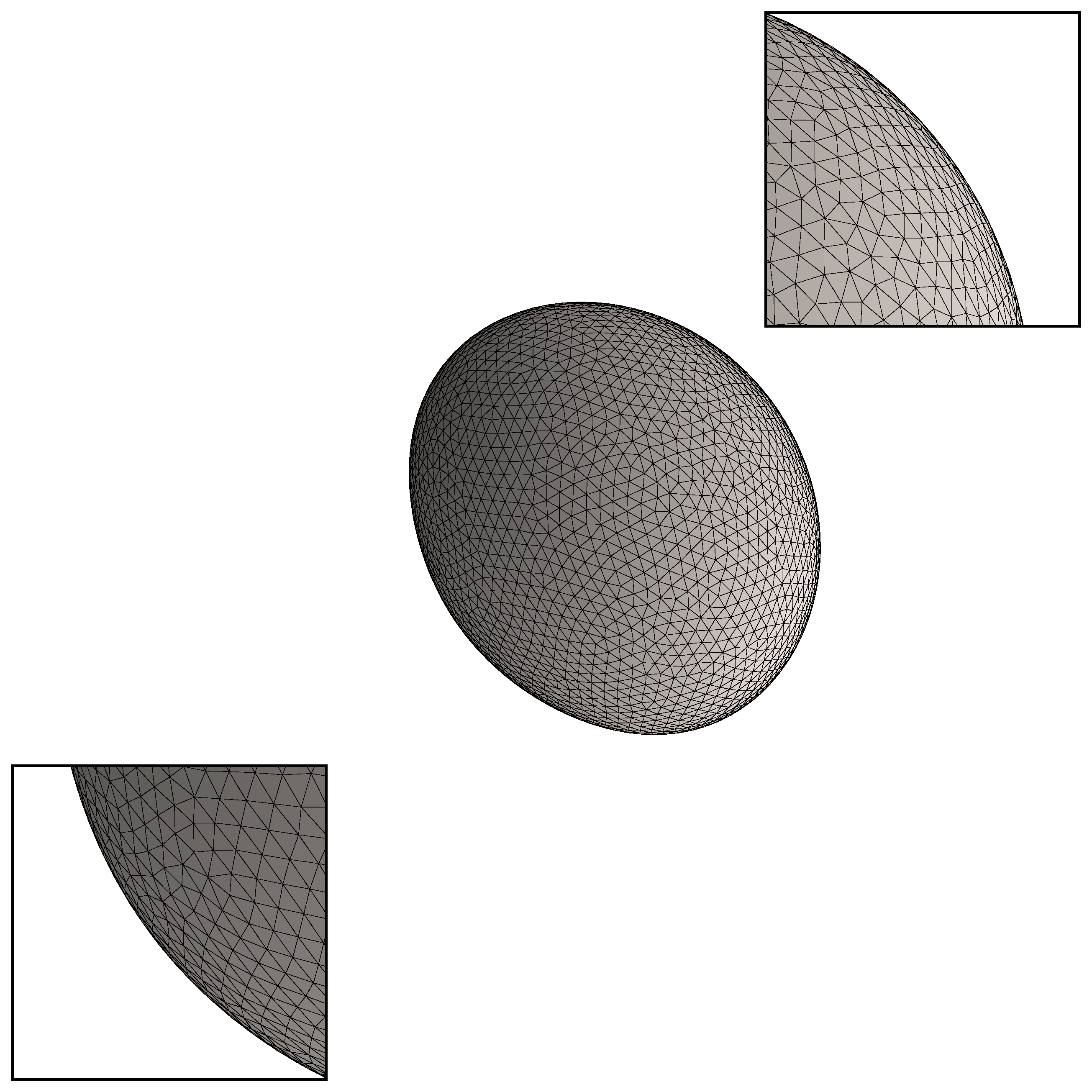}}
    \subfigure[$\eta=100$]{\includegraphics[width=3.2cm]{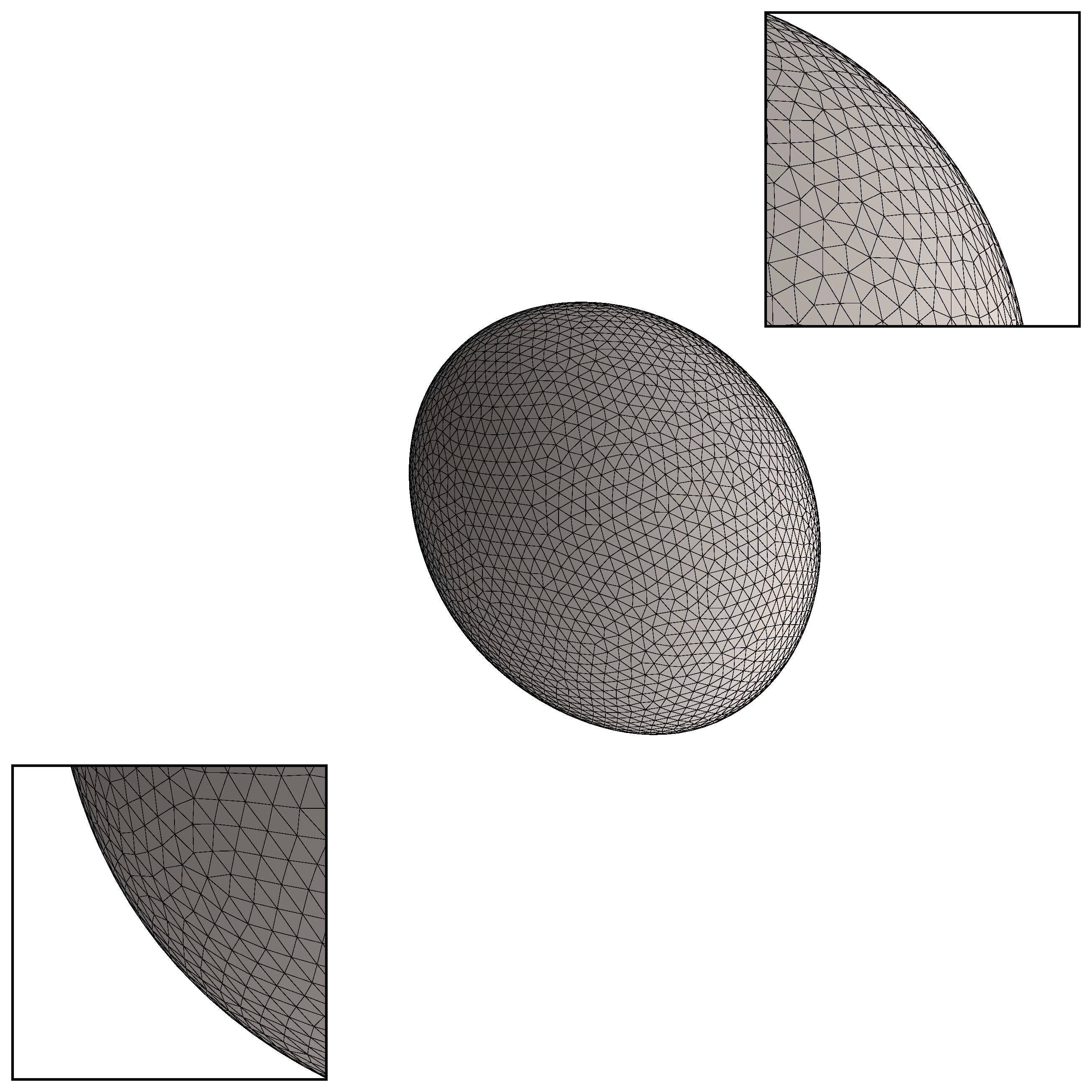}}\\
    \caption{Surfaces at $T=2$ computed by BDF2 method (Algorithm~\ref{BDF2-Algorithm}) with different $\eta$.}
    \label{imp-shpere}
    \end{center}
\end{figure}

{\bf Spherical surface with simple velocity.} 
We consider the same velocity field in the surface evolution as the example in Section \ref{sec:converge}, that is:
\begin{equation*}\label{imp-sim-v-sphere}
    \begin{aligned}
    \begin{aligned}
        \boldsymbol{v}(\boldsymbol{X}(\boldsymbol{x}, t),t) = \boldsymbol{X}(\boldsymbol{x}, t)(1 - |\boldsymbol{X}(\boldsymbol{x}, t)|).
    \end{aligned}
    \end{aligned}
\end{equation*}
Visually it is displayed in Figure~\ref{imp-shpere}(b). 

However, the initial surface is different to the one in Section \ref{sec:converge}. Here $\Gamma(t)=\{\boldsymbol{x}=(x,y,z)\in\mathbb{R}^{3}: \varphi(\boldsymbol{x}, t)=\frac{1}{2}\}$ is described by a level set function
\begin{equation*}\label{imp-levelset-shpere}
    \begin{aligned}
        \varphi(\boldsymbol{x}, t) = (x-\frac{1}{4})^2 + (y-\frac{1}{4})^2 +(z-\frac{1}{4})^2.
    \end{aligned}
\end{equation*}
It is a sphere with a radius of $r = 0.5$ and the center is shifted from the origin to $(0.25, 0.25, 0.25)$. 
Under such a configuration, due to the different magnitudes and directions of the velocities on both sides of the ball, there will be the phenomenon of point aggregation and sparsity without interaction to the surface evolution. 
We use the BDF2 scheme with $N_{\boldsymbol{x}}=2904$ and $\tau=1\times10^{-4}$ in our numerical tests. 

\begin{figure}
    \centering
    \includegraphics[width=0.55\linewidth]{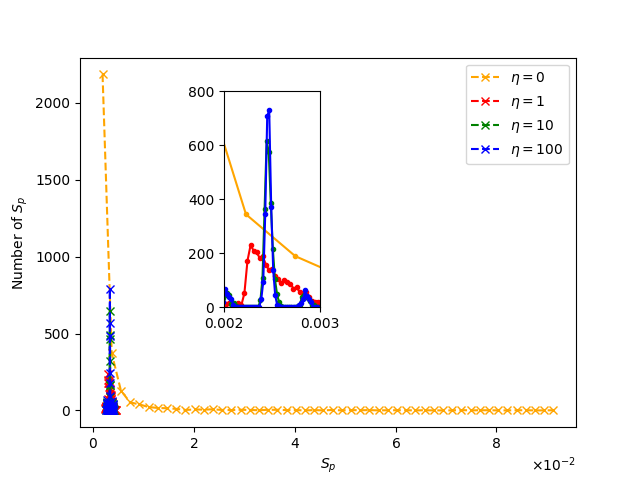}
    \caption{Distribution of area elements of each point, \eqref{eq:Sp}, for Algorithm~\ref{BDF2-Algorithm} with different $\eta$. The surface is evolved up to time $T = 2$.}
    \label{imp-shpere-s}
\end{figure}

Figure~\ref{imp-shpere} illustrates the distribution of spatial points $\boldsymbol{x}$ under various values of $\eta=0, 1, 10, 100$. 
In the cases of $\eta=1, 10, 100$, the points distribute more evenly on the sphere, and the point clouds show better shapes of the sphere comparing to the case of $\eta=0$, that is without the artificial tangential velocity introduced to the evolution.
Furthermore, we present the histogram of $S_{p}$ in Figure~\ref{imp-shpere-s}. 
It shows that adding the tangential velocity is able to have more desirable distribution of $S_{p}$, and therefore a more uniform density of the point cloud.

\begin{figure}[htbp]
    \begin{center}
    \subfigure[Initial surface]{\includegraphics[width=3.2cm]{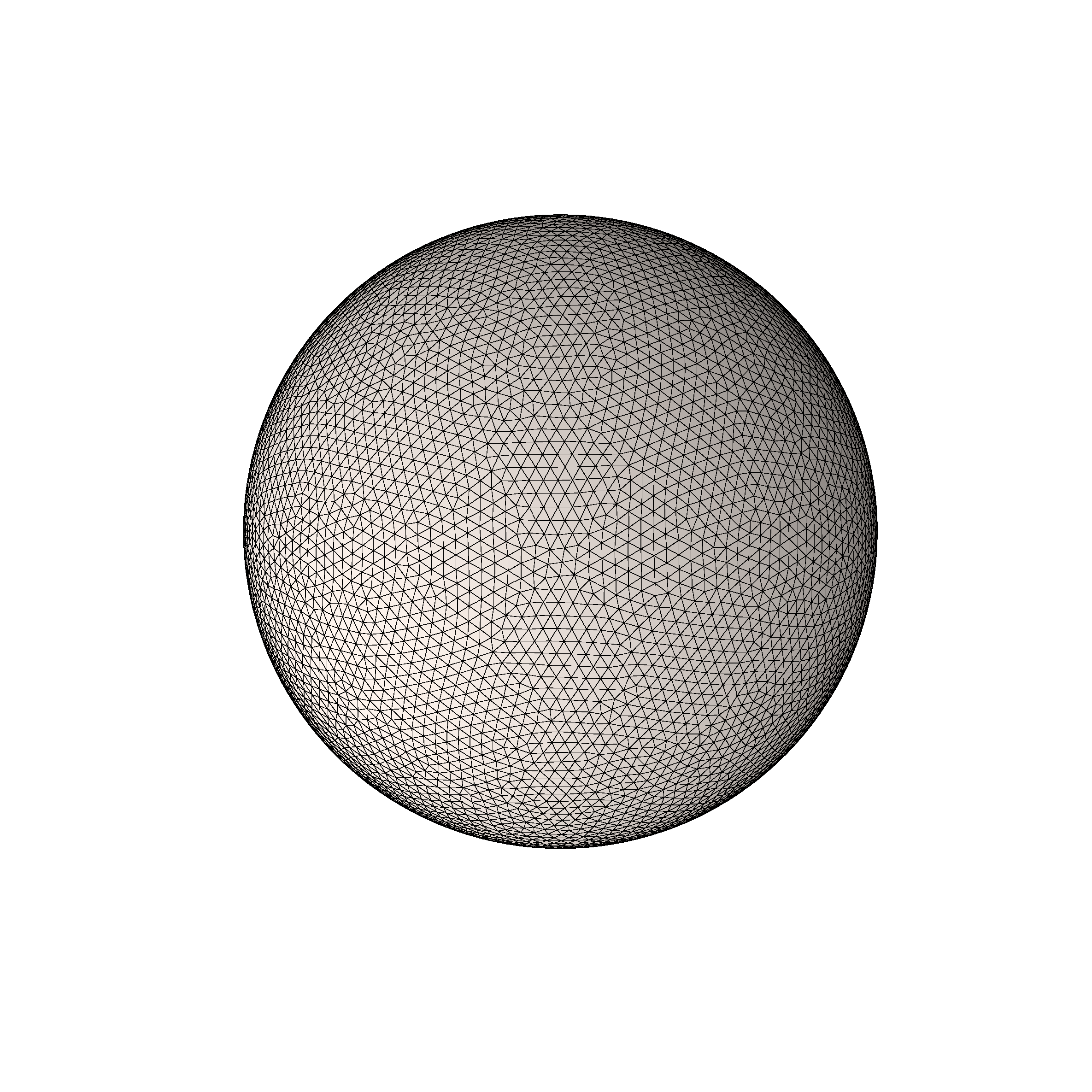}}\\
    \subfigure[$\eta=0$]{\includegraphics[width=3.2cm]{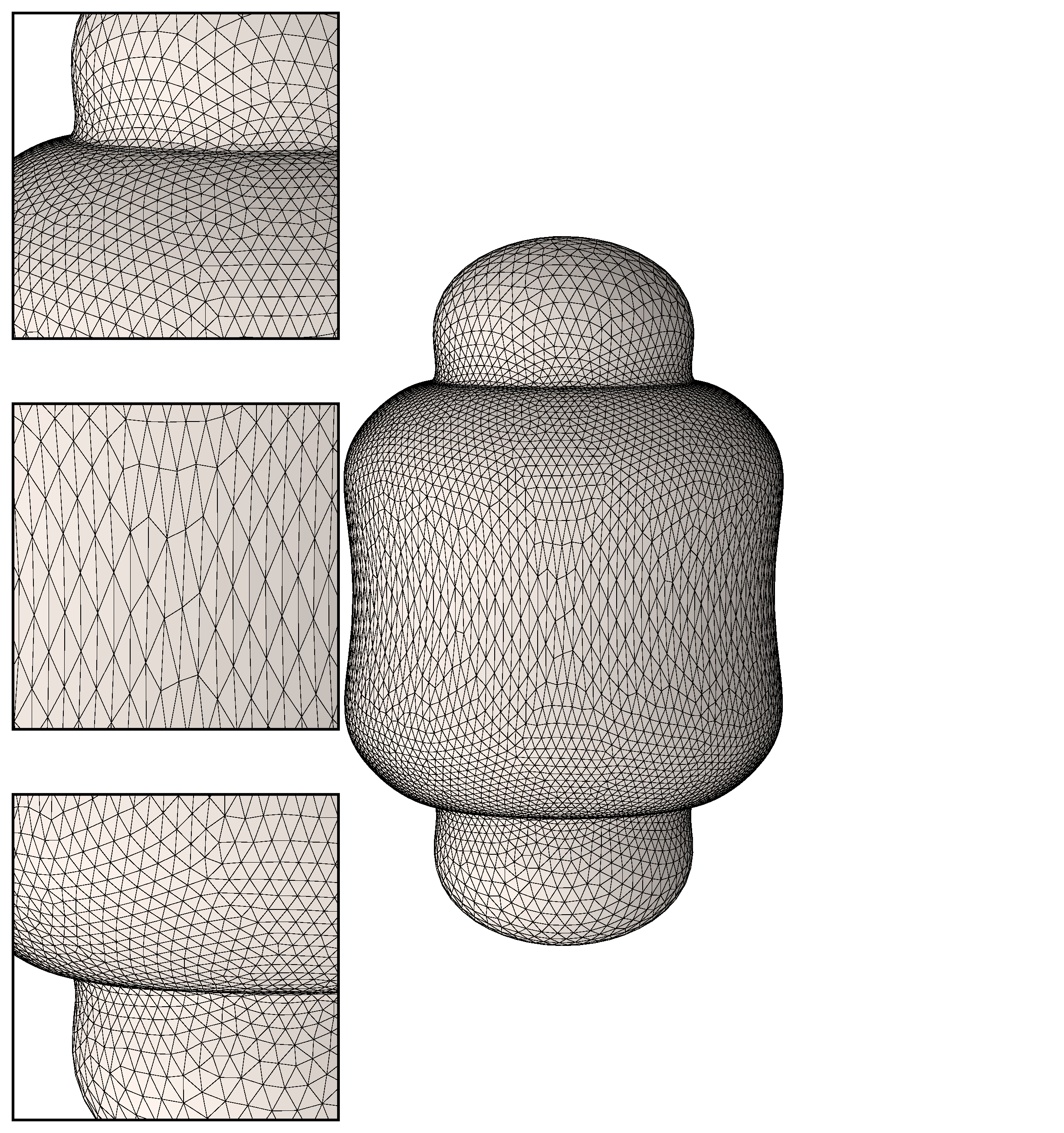}}
    \subfigure[$\eta=1$]{\includegraphics[width=3.2cm]{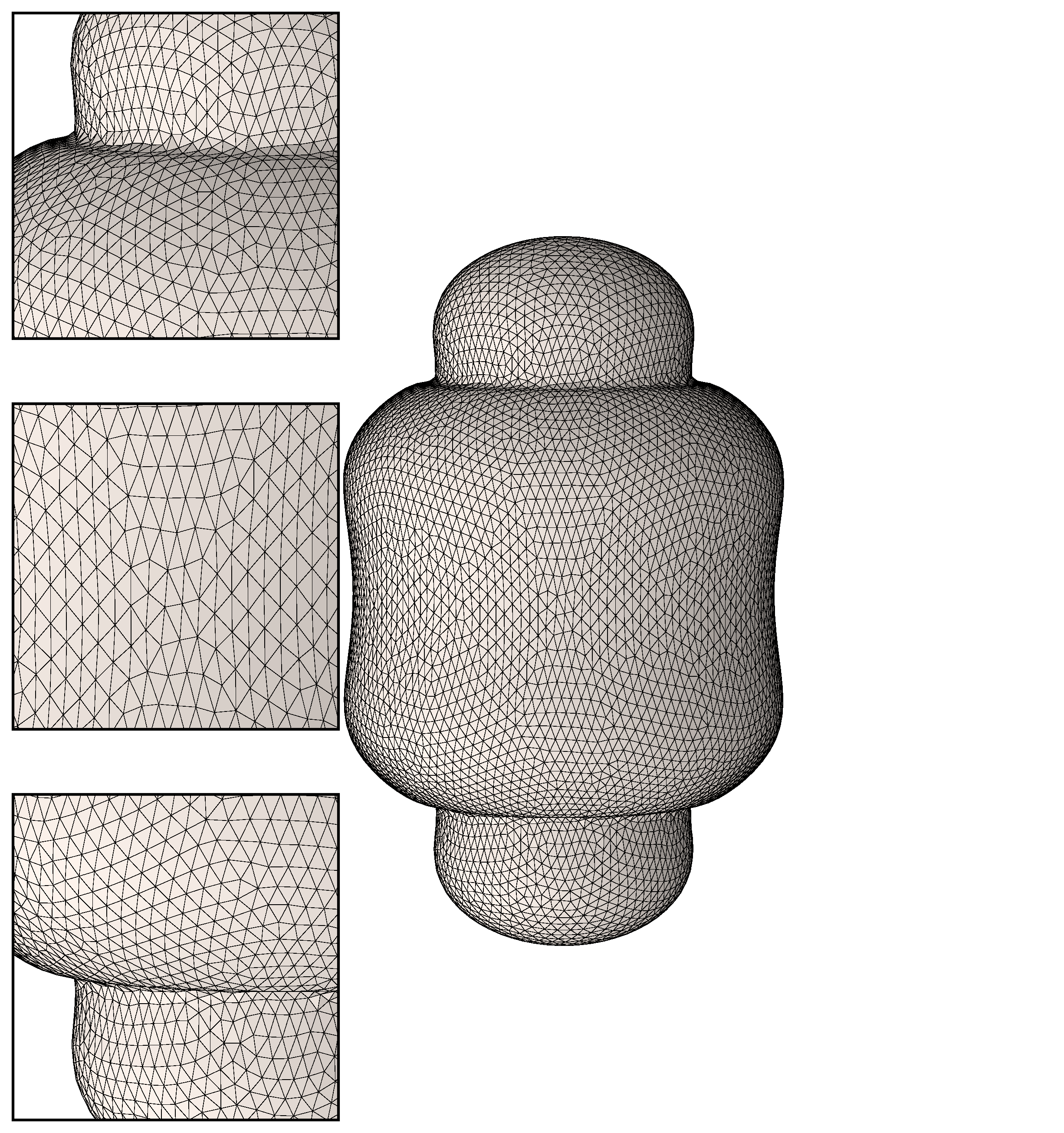}}
    \subfigure[$\eta=10$]{\includegraphics[width=3.2cm]{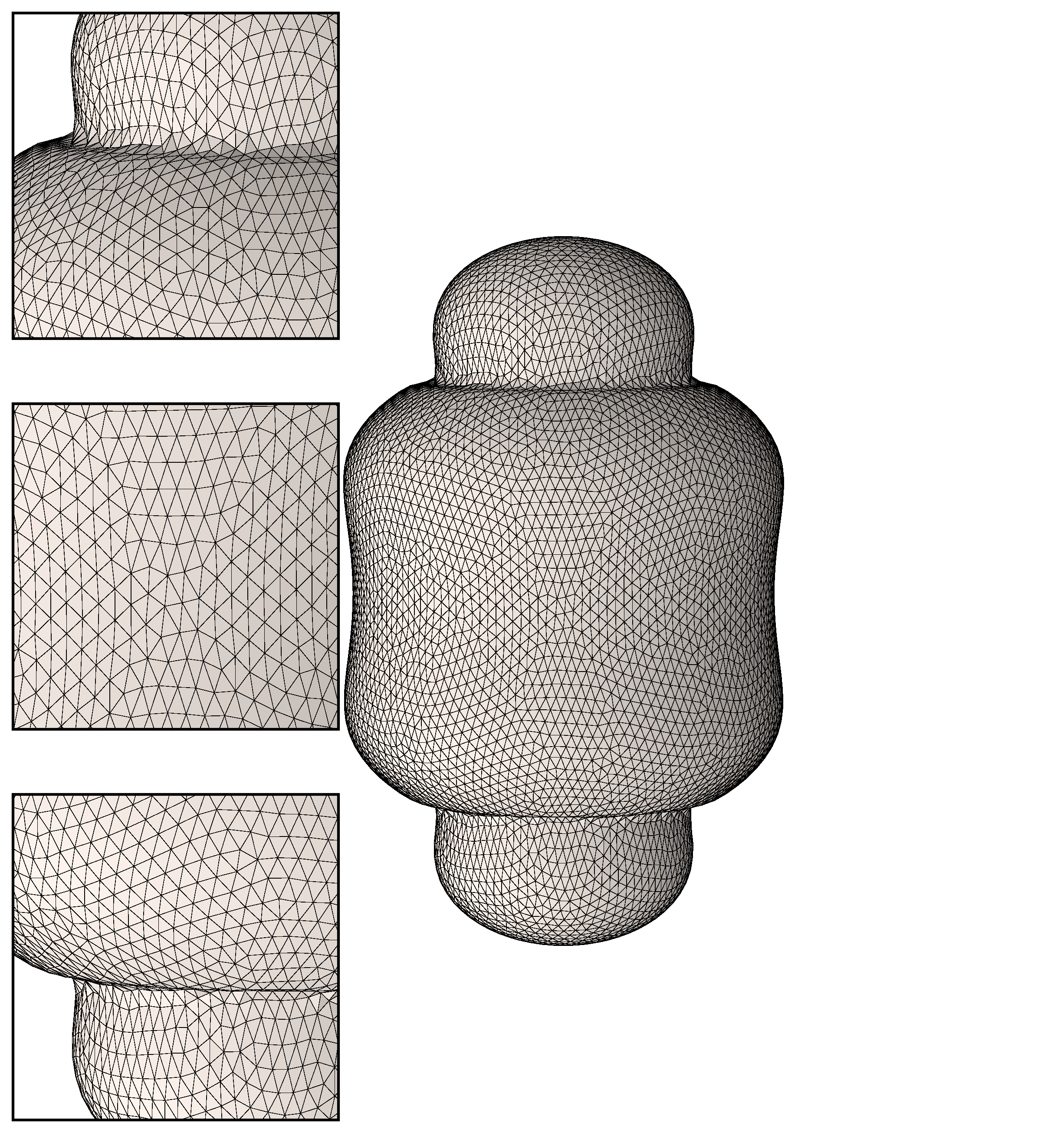}}
    \subfigure[$\eta=100$]{\includegraphics[width=3.2cm]{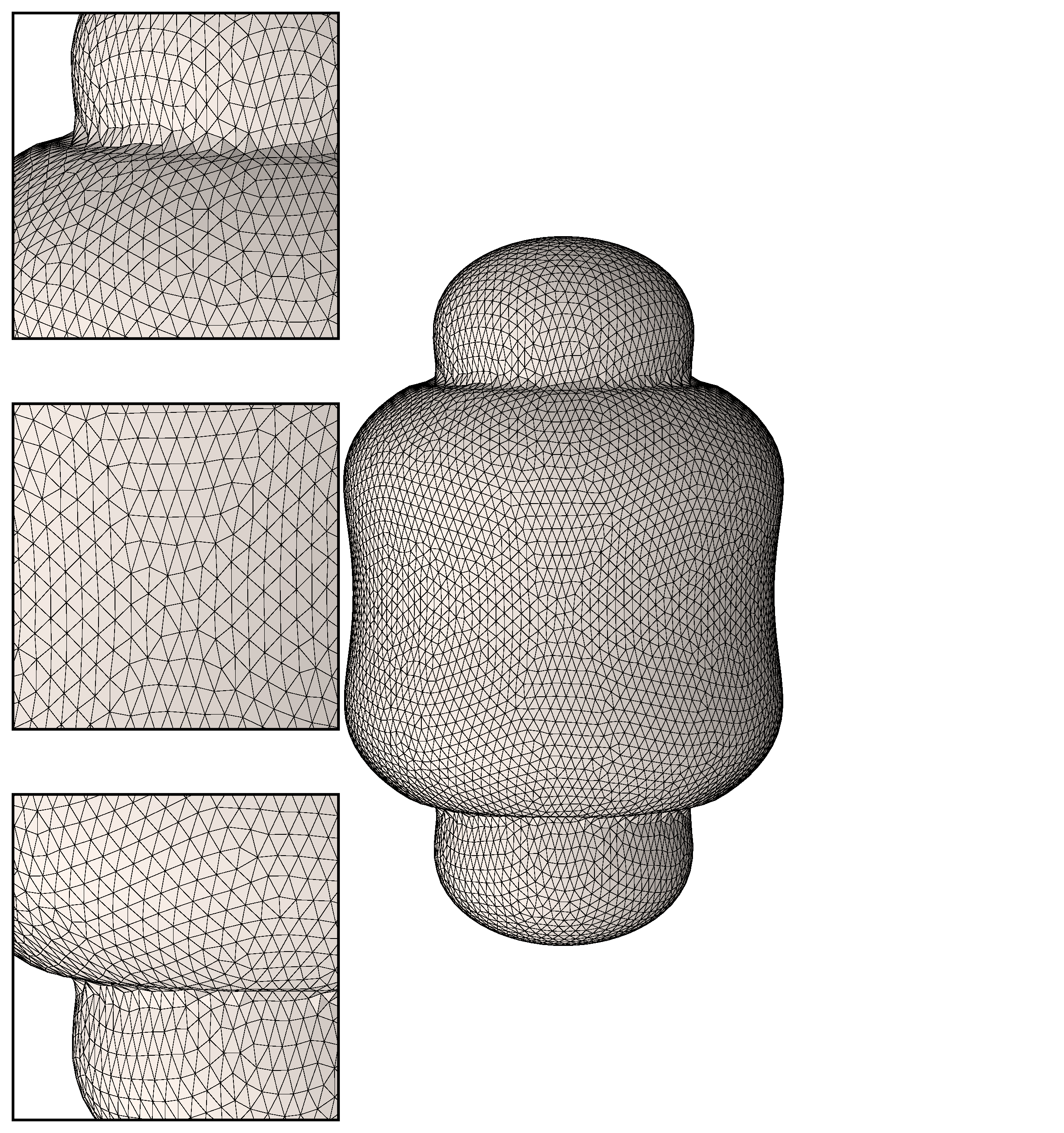}}\\
    \caption{Surfaces at $T=0.6$ computed by BDF2 method (Algorithm~\ref{BDF2-Algorithm}) with different $\eta$.}
    \label{imp-sphere-new}
    \end{center}
\end{figure}

{\bf Spherical surface with complex velocity.} 
We consider a scenario involving a complex velocity field coupled with a relatively simple surface evolution. 
The prescribed velocity field is defined as follows:
\begin{equation*}\label{imp-com-v-sphere}
    \begin{aligned}
        \boldsymbol{v}(\boldsymbol{X}(\boldsymbol{x}, t),t) = -\frac{\xi_{t}\nabla\xi}{|\nabla\xi|^2}.
    \end{aligned}
\end{equation*}
where $\xi(\boldsymbol{x}, t) = \frac{x^2}{a^2(t)} + \frac{y^2}{a^2(t)} + G(\frac{z^2}{b^2(t)})$, $G(s) = 200s(s-\frac{199}{200})$, $a(t) = 0.1+0.05\sin{2\pi t}$ and $b(t) = 1+0.2\sin{4\pi t}$. 
Evolving surface $\Gamma(t)=\{\boldsymbol{x}=(x,y,z)\in\mathbb{R}^{3}: \varphi(\boldsymbol{x}, t)=1\}$ is described by a level set function
\begin{equation*}\label{imp-com-levelset-shpere}
    \begin{aligned}
        \varphi(\boldsymbol{x}, t) = x^2 + y^2 + z^2.
    \end{aligned}
\end{equation*}
It is a unit sphere centered at the origin. 
If the tangential velocity component is omitted, spatial points exhibit significant clustering near the top and bottom regions by time $T = 0.6$, while points around the middle of the surface become sparsely distributed. 
These artifacts are clearly illustrated in Figure~\ref{imp-sphere-new}(b). 
We use again BDF2 scheme with $N_{\boldsymbol{x}}=7446$ and $\tau=1\times10^{-4}$ for surface evolution. 

\begin{figure}
    \centering
    \includegraphics[width=0.55\linewidth]{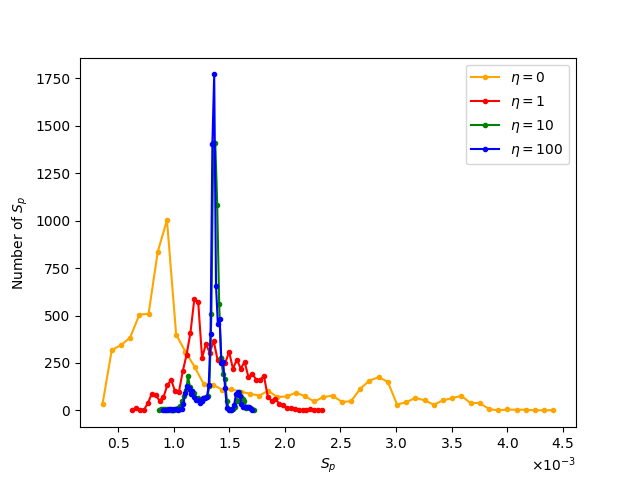}
    \caption{Distribution of area element of each point, \eqref{eq:Sp}, for Algorithm~\ref{BDF2-Algorithm} with different $\eta$. Surface is evolved up to time $T = 0.6$.}
    \label{imp-shpere-new-s}
\end{figure}

In Figure~\ref{imp-sphere-new}, we present the results computed using Algorithm~\ref{BDF2-Algorithm}. 
The results indicate that Algorithm~\ref{BDF2-Algorithm} effectively evolves the surface geometry while maintaining reasonable point distribution. 
Furthermore, we present the histogram of $S_{p}$ in Figure~\ref{imp-shpere-new-s}, where similar performance to the previous experiment is reported.

\begin{figure}[htbp]
    \begin{center}
    \subfigure[Initial surface]{\includegraphics[width=3.2cm]{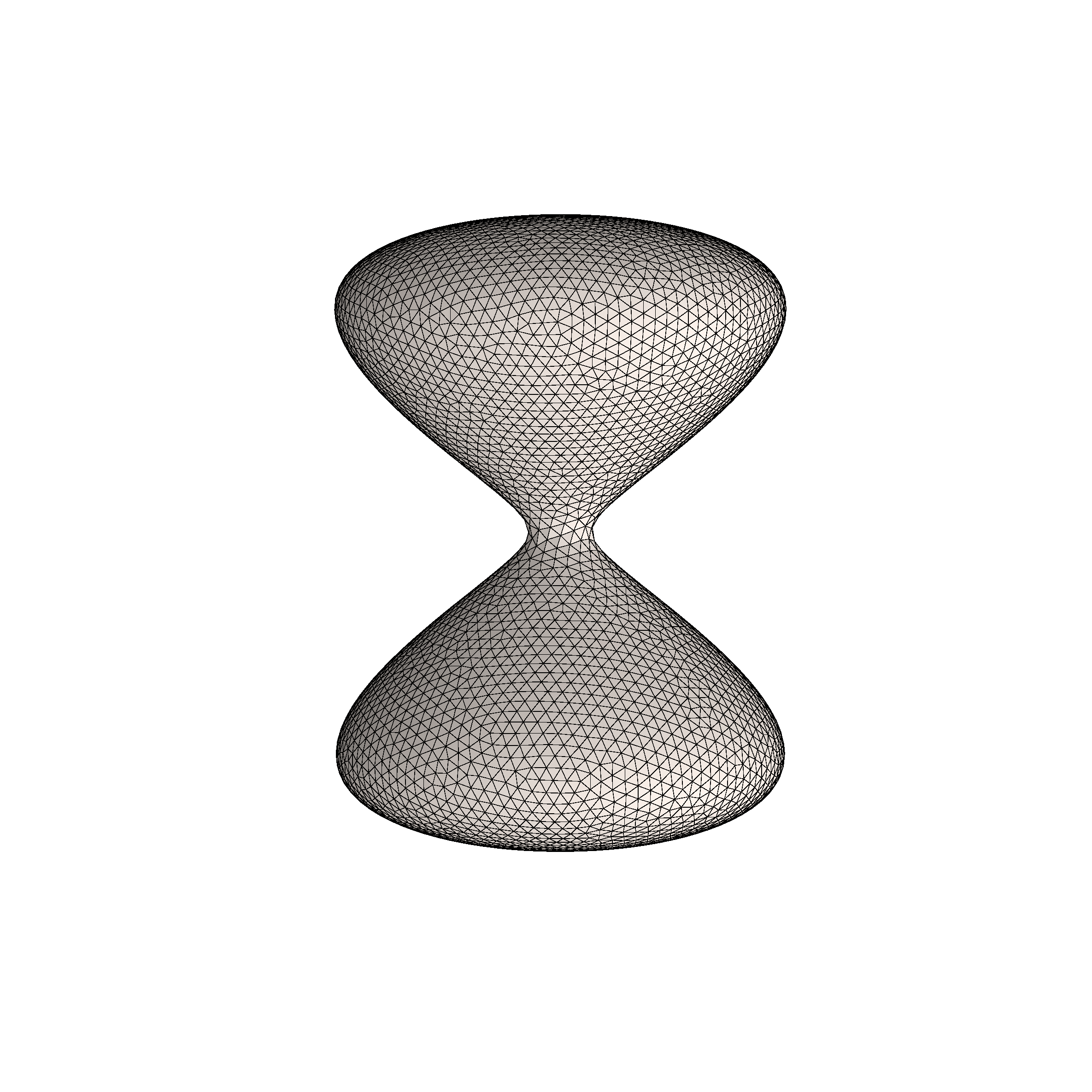}}\\
    \rotatebox{90}{~~~~~~~~Coupled model}
    \subfigure[$\eta=0$]{\includegraphics[width=3.2cm]{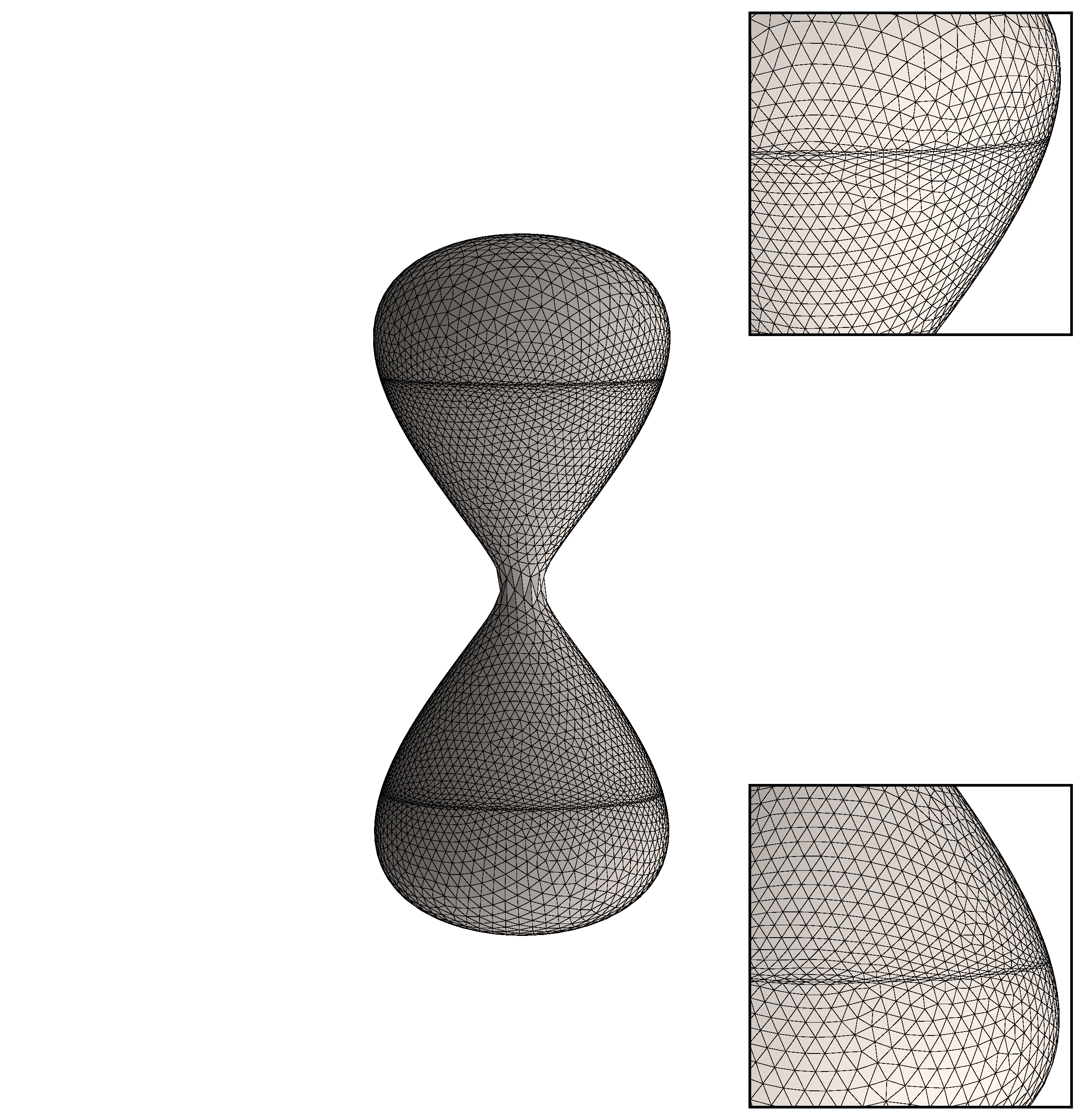}}
    \subfigure[$\eta=1$]{\includegraphics[width=3.2cm]{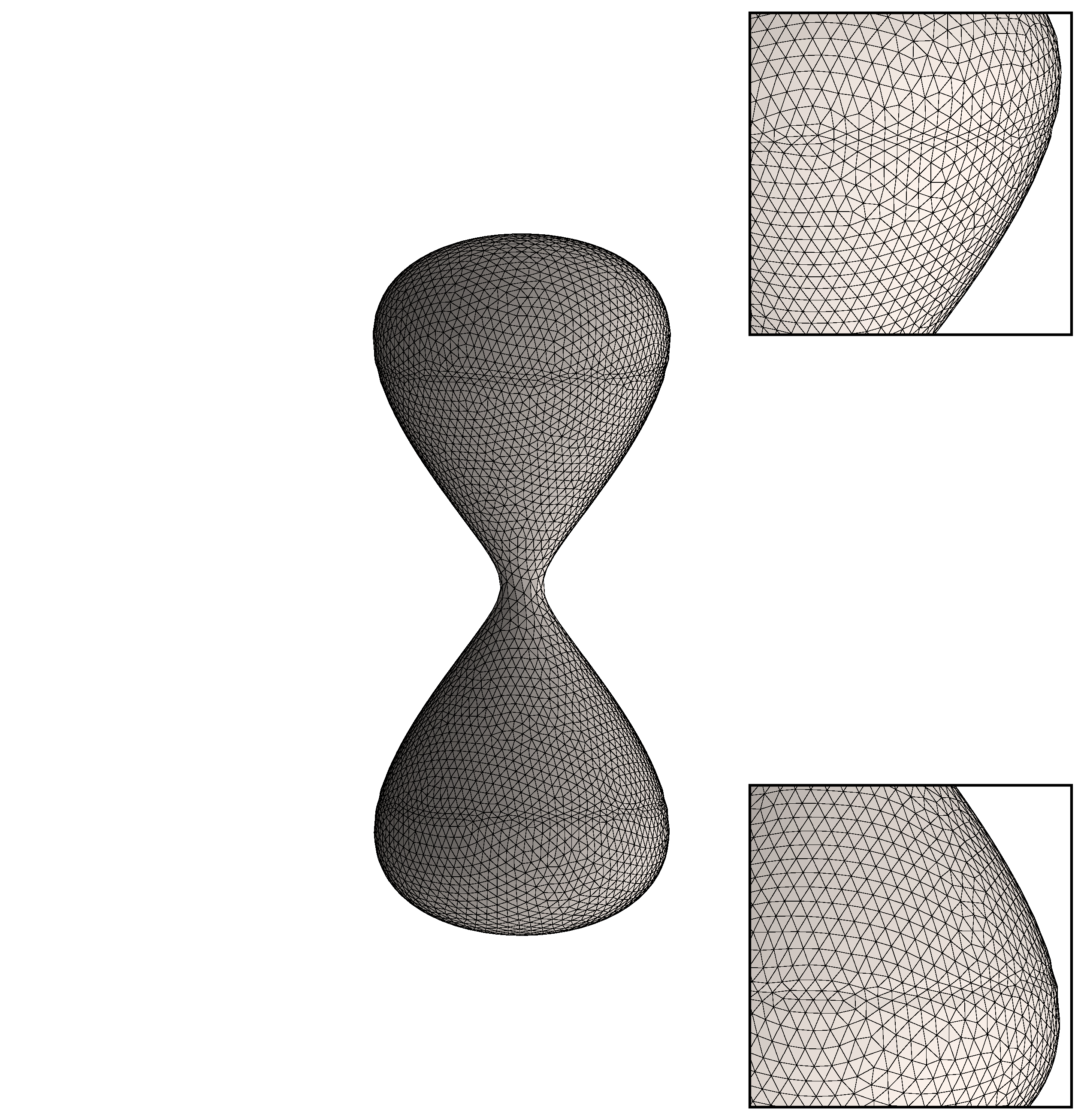}}
    \subfigure[$\eta=10$]{\includegraphics[width=3.2cm]{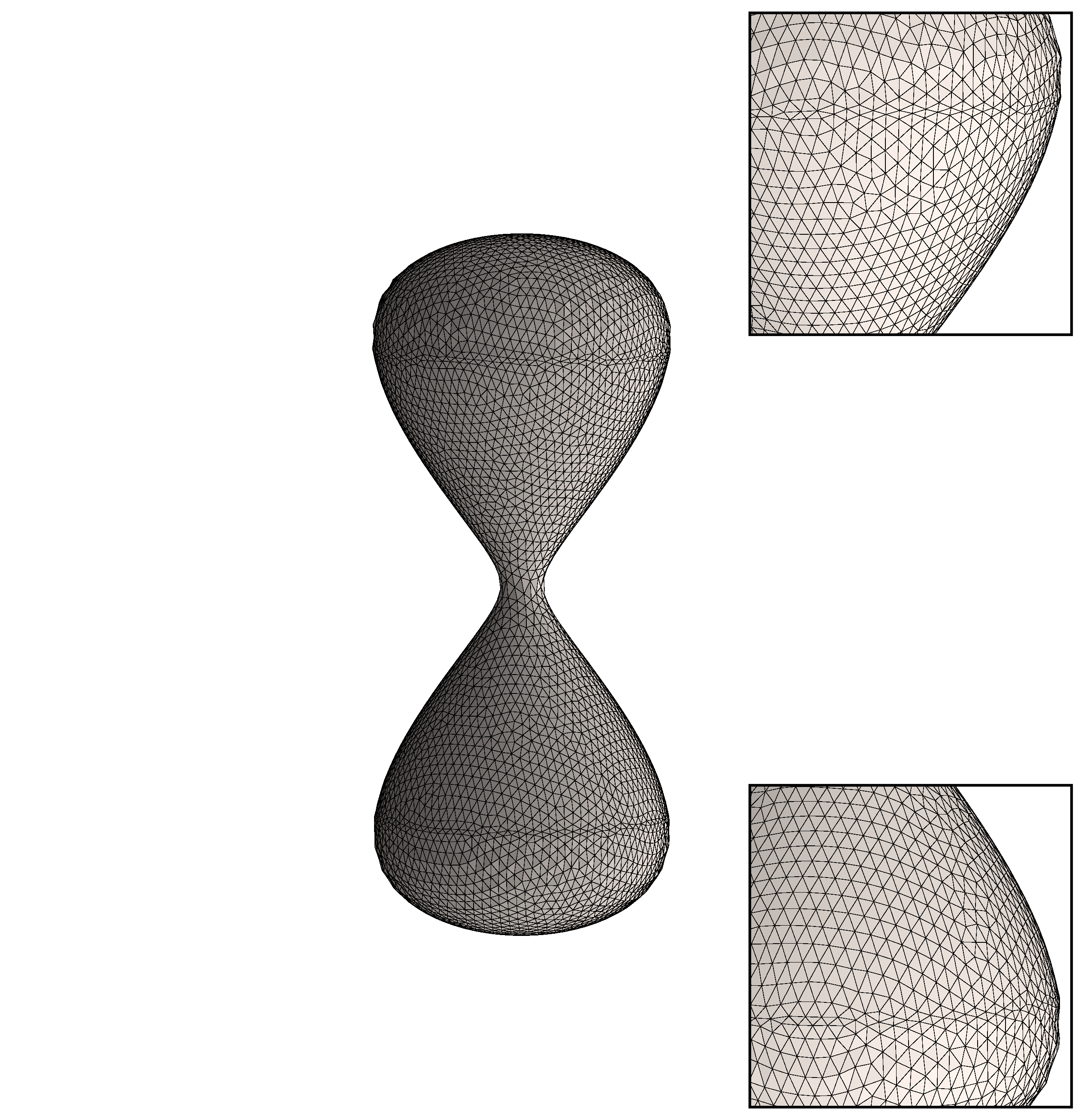}}
    \subfigure[$\eta=100$]{\includegraphics[width=3.2cm]{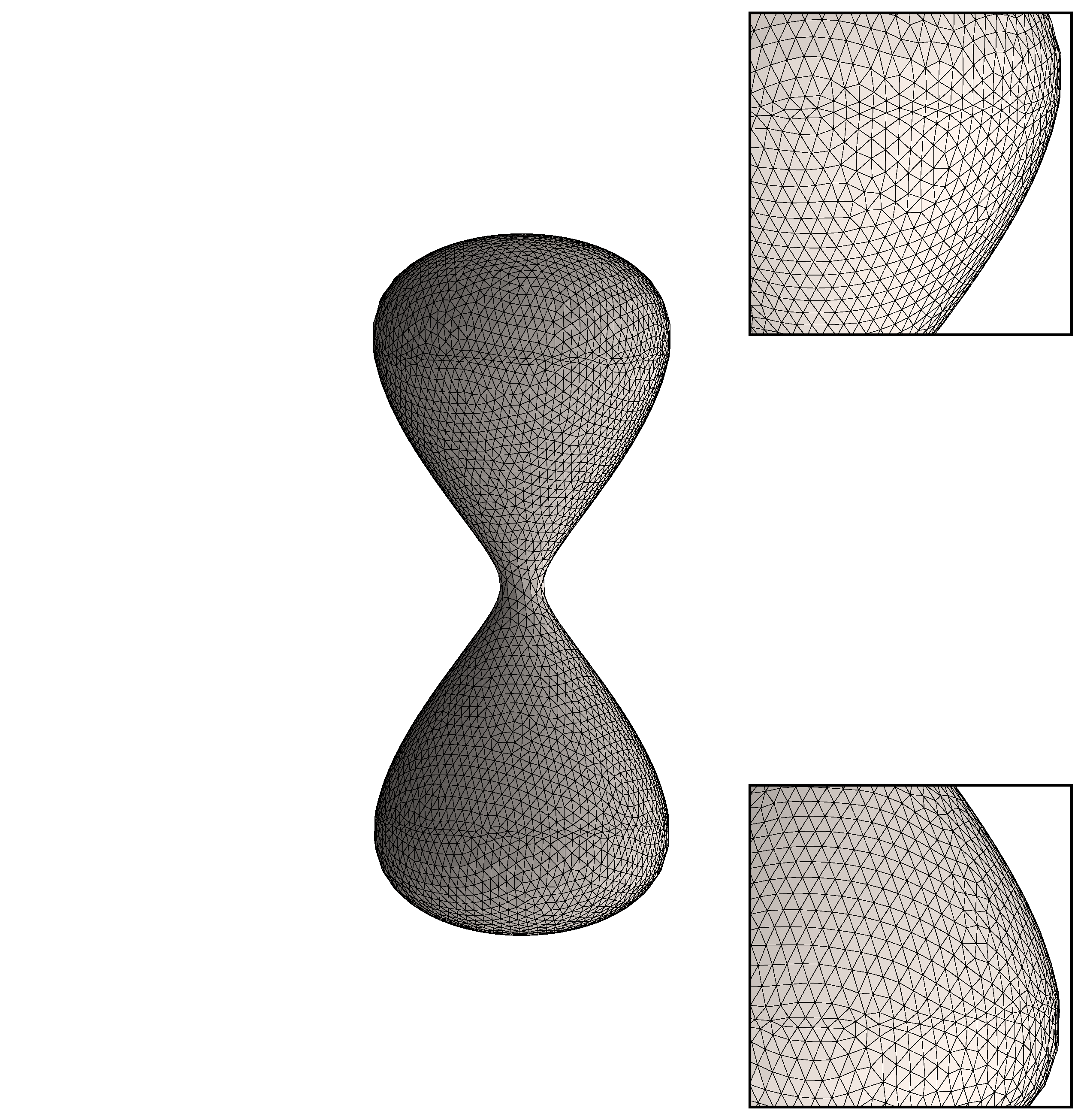}}\\
    \rotatebox{90}{~~~~After re-distribution}
    \subfigure[$\eta=0$]{\includegraphics[width=3.2cm]{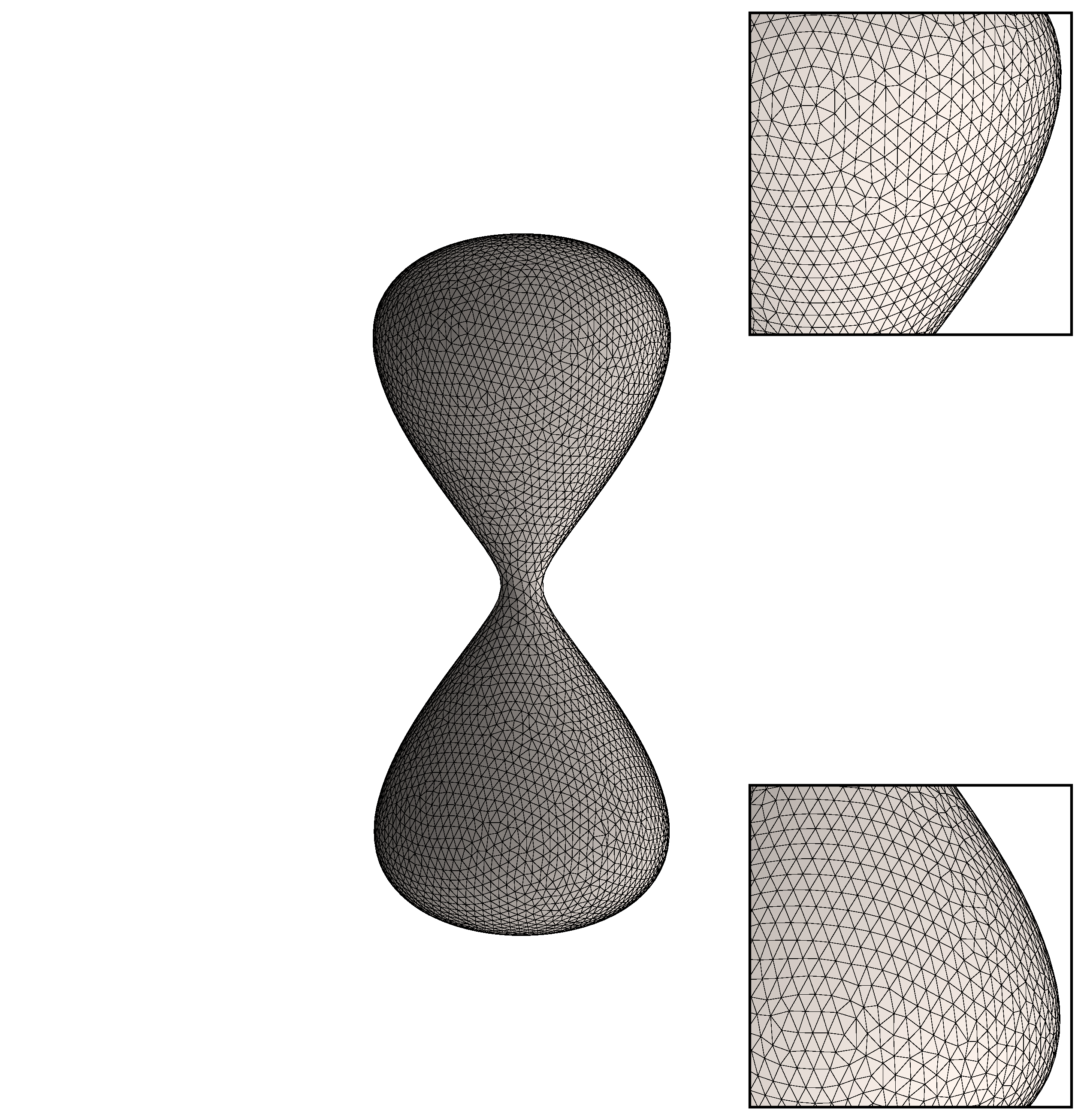}}
    \subfigure[$\eta=1$]{\includegraphics[width=3.2cm]{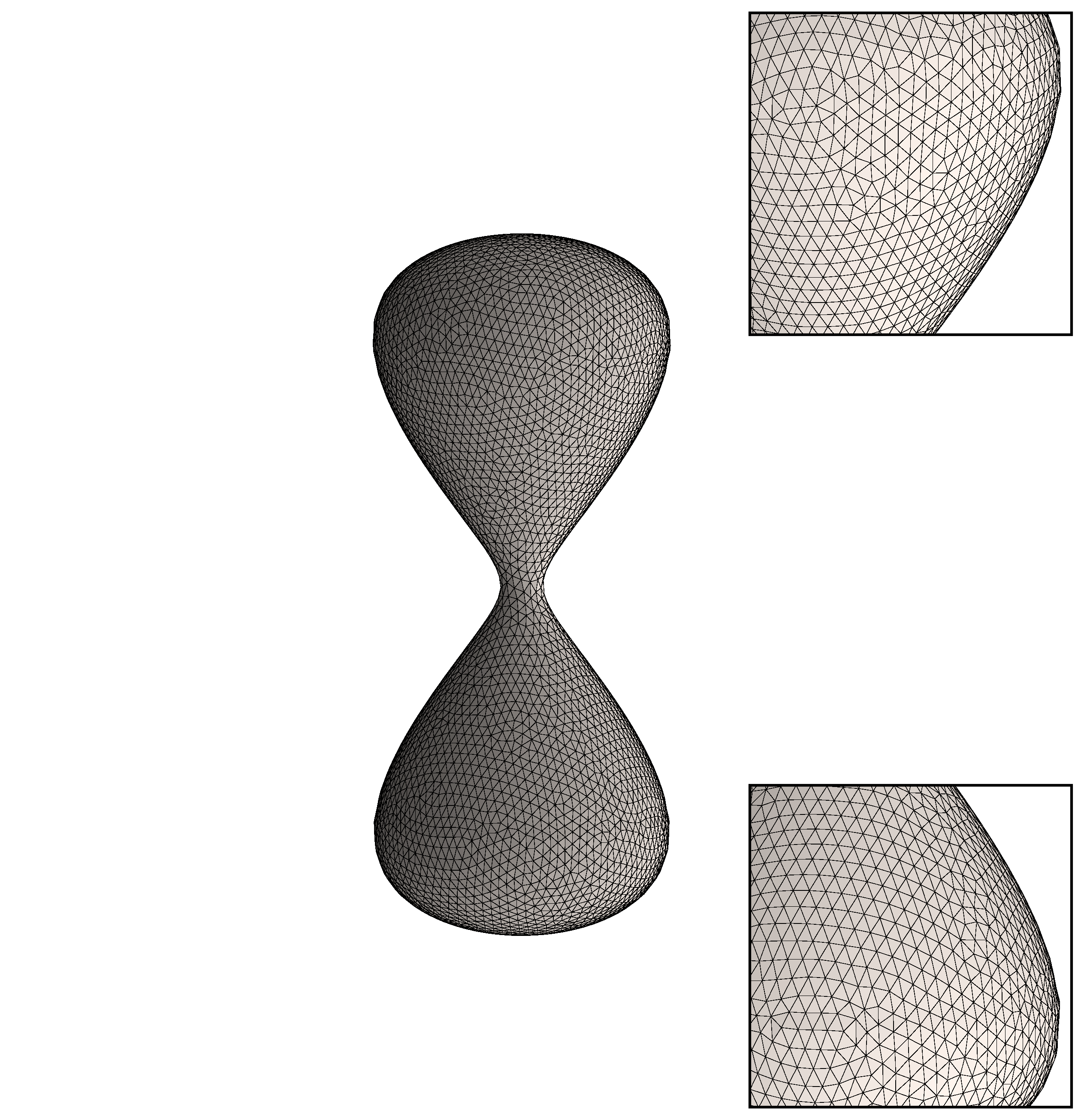}}
    \subfigure[$\eta=10$]{\includegraphics[width=3.2cm]{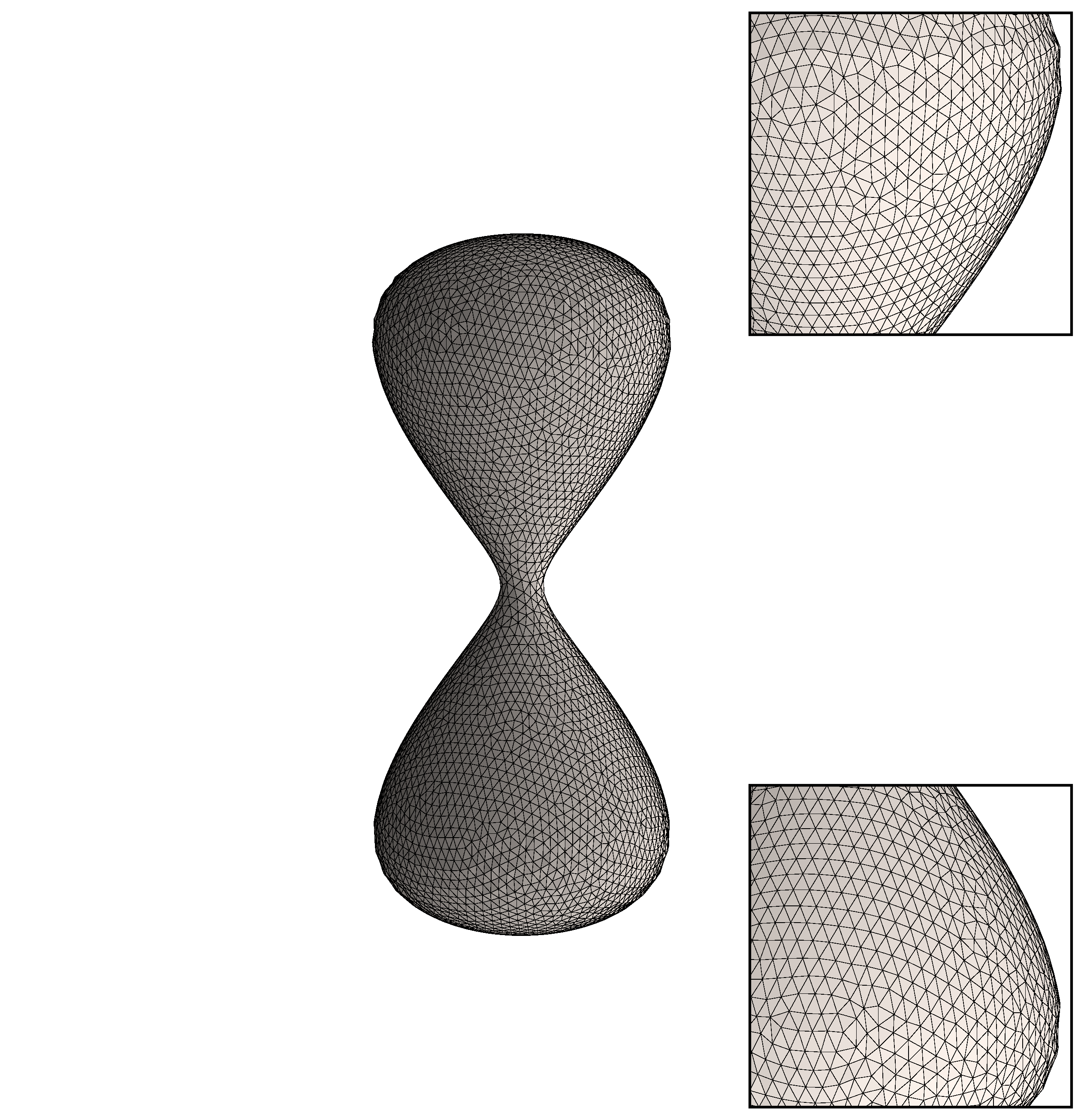}}
    \subfigure[$\eta=100$]{\includegraphics[width=3.2cm]{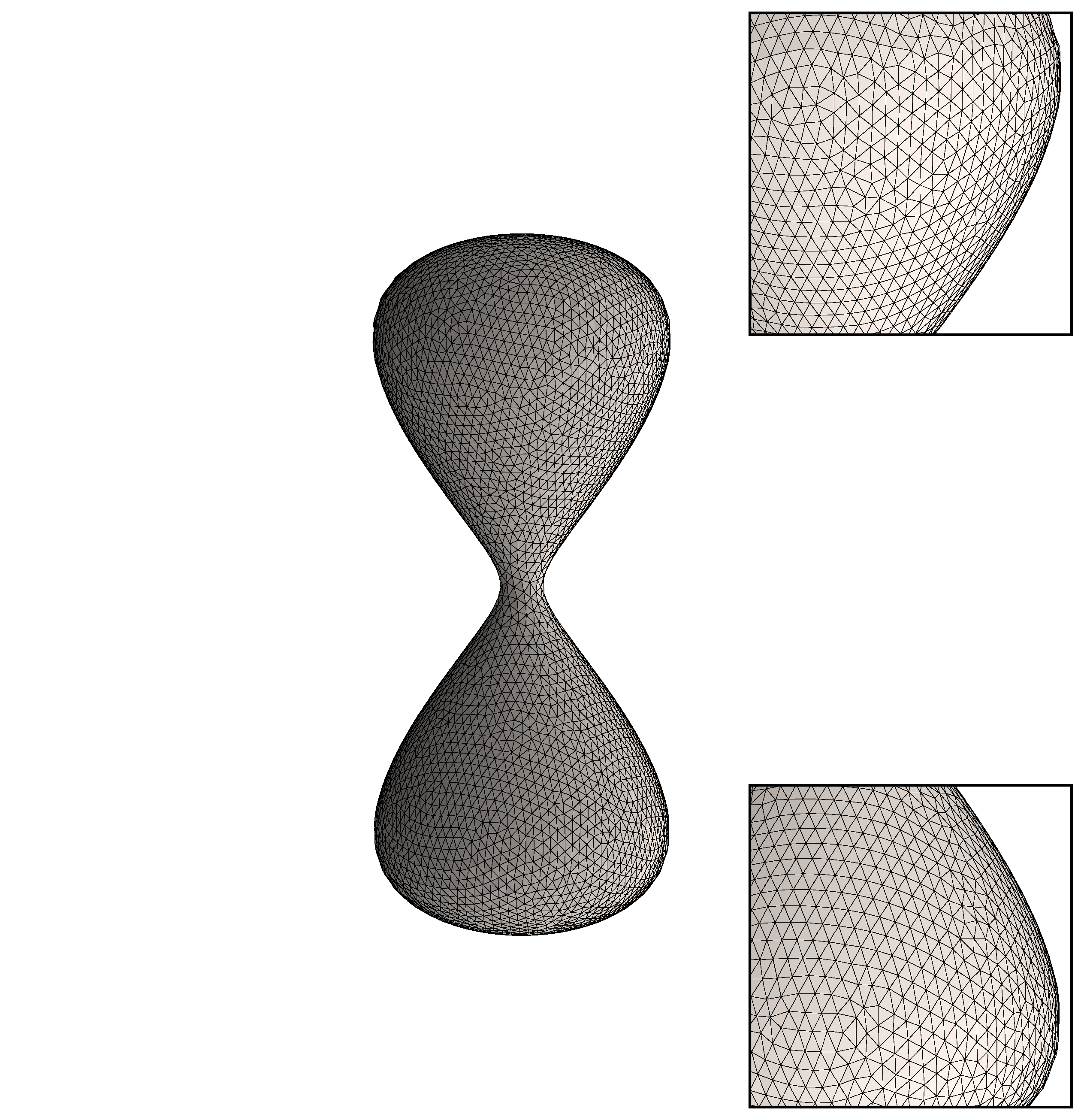}}\\
    \caption{Surfaces at $T=0.6$ computed by coupled model (Algorithm~\ref{BDF2-Algorithm}) and re-distribution method (Algorithm~\ref{RA-Algorithm}) with different $\eta$.}
    \label{imp-dumbbell}
    \end{center}
\end{figure}

{\bf Dumbbell-shaped surface with complex velocity.} 
We consider an evolving surface $\Gamma(t)=\{\boldsymbol{x}=(x,y,z)\in\mathbb{R}^{3}: \varphi(\boldsymbol{x}, t)=1\}$  described by a level set function
\begin{equation*}\label{imp-levelset-dumbbell}
    \begin{aligned}
        \varphi(\boldsymbol{x}, t) = \frac{x^2}{a^2(t)} + \frac{y^2}{a^2(t)} + G(\frac{z^2}{b^2(t)}),
    \end{aligned}
\end{equation*}
where $G(s) = 200s(s-\frac{199}{200})$, $a(t) = 0.1+0.05\sin{2\pi t}$ and $b(t) = 1+0.2\sin{4\pi t}$. 
The surface evolves under the following velocity field:
\begin{equation*}\label{imp-v-dumbbell}
    \begin{aligned}
        \boldsymbol{v}(\boldsymbol{X}(\boldsymbol{x}, t),t) = -\frac{\varphi_{t}\nabla\varphi}{|\nabla\varphi|^2}.
    \end{aligned}
\end{equation*}
This example had been considered in \cite{elliott2012ale} for surface finite element methods, where tangential velocity method was used to improve mesh quality. 
When the evolution time $T \geq  0.6$, the vertices are clustering and aligned along narrow bands, leading to anisotropic point distribution. 
Such degeneration is a bottleneck for long-term simulations and compromises numerical stability. 
This phenomenon has been observed in our point cloud algorithm (e.g. Algorithm~\ref{BDF2-Algorithm}) as well, as
clearly illustrated in Figure~\ref{imp-dumbbell}, where the point accumulation is visually evident. 
We use BDF2 scheme with $N_{\boldsymbol{x}}=5606$ and $\tau=1\times10^{-4}$ for surface evolution. 

\begin{figure}
    \centering
    \subfigure[Distribution of area elements with different $\eta$]{\includegraphics[width=0.45\linewidth]{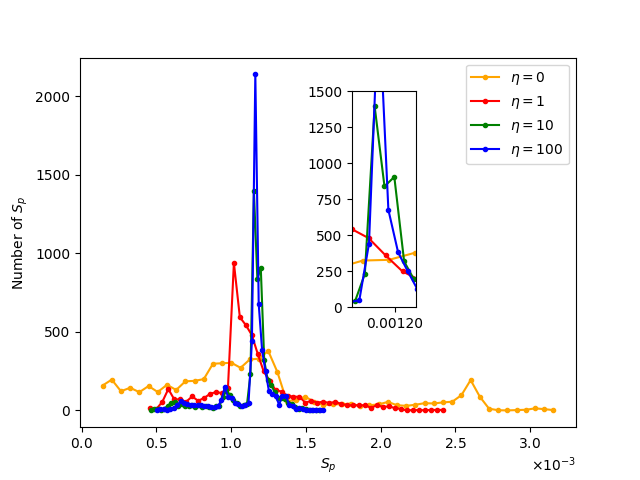}}
    \subfigure[Distribution of area elements after point redistribution]{\includegraphics[width=0.45\linewidth]{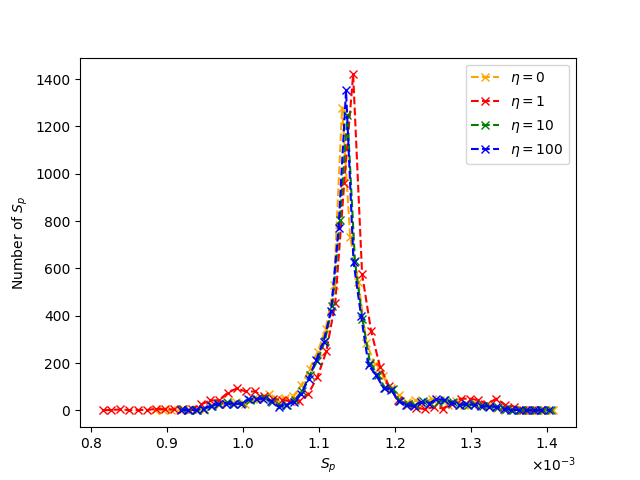}}
    \caption{Distribution of area element of each point, \eqref{eq:Sp}, for Algorithm~\ref{BDF2-Algorithm} and Algorithm~\ref{RA-Algorithm}, with different $\eta$. Surface is evolved up to time $T = 0.6$.}
    \label{imp-dumbbell-s}
\end{figure}
Similarly, in the first line of Figure~\ref{imp-dumbbell}, we present the results obtained using Algorithm~\ref{BDF2-Algorithm}, while the second line shows the improved results after applying Algorithm~\ref{RA-Algorithm}. 
In this example, it is evident that due to the complexity of the prescribed velocity field, by adjusting the scalar parameter $\eta$ alone appears insufficient to fully prevent point aggregation. 
Though, it does play an important role in mitigating the clustering effect. 
As expected, the application of Algorithm~\ref{RA-Algorithm} successfully achieves the desired spatial redistribution. 
In addition, similar results are also presented in Figure~\ref{imp-dumbbell-s}, where we see that Algorithm~\ref{RA-Algorithm} makes the distribution of points more uniform.

{\bf Nonuniform distribution on surfaces.} 
In practical applications, there might be a need for a non-uniform distribution of points, with larger density in regions of particular interest. For instance, the areas exhibit high curvature on the surface. To accommodate such requirements, we use Algorithm~\ref{Target-RA-Algorithm} to have targeted point redistribution. Specifically, we introduce a prescribed target distribution $p(\boldsymbol{x})$ to guide the evolution of the density function. For simplicity, we define the target distributions $p(\boldsymbol{x})$ for both the ellipsoid surface and the dumbbell-shaped surface as follows:
\begin{equation*}\label{target-p}
    \begin{aligned}
        p_{e}(\boldsymbol{x}) &\propto
        \exp{(\theta(\sin{z}+1))}, \\
        p_{d}(\boldsymbol{x}) &\propto \exp{(\theta(\cos{z}+1))}.
    \end{aligned}
\end{equation*}
\begin{figure}[htbp]
    \begin{center}
    \subfigure[$\theta=0$]{\includegraphics[width=4cm]{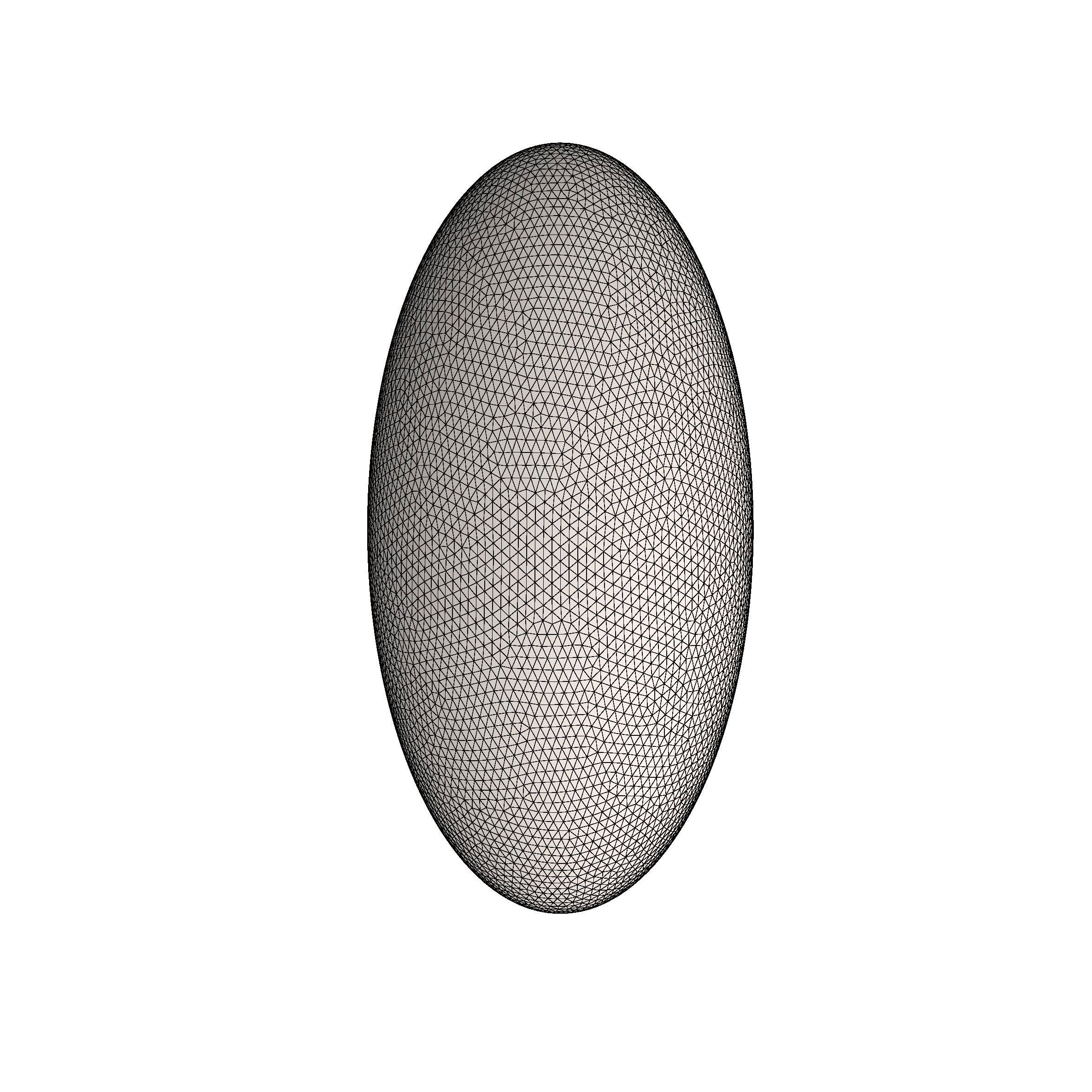}}
    \subfigure[$\theta=2$]{\includegraphics[width=4cm]{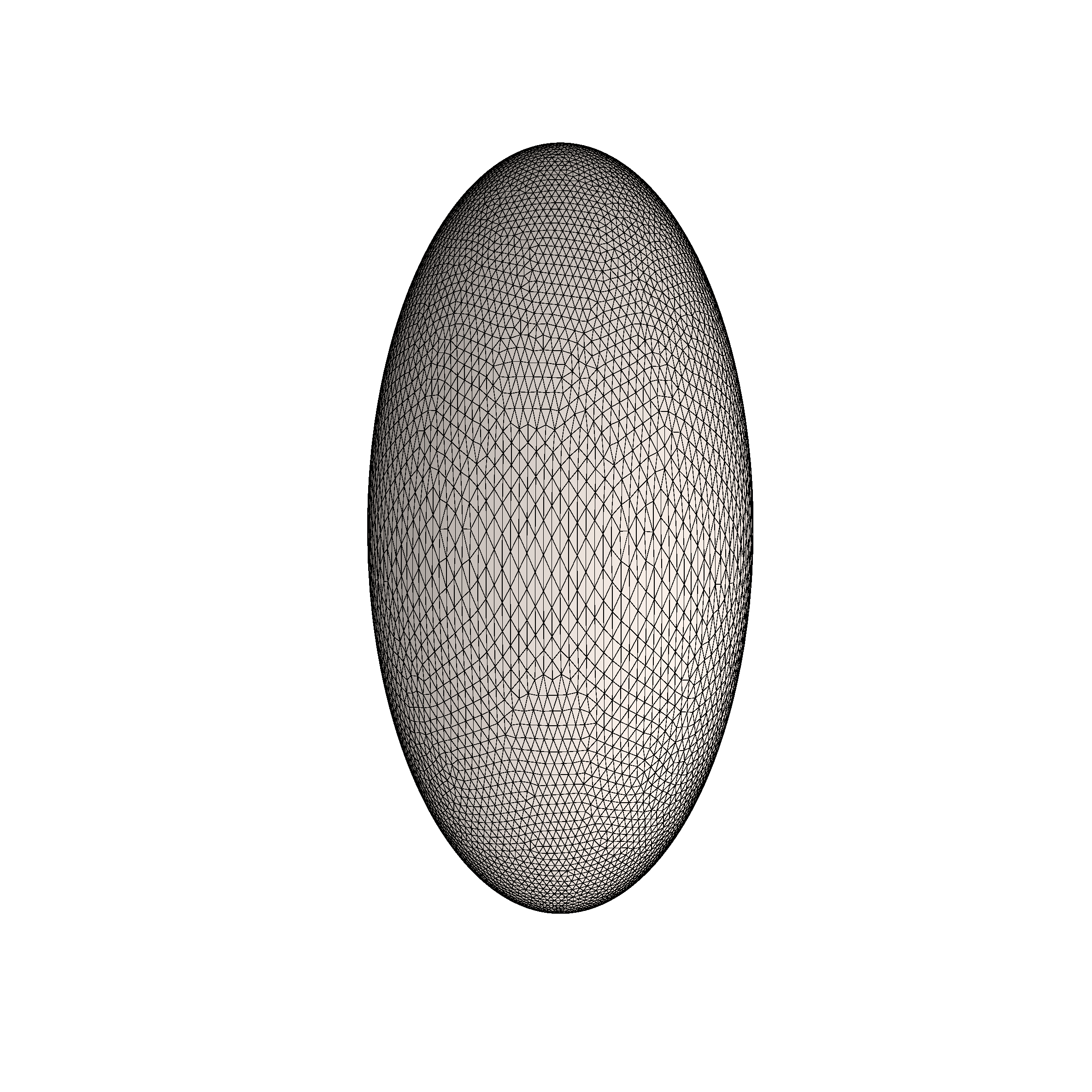}}
    \subfigure[$\theta=4$]{\includegraphics[width=4cm]{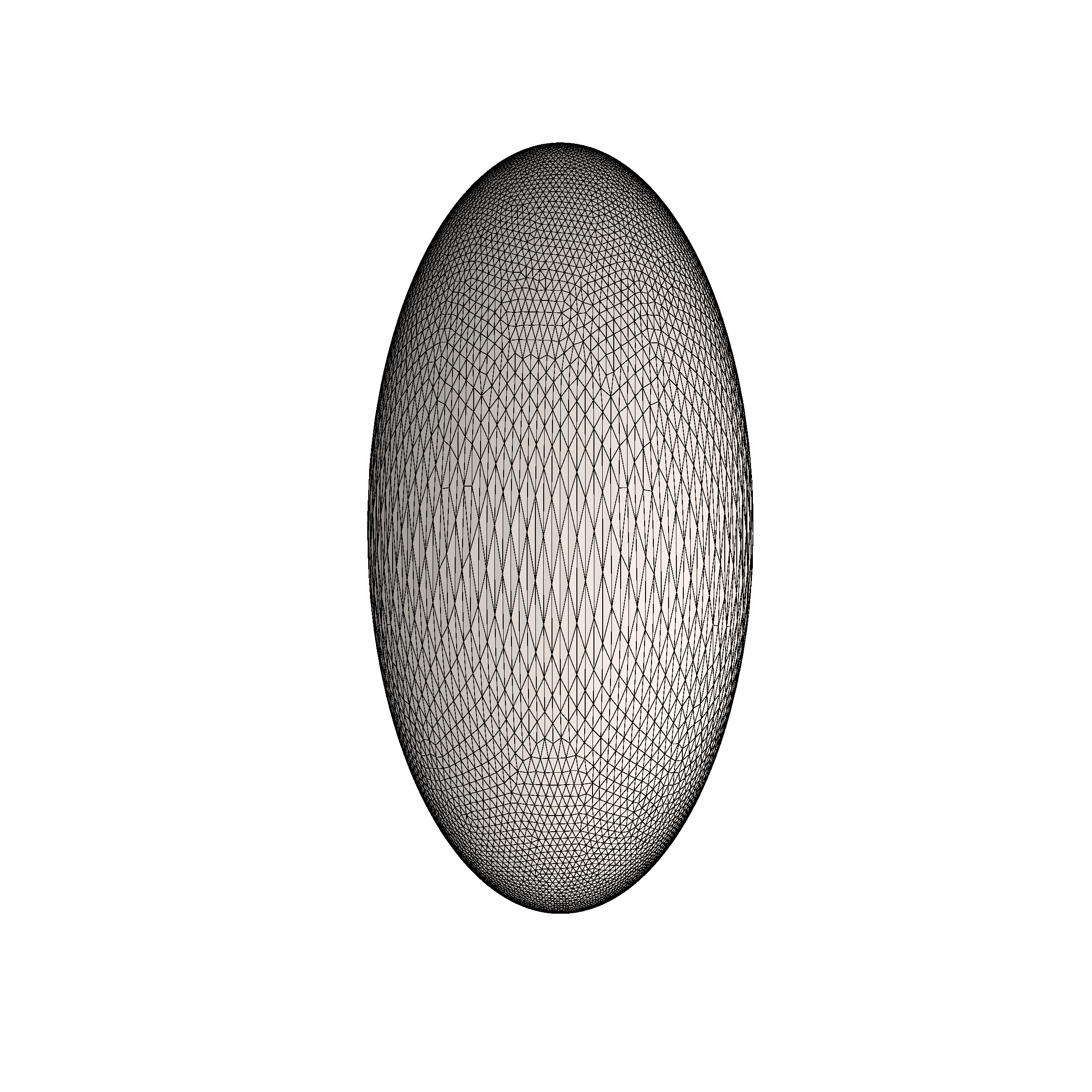}}
    
    \subfigure[$\theta=0$]{\includegraphics[width=4cm]{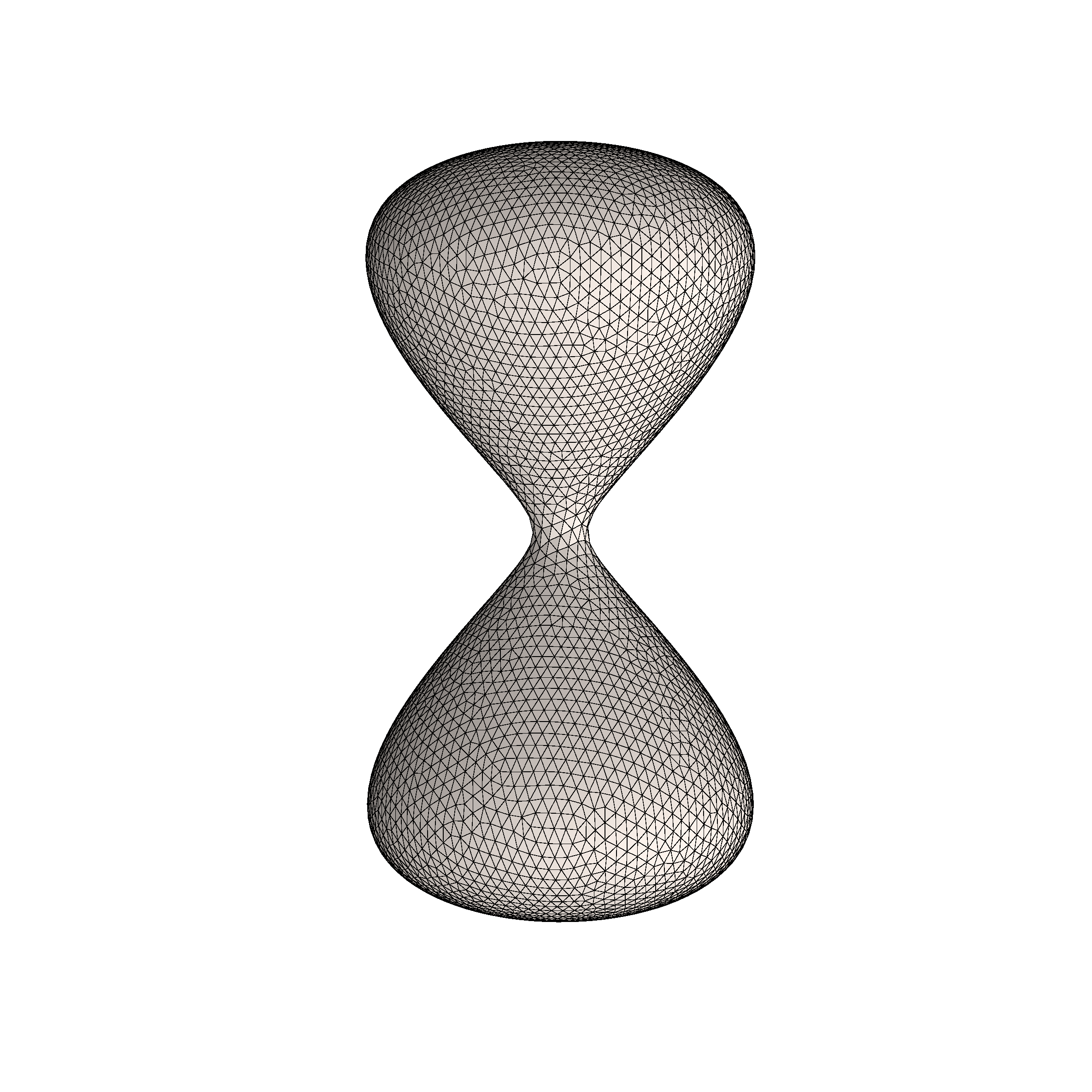}}
    \subfigure[$\theta=2$]{\includegraphics[width=4cm]{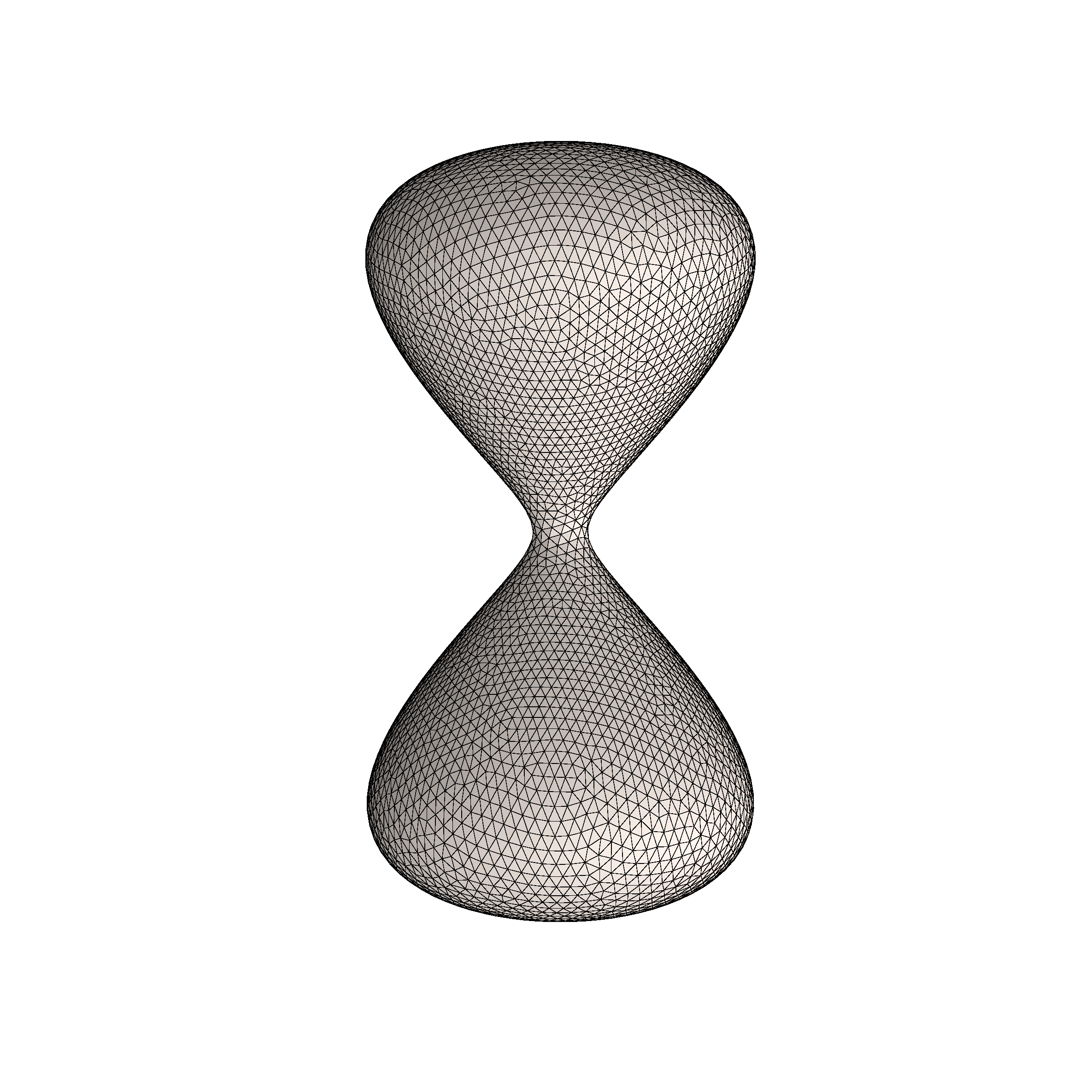}}
    \subfigure[$\theta=4$]{\includegraphics[width=4cm]{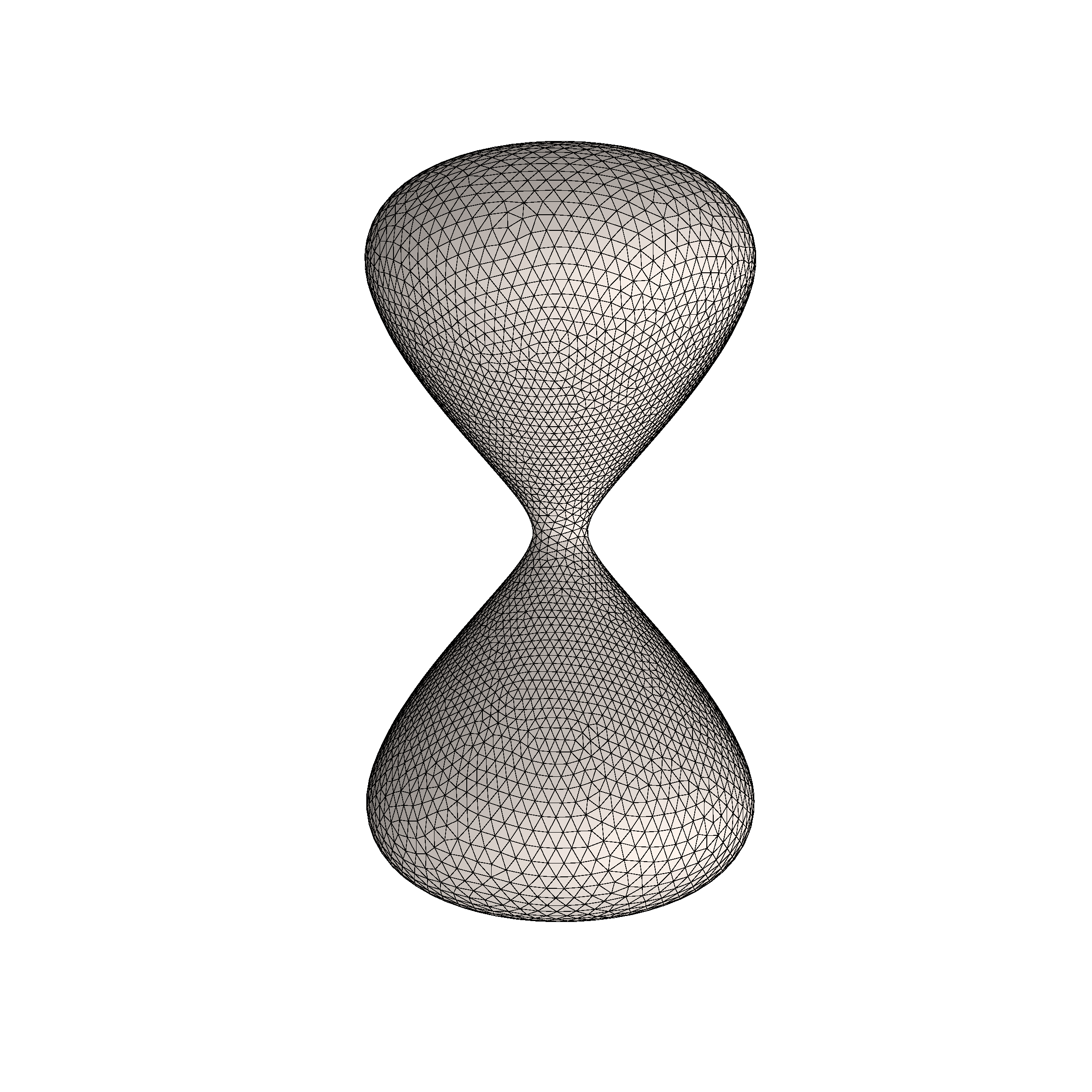}}
    \caption{Results of point redistribution with assigned target distribution $p(\boldsymbol{x})$ (Algorithm \ref{Target-RA-Algorithm}).}
    \label{diff-p}
    \end{center}
\end{figure}

From Figure~\ref{diff-p}, it is evident that Algorithm~\ref{Target-RA-Algorithm} effectively redistributes points in accordance with the prescribed target distributions. On the ellipsoidal surface, the points are concentrated toward both polar points, aligning with regions of higher target density. In contrast, for the dumbbell-shaped surface, the points accumulate predominantly in the central region. Both examples demonstrate the algorithm’s capability of adaptively concentrating points as it is designed for. These results confirm that Algorithm~\ref{Target-RA-Algorithm} can flexibly accommodate spatially varying distribution requirements in complex geometries.

\subsection{Mean curvature flow}
In this section, we apply our method to the renowned mean curvature flow (MCF) \cite{mantegazza2011lecture}. 
The MCF describes the geometric motion along the normal field scaled with the scaler mean curvature at each point of a surface as it evolves over time. 
This means the velocity satisfies the following equation:
\begin{equation*}\label{v-MCF}
    \begin{aligned}
        \boldsymbol{v}(\boldsymbol{X}(\boldsymbol{x}, t),t) = -H \boldsymbol{n},
    \end{aligned}
\end{equation*}
where $H$ represents the scalar mean curvature ($H\boldsymbol{n}=-\Delta_{\Gamma}\boldsymbol{X}$) and $\boldsymbol{n}$ is the outward unit normal vector field. 
Then the problem \eqref{ES-rho-eq} can be changed to
\begin{equation*}\label{ES-rho-eq-MFC}
    \begin{aligned}
        \frac{d }{d t}\boldsymbol{X}(\boldsymbol{x}, t) &= \boldsymbol{v}(\boldsymbol{X}(\boldsymbol{x}, t), t) -\eta\nabla_{\Gamma(t)}s(\boldsymbol{X}(\boldsymbol{x}, t), t), && \boldsymbol{x}\in \Gamma_{0}, \quad t\geq 0,\\
        \boldsymbol{v}(\boldsymbol{X}(\boldsymbol{x}, t), t) &= \Delta_{\Gamma(t)}\boldsymbol{X}(\boldsymbol{x}, t), && \boldsymbol{x}\in \Gamma_{0}, \quad t\geq 0,\\
        \frac{d }{d t} s(\boldsymbol{X}(\boldsymbol{x}, t), t) &= \eta\Delta_{\Gamma(t)}s(\boldsymbol{X}(\boldsymbol{x}, t), t) - \text{div}_{\Gamma(t)}\boldsymbol{v}(\boldsymbol{X}(\boldsymbol{x}, t), t), && \boldsymbol{x}\in \Gamma_{0}, \quad t\geq 0,\\
        \boldsymbol{X}(\boldsymbol{x}, 0) &= \boldsymbol{x}, \quad s(\boldsymbol{x}, 0) = \log\rho_{0}(\boldsymbol{x}), && \boldsymbol{x}\in \Gamma_{0}. 
    \end{aligned}
\end{equation*}
Using the discrete and notation in Section \ref{Sec-Discret}, we obtain the following fully discrete scheme 
\begin{equation}\label{fullDis-ES-rho-eq-MFC}
    \begin{aligned}
        \frac{a}{\tau}\boldsymbol{X}_{k} - \boldsymbol{v}_{k} + \eta M_{G}^{T}(\bar{\boldsymbol{X}}_{k})s_{k} &= \frac{1}{\tau}\hat{\boldsymbol{X}}_{k-1}, \\
        M_{L}(\bar{\boldsymbol{X}}_{k})\boldsymbol{X}_{k} - \boldsymbol{v}_{k} &= 0, \\
        M_{G}(\bar{\boldsymbol{X}}_{k})\boldsymbol{v}_{k} + \frac{a}{\tau}s_{k} - \eta M_{L}(\bar{\boldsymbol{X}}_{k})s_{k} &= \frac{1}{\tau}\hat{s}_{k-1}, \\
        \boldsymbol{X}(\boldsymbol{x}, 0) = \boldsymbol{x}, \quad s(\boldsymbol{x}, 0) &= \log\rho_{0}(\boldsymbol{x}).
    \end{aligned}
\end{equation}
Then, by solving \eqref{fullDis-ES-rho-eq-MFC}, we obtain the discrete trajectory for the MCF.

We begin with a specific example. 
The parameterization of the two-dimensional surface under consideration is as follows:
\begin{equation}\label{levelset-MCF}
    \begin{aligned}
        \boldsymbol{x}(\theta, \phi) = \left(
        \begin{array}{c}
             \cos{\phi}  \\
             (\frac{3}{5}\cos^2\phi + \frac{2}{5})\cos{\theta}\sin{\phi} \\ 
             (\frac{3}{5}\cos^2\phi + \frac{2}{5})\sin{\theta}\sin{\phi}
        \end{array}
        \right), \quad \theta\in[0, 2\pi), \quad\phi\in[0, \pi],
    \end{aligned}
\end{equation}
which is a benchmark example proposed in \cite{m2017approximations}.
It is well known that without incorporating tangential velocity, the approximation of this problem typically faces mesh distortion and node clustering, leading to numerical breakdown before the surface evolving into a sphere. This typically necessitates mesh redistribution techniques, as discussed in \cite{duan2021high}. 
Therefore, the using of extra tangential velocity is a cure for numerical treatment. 
Several mesh-based approaches have been proposed to enhance the mesh quality by introducing artificial tangential motion. Notably, the BGN algorithm in \cite{barrett2007parametric, barrett2008parametric, barrett2008parametric} and DeTurck’s reparameterization technique in \cite{m2017approximations} are those representatives. These methods rely on an appropriate initial parameterization, such as that in \eqref{levelset-MCF}, and are capable of evolving the surface smoothly into a sphere. 

In the following, we denote by $E_h$ the surface area of the numerically computed surface. 
We use the rescaling trick to maintain the area of $\boldsymbol{X}(t)$ in a reasonable range during the computation.
\begin{equation*}\label{scaling}
    \begin{aligned}
        \boldsymbol{X}^{*}=\lambda\boldsymbol{X}, \quad s^{*}=s, \quad t^{*}=\lambda^2t,
    \end{aligned}
\end{equation*}
where $\lambda$ is a scaling parameter, which is adjusted adaptively in the computation. 
We use the BDF2 scheme, $N_{\boldsymbol{x}}=5606$ and $\tau=1\times10^{-3}$ as initial value for the surface evolution. 
Moreover, during the simulation, point redistribution (Algorithm~\ref{RA-Algorithm}) is invoked every $100$ time steps. This periodic correction ensures that the point cloud remains well-distributed throughout the evolution process, effectively preventing excessive concentration in localized regions.
\begin{figure}[htbp]
    \begin{center}
    \subfigure[Initial surface]{\includegraphics[width=4.5cm]{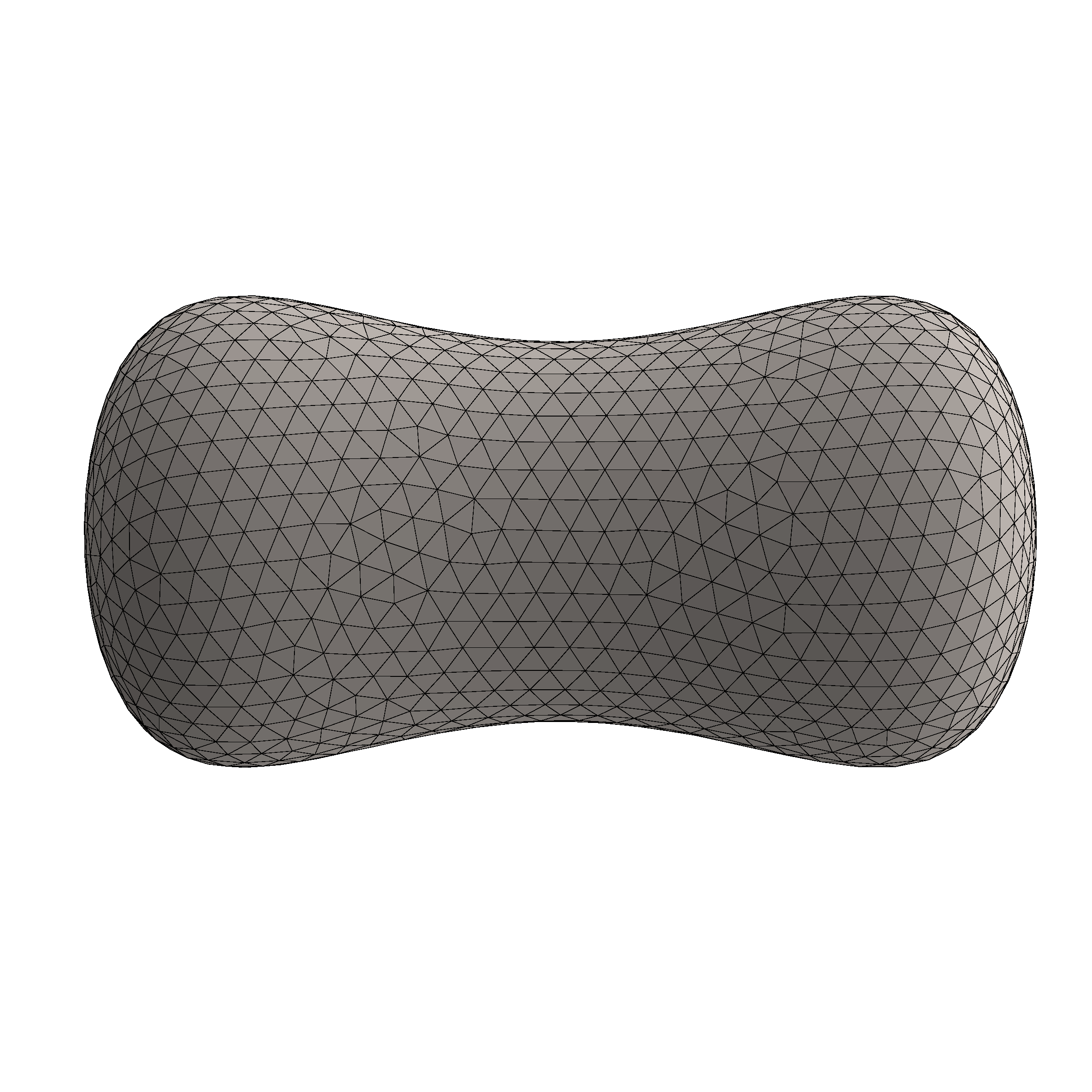}} \\
    \rotatebox{0}{~~~~~~~~$t=0.05$~~~~~~~~~~~~~~~~~~~~~~~~~~~~$t=0.087625$~~~~~~~~~~~~~~~~~$t=0.0896913948$}\\
    \rotatebox{90}{~~~~~~~~~~~~~~~~$\eta=0$}
    \subfigure[$E_{h}=3.91E-4$]{\includegraphics[width=4.5cm]{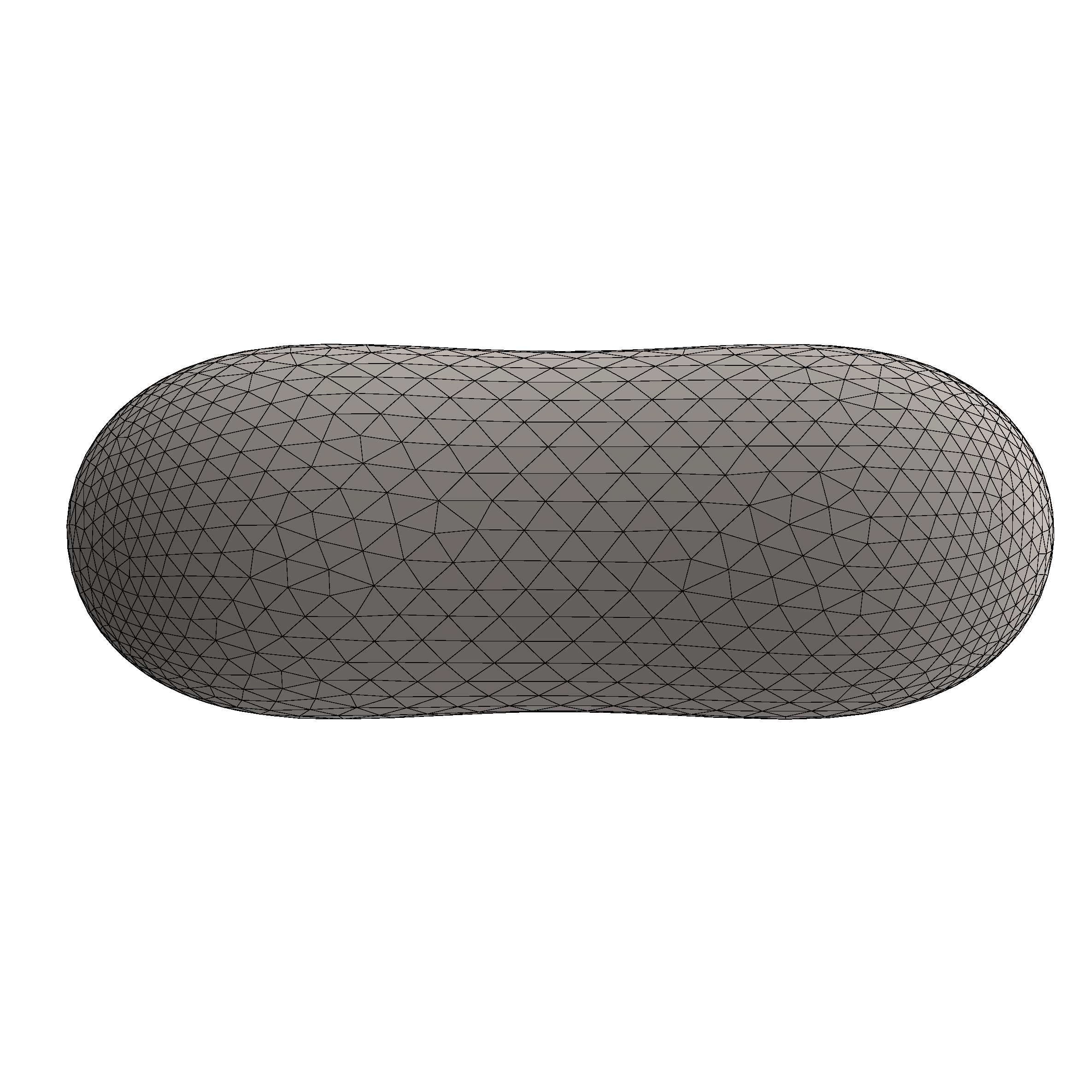}}
    \subfigure[$E_{h}=3.67E-10$]{\includegraphics[width=4.5cm]{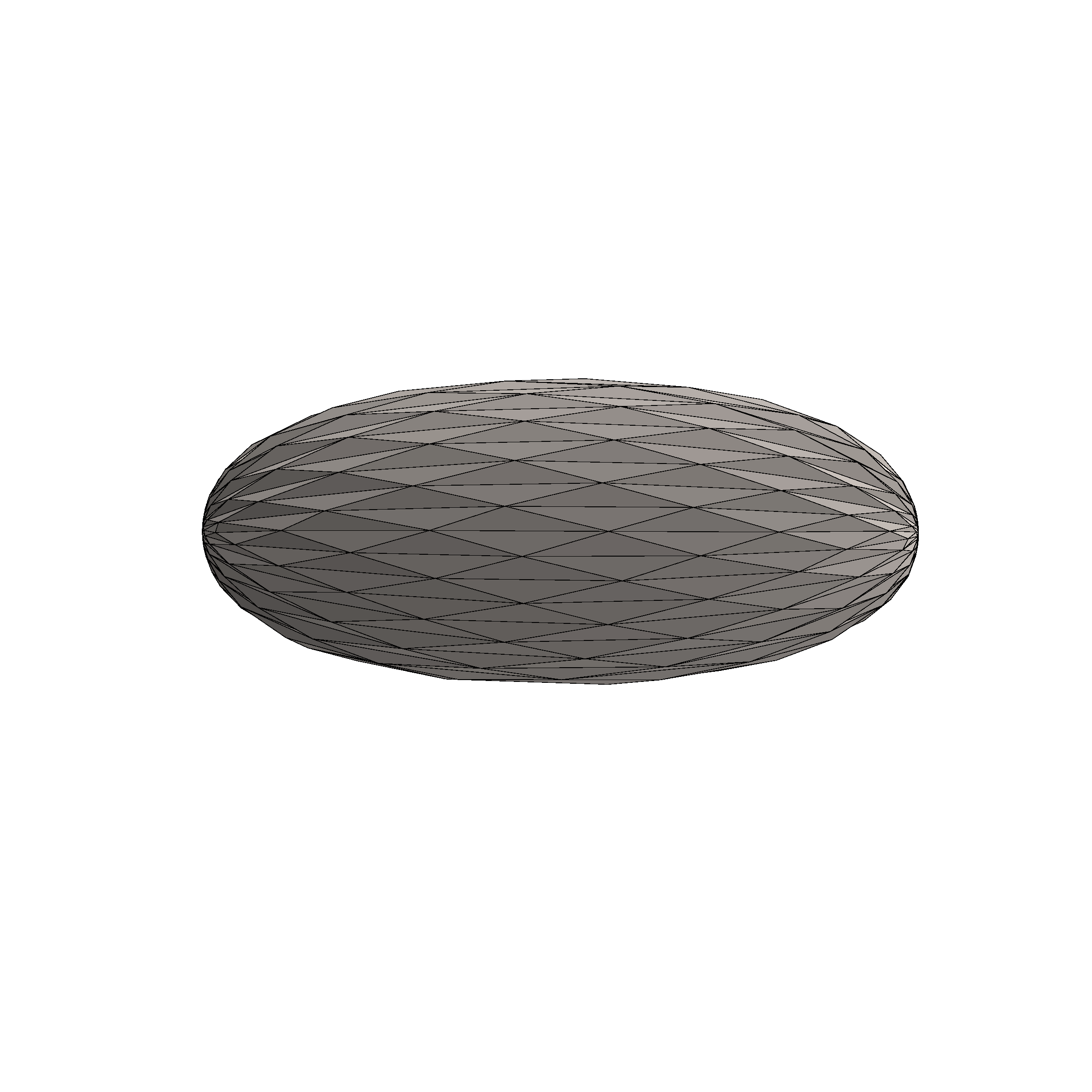}}
    \subfigure[none]{\includegraphics[width=4.5cm]{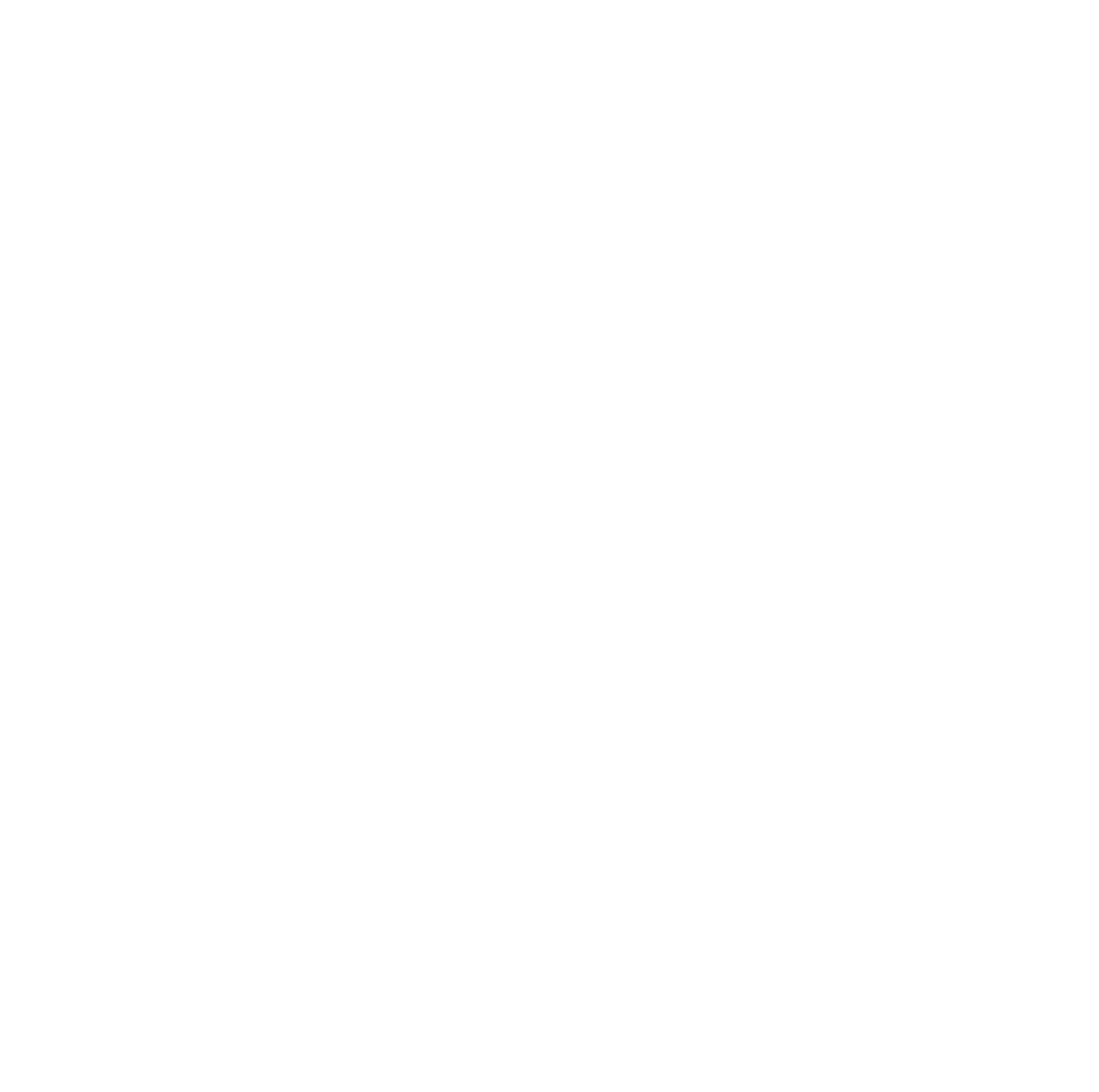}}\\
    \rotatebox{90}{~~~~~~~~~~~~~~$\eta=100$}
    \subfigure[$E_{h}=8.98E-4$]{\includegraphics[width=4.5cm]{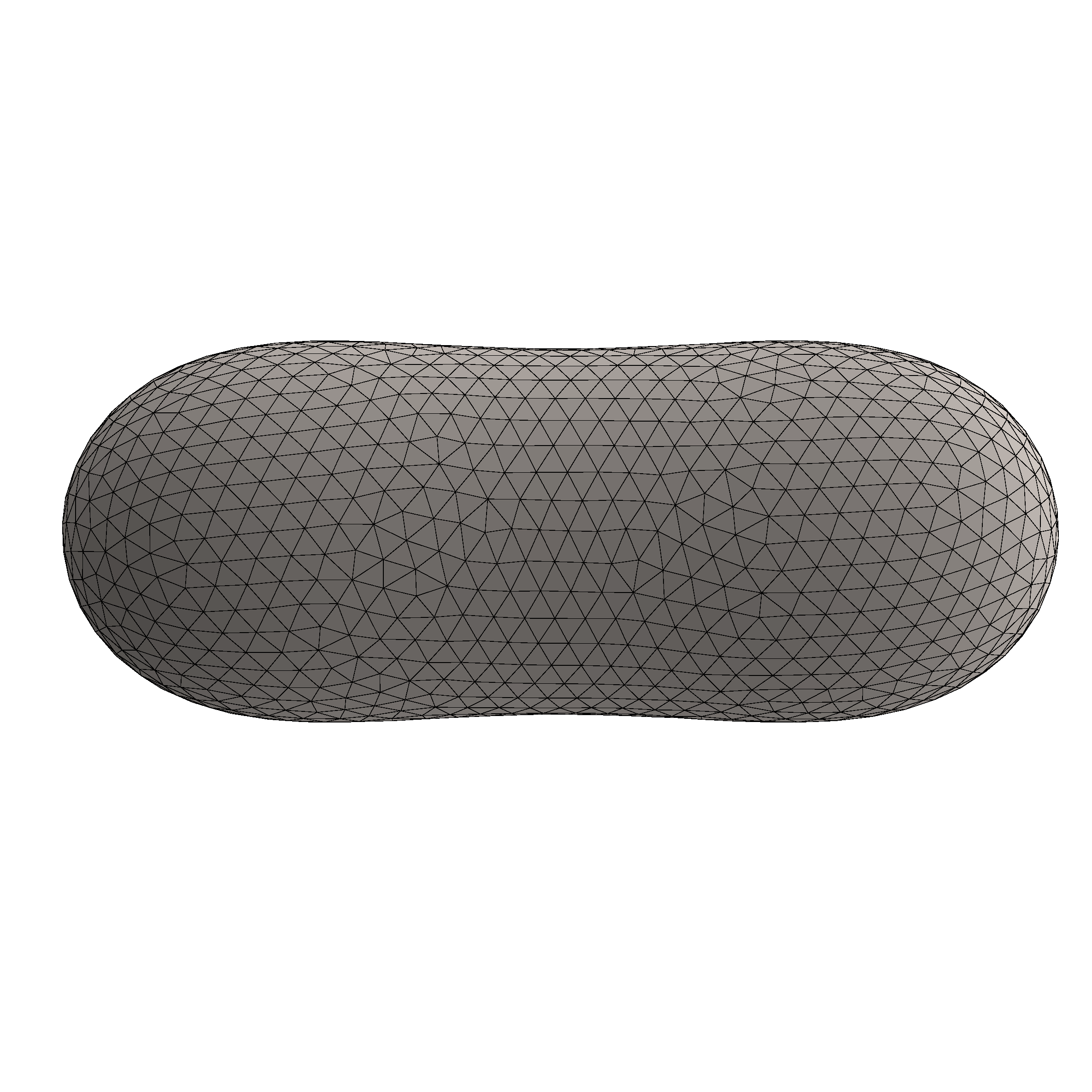}}
    \subfigure[$E_{h}=2.94E-5$]{\includegraphics[width=4.5cm]{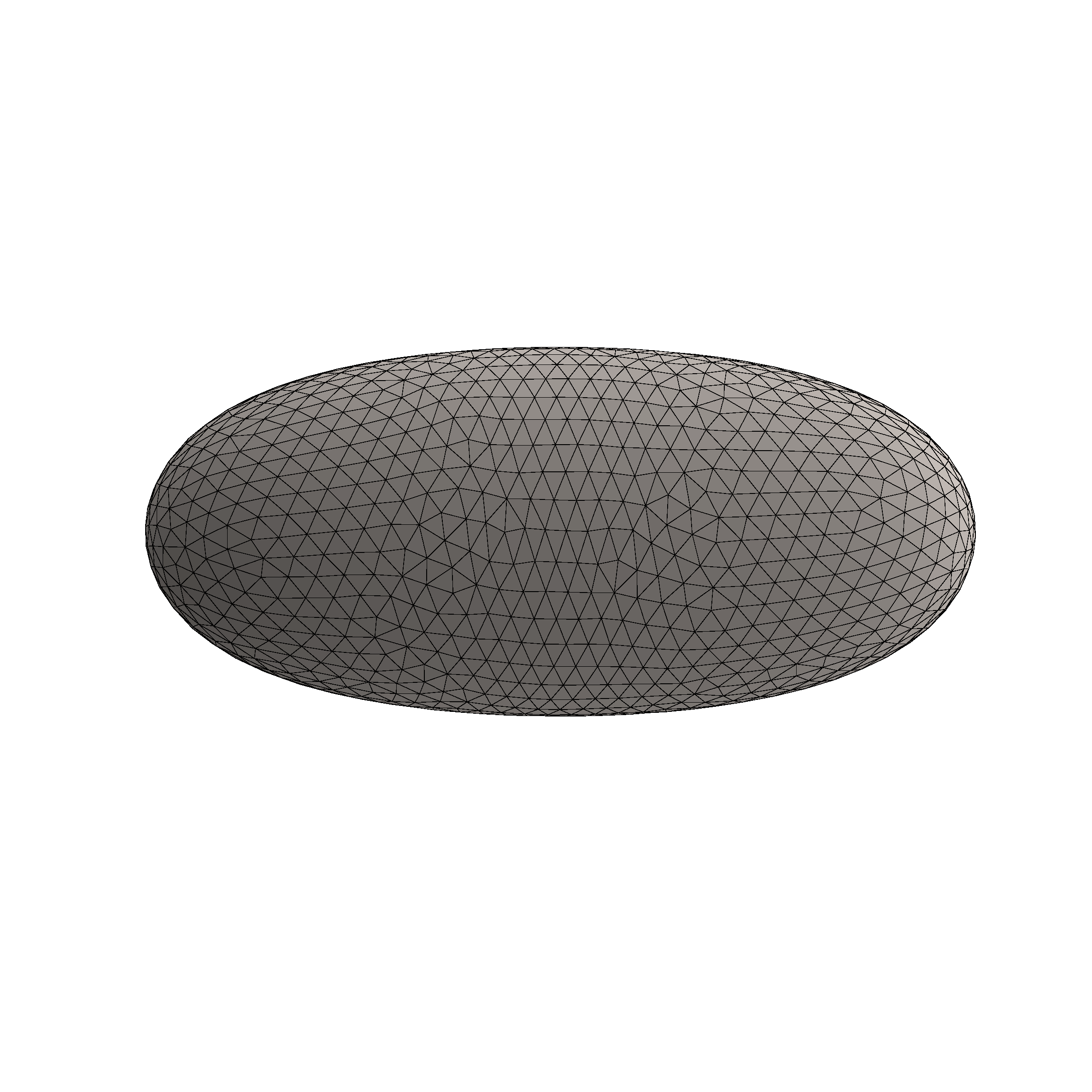}}
    \subfigure[$E_{h}=8.18E-39$]{\includegraphics[width=4.5cm]{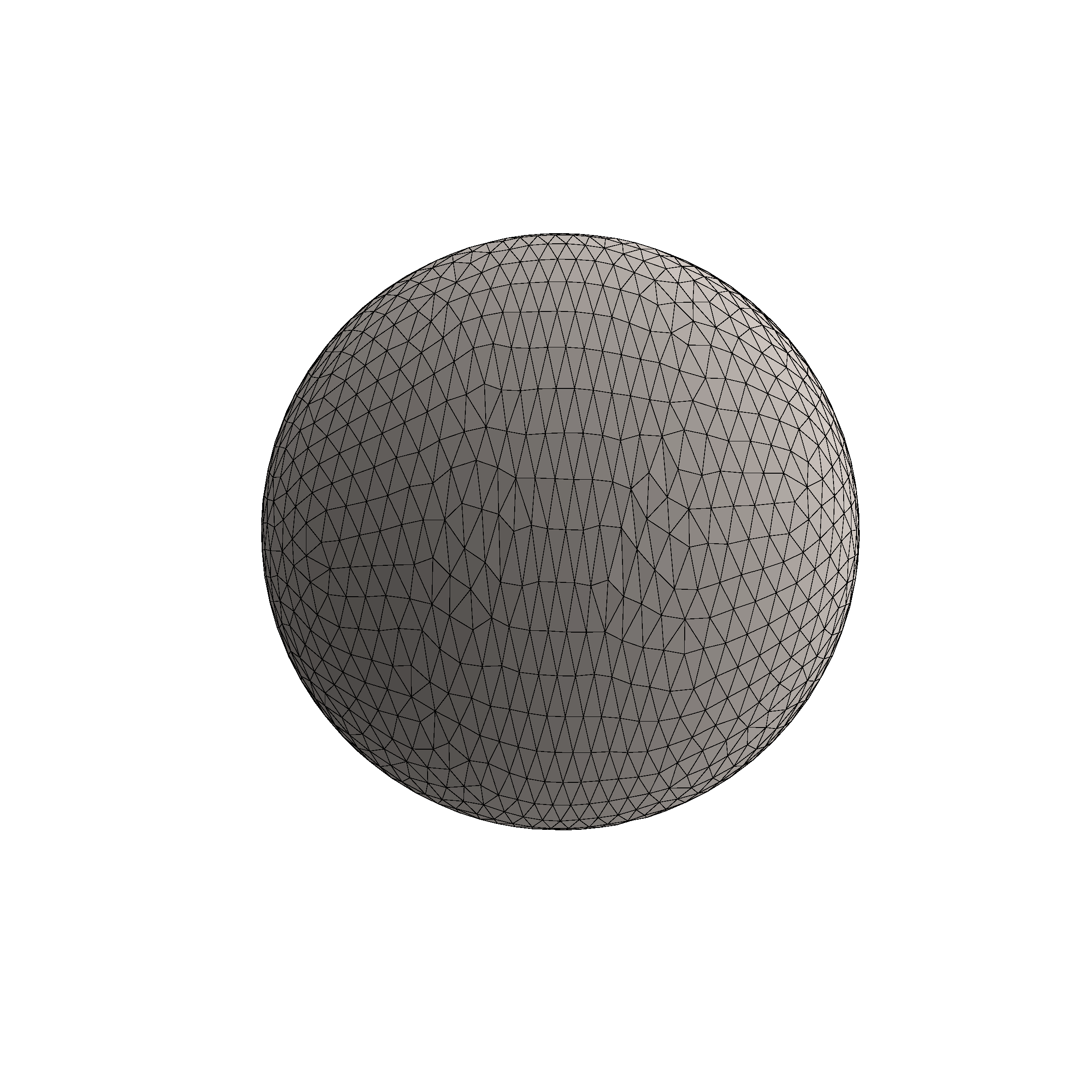}}\\
    \caption{Algorithm~\ref{RA-Algorithm} evolves the MCF problem (images are rescaled). (d):Without adding the artificial tangential velocity, the numerical method becomes unfeasible gradually.}
    \label{EMSL-DAC-MCF}
    \end{center}
\end{figure}

In the first line of Figure~\ref{EMSL-DAC-MCF}, we present the results of surface evolution without applying tangential velocity correction. 
It is evident that at time $t=0.087625$, the point cloud has accumulated significantly at both ends of the ellipsoid. 
At this stage, further evolution becomes numerically infeasible due to the severe point clustering. 
In contrast, the second line shows the results computed using Algorithm~\ref{RA-Algorithm}. 
At the same time step (Figure~\ref{EMSL-DAC-MCF}(f)), the point cloud remains uniformly distributed, allowing the surface to continue evolving. 
As time progresses, the surface clearly converges toward a spherical shape—visually enhanced by scaling the radius.
At $t=0.0896913948$, $E_h=8.18\times10^{-39}$ indicate that the surface has nearly collapsed to a single point, as one expects for mean curvature flow (Figure~\ref{EMSL-DAC-MCF}(g)). 
These results clearly demonstrate the effectiveness of our algorithm while accurately capturing MCF dynamics on point clouds.

\section{Conclusion}\label{Sec-Con}
In this paper, we proposed a novel artificial tangential velocity method for the evolution of surfaces represented by point clouds. 
To address the issue of unexpected point clustering during surface evolution, we introduced a surface density field and derived a Fokker–Planck equation to govern its evolution. 
Then we formulated an evolutionary flow with an artificial tangential velocity coupled with the Fokker–Planck equation to dynamically guide the distribution of points. 
Redistribution and target distribution have been introduced in the case when extra configurations are needed in the algorithms.
Extensive numerical experiments demonstrated the robustness and effectiveness of our method across a variety of scenarios. 
Overall, the proposed method and the idea provides foundational approach for accurate and stable simulation of surface evolution on point clouds. 
Hopefully, this approach is also possible to tackle more complex surface evolution problems. 

It is worth noting that the trigger of the redistribution algorithm is pulled empirically in the current work.
A more intelligent switch would be desirable, as well there are many theoretical questions open for future investigations.

\section*{Acknowledgement}\label{Acknow}
The work of JP and ZS was supported by the National Natural Science Foundation of China (NSFC) No. 92370125.
The work of GD was supported by NSFC No. 12471402. 
The work of HG was partially supported by the Andrew Sisson Fund and the Faculty Science Researcher Development Grant of The University of Melbourne.  

\bibliographystyle{elsarticle-num} 
\bibliography{ref.bib}

\begin{thebibliography}{10}
\expandafter\ifx\csname url\endcsname\relax
  \def\url#1{\texttt{#1}}\fi
\expandafter\ifx\csname urlprefix\endcsname\relax\def\urlprefix{URL }\fi
\expandafter\ifx\csname href\endcsname\relax
  \def\href#1#2{#2} \def\path#1{#1}\fi

\bibitem{bansch2005finite}
E.~B{\"a}nsch, P.~Morin, R.~H. Nochetto, A finite element method for surface diffusion: the parametric case, Journal of Computational Physics 203~(1) (2005) 321--343.

\bibitem{marchandise2011high}
E.~Marchandise, C.~C. de~Wiart, W.~Vos, C.~Geuzaine, J.-F. Remacle, High-quality surface remeshing using harmonic maps—{Part II}: Surfaces with high genus and of large aspect ratio, International Journal for Numerical Methods in Engineering 86~(11) (2011) 1303--1321.

\bibitem{remacle2010high}
J.-F. Remacle, C.~Geuzaine, G.~Compere, E.~Marchandise, High-quality surface remeshing using harmonic maps, International Journal for Numerical Methods in Engineering 83~(4) (2010) 403--425.

\bibitem{dziuk1990algorithm}
G.~Dziuk, An algorithm for evolutionary surfaces, Numerische Mathematik 58~(1) (1990) 603--611.

\bibitem{bonito2010parametric}
A.~Bonito, R.~H. Nochetto, M.~S. Pauletti, Parametric {FEM} for geometric biomembranes, Journal of Computational Physics 229~(9) (2010) 3171--3188.

\bibitem{dziuk2008computational}
G.~Dziuk, Computational parametric {Willmore} flow, Numerische Mathematik 111 (2008) 55--80.

\bibitem{cheung2015localized}
K.~C. Cheung, L.~Ling, S.~J. Ruuth, A localized meshless method for diffusion on folded surfaces, Journal of Computational Physics 297 (2015) 194--206.

\bibitem{deckelnick2018stability}
K.~Deckelnick, V.~Styles, Stability and error analysis for a diffuse interface approach to an advection-diffusion equation on a moving surface, Numerische mathematik 139 (2018) 709--741.

\bibitem{dziuk2007finite}
G.~Dziuk, C.~M. Elliott, Finite elements on evolving surfaces, IMA Journal of numerical analysis 27~(2) (2007) 262--292.

\bibitem{dziuk2012fully}
G.~Dziuk, C.~M. Elliott, A fully discrete evolving surface finite element method, SIAM Journal on Numerical Analysis 50~(5) (2012) 2677--2694.

\bibitem{lehrenfeld2018stabilized}
C.~Lehrenfeld, M.~A. Olshanskii, X.~Xu, A stabilized trace finite element method for partial differential equations on evolving surfaces, SIAM Journal on Numerical Analysis 56~(3) (2018) 1643--1672.

\bibitem{li2018direct}
Y.~Li, X.~Qi, J.~Kim, Direct discretization method for the {Cahn--Hilliard} equation on an evolving surface, Journal of Scientific Computing 77 (2018) 1147--1163.

\bibitem{petras2019least}
A.~Petras, L.~Ling, C.~Piret, S.~J. Ruuth, A least-squares implicit {RBF-FD} closest point method and applications to {PDEs} on moving surfaces, Journal of Computational Physics 381 (2019) 146--161.

\bibitem{barrett2020review}
J.~W. Barrett, H.~Garcke, R.~N{\"u}rnberg, Parametric finite element approximations of curvature-driven interface evolutions, in: Handbook of numerical analysis, Vol.~21, Elsevier, 2020, pp. 275--423.

\bibitem{barrett2007parametric}
J.~W. Barrett, H.~Garcke, R.~N{\"u}rnberg, A parametric finite element method for fourth order geometric evolution equations, Journal of Computational Physics 222~(1) (2007) 441--467.

\bibitem{barrett2008parametric}
J.~W. Barrett, H.~Garcke, R.~N{\"u}rnberg, On the parametric finite element approximation of evolving hypersurfaces in $\mathbb{R}^{3}$, Journal of Computational Physics 227~(9) (2008) 4281--4307.

\bibitem{barrett2008parametric1}
J.~W. Barrett, H.~Garcke, R.~N{\"u}rnberg, Parametric approximation of {Willmore} flow and related geometric evolution equations, SIAM Journal on Scientific Computing 31~(1) (2008) 225--253.

\bibitem{bao2022volume}
W.~Bao, H.~Garcke, R.~N{\"u}rnberg, Q.~Zhao, Volume-preserving parametric finite element methods for axisymmetric geometric evolution equations, Journal of Computational Physics 460 (2022) 111180.

\bibitem{bao2021structure}
W.~Bao, Q.~Zhao, A structure-preserving parametric finite element method for surface diffusion, SIAM Journal on Numerical Analysis 59~(5) (2021) 2775--2799.

\bibitem{zhao2021energy}
Q.~Zhao, W.~Jiang, W.~Bao, An energy-stable parametric finite element method for simulating solid-state dewetting, IMA Journal of Numerical Analysis 41~(3) (2021) 2026--2055.

\bibitem{barrett2013eliminating}
J.~W. Barrett, H.~Garcke, R.~N{\"u}rnberg, Eliminating spurious velocities with a stable approximation of viscous incompressible two-phase {Stokes} flow, Computer Methods in Applied Mechanics and Engineering 267 (2013) 511--530.

\bibitem{barrett2015stable}
J.~W. Barrett, H.~Garcke, R.~N{\"u}rnberg, A stable parametric finite element discretization of two-phase {Navier-Stokes} flow, Journal of Scientific Computing 63 (2015) 78--117.

\bibitem{fu2020arbitrary}
G.~Fu, Arbitrary {Lagrangian-Eulerian} hybridizable discontinuous galerkin methods for incompressible flow with moving boundaries and interfaces, Computer Methods in Applied Mechanics and Engineering 367 (2020) 113158.

\bibitem{ganesan2017ale}
S.~Ganesan, A.~Hahn, K.~Simon, L.~Tobiska, {ALE-FEM} for two-phase and free surface flows with surfactants, Transport Processes at Fluidic Interfaces (2017) 5--31.

\bibitem{m2017approximations}
C.~M.~Elliott, H.~Fritz, On approximations of the curve shortening flow and of the mean curvature flow based on the deturck trick, IMA Journal of Numerical Analysis 37~(2) (2017) 543--603.

\bibitem{hu2022evolving}
J.~Hu, B.~Li, Evolving finite element methods with an artificial tangential velocity for mean curvature flow and {Willmore} flow, Numerische Mathematik 152~(1) (2022) 127--181.

\bibitem{duan2024new}
B.~Duan, B.~Li, New artificial tangential motions for parametric finite element approximation of surface evolution, SIAM Journal on Scientific Computing 46~(1) (2024) A587--A608.

\bibitem{li2021convergence}
B.~Li, Convergence of {Dziuk's} semidiscrete finite element method for mean curvature flow of closed surfaces with high-order finite elements, SIAM Journal on Numerical Analysis 59~(3) (2021) 1592--1617.

\bibitem{bai2024convergent}
G.~Bai, J.~Hu, B.~Li, A convergent evolving finite element method with artificial tangential motion for surface evolution under a prescribed velocity field, SIAM Journal on Numerical Analysis 62~(5) (2024) 2172--2195.

\bibitem{bai2024new}
G.~Bai, B.~Li, A new approach to the analysis of parametric finite element approximations to mean curvature flow, Foundations of Computational Mathematics 24~(5) (2024) 1673--1737.

\bibitem{bai2024convergence}
G.~Bai, B.~Li, Convergence of a stabilized parametric finite element method of the {Barrett-Garcke-N{\"u}rnberg} type for curve shortening flow, Mathematics of Computation (2024).

\bibitem{risken1991fokker}
H.~Risken, T.~Caugheyz, The {F}okker-{P}lanck equation: Methods of solution and application, Journal of Applied Mechanics 58~(3) (1991) 860.

\bibitem{do1976differential}
M.~P. do~Carmo, Differential geometry of curves and surfaces, Englewood Cliffs, New Jersey (1976).

\bibitem{liang2013solving}
J.~Liang, H.~Zhao, Solving partial differential equations on point clouds, SIAM Journal on Scientific Computing 35~(3) (2013) A1461--A1486.

\bibitem{rosser1967runge}
J.~B. Rosser, A {Runge-Kutta} for all seasons, SIAM Review 9~(3) (1967) 417--452.

\bibitem{lai2013local}
R.~Lai, J.~Liang, H.~Zhao, A local mesh method for solving {PDEs} on point clouds., Inverse Problems \& Imaging 7~(3) (2013).

\bibitem{elliott2012ale}
C.~M. Elliott, V.~Styles, An {ALE ESFEM} for solving {PDEs} on evolving surfaces, Milan Journal of Mathematics 80~(2) (2012) 469--501.

\bibitem{mantegazza2011lecture}
C.~Mantegazza, Lecture Notes on {Mean Curvature Flow}, Vol. 290, Springer Science \& Business Media, 2011.

\bibitem{duan2021high}
B.~Duan, B.~Li, Z.~Zhang, High-order fully discrete energy diminishing evolving surface finite element methods for a class of geometric curvature flows, Annals of Applied Mathematics 37~(4) (2021) 405--436.

\end{thebibliography}

\end{document}